\documentclass[aop,MSNbibl,seceqn,dvips]{arximspdf}
\usepackage{tipa}
\doi{10.1214/13-AOP849} %kopijuoti is PTS
\volume{42}
\issue{4}
\pubyear{2014}
\firstpage{1724}
\lastpage{1768}

\makeatletter
\protected\def\tp{\textpolhook}

\def\cal{\mathcal}
\def\inf{\mathop{\operatorname{inf}}}
\def\esssup{\mathop{\operatorname{ess\,sup}}}
\def\essinf{\mathop{\operatorname{ess\,inf}}}
\def\sup{\mathop{\operatorname{sup}}}
\def\inf{\mathop{\operatorname{inf}}}

\newtheorem{lemma}{Lemma}[section]
\newproclaim{remark}{Remark}[section]
\newproclaim{example}{Example}[section]
\newtheorem{theorem}{Theorem}[section]
\newtheorem{corollary}{Corollary}[section]

\newproclaim{definition}{Definition}[section]
\newtheorem{proposition}{Proposition}[section]
\newtheorem{condition}{Condition}[section]
\makeatother

\begin{document}
\begin{frontmatter}

\title{Value in mixed strategies for zero-sum stochastic differential
games without Isaacs condition}
\runtitle{Value in mixed strategies for zero-sum SDGs}

\begin{aug}
\author[a]{\fnms{Rainer} \snm{Buckdahn}\ead[label=e1]{rainer.buckdahn@univ-brest.fr}\thanksref{t1}},
\author[b]{\fnms{Juan} \snm{Li}\corref{}\ead[label=e2]{juanli@sdu.edu.cn}\thanksref{t2}}
\and
\author[a]{\fnms{Marc}~\snm{Quincampoix}\ead[label=e3]{marc.quincampoix@univ-brest.fr}\thanksref{t1}}

\thankstext{t1}{Supported in part by the Commission of
the European Communities under the 7th Framework Programme Marie
Curie Initial Training Networks
Project ``Deterministic and Stochastic Controlled Systems and
Applications'' FP7-PEOPLE-2007-1-1-ITN, No. 213841-2 and project
SADCO, FP7-PEOPLE-2010-ITN, No. 264735 and by the French National Research Agency ANR-10-BLAN
0112.}
\thankstext{t2}{Supported by the NSF of P.R. China (Nos.
11071144, 11171187, 11222110), Shandong Province (Nos. BS2011SF010,
JQ201202), Program for New Century Excellent Talents in University (No.~NCET-12-0331), 111 Project (No. B12023).}
%%\and
%%\author[B]{\fnms{} \snm{}}
\runauthor{R. Buckdahn, J. Li and M. Quincampoix}
\affiliation{Universit\'{e} de Bretagne Occidentale and Shandong University,
Shandong University, Weihai, and Universit\'{e}
de Bretagne Occidentale}
\address[a]{R. Buckdahn\\
M. Quincampoix\\
Laboratoire de Math\'{e}matiques\\
CNRS-UMR 6205\\
Universit\'{e} de Bretagne Occidentale\\
6, avenue Victor-le-Gorgeu\\
B.P. 809, 29285 Brest cedex\\
France\\
%and\\
%School of Mathematics and Statistics\\
%Shandong University, Weihai\\
%NO 180 Wenhua Xilu\\
%Weihai, Shandong Province, 264209\\
%P. R. China\\
\printead{e1}\\
\phantom{E-mail:\ }\printead*{e3}}
\address[b]{J. Li\\
School of Mathematics and Statistics\\
Shandong University, Weihai\\
NO 180 Wenhua Xilu\\
Weihai, Shandong Province, 264209\\
P. R. China\\
\printead{e2}}

%Laboratoire de Math\'{e}matiques\\
%CNRS-UMR 6205\\
%Universit\'{e} de Bretagne Occidentale\\
%6, avenue Victor-le-Gorgeu\\
%B.P. 809, 29285 Brest cedex\\
%France\\
\end{aug}

% HISTORY:
\received{\smonth{6} \syear{2012}}
\revised{\smonth{3} \syear{2013}}

% ABSTRACT
%
\begin{abstract}
In the present work, we consider
2-person zero-sum stochastic differential games with a nonlinear
pay-off functional which is defined through a backward stochastic
differential equation. Our main objective is to study for such a
game the problem of the existence of a value without Isaacs
condition. Not surprising, this requires a suitable concept of mixed
strategies which, to the authors' best knowledge, was not known in the
context of stochastic differential games. For this, we consider
nonanticipative strategies with
a delay defined through a partition $\pi$ of the time interval $[0,T]$.
The underlying stochastic controls for the both players are randomized
along $\pi$ by a hazard which is independent of the governing
Brownian motion, and knowing the information available at the
left time point $t_{j-1}$ of the subintervals generated by~$\pi$,
the controls of Players~1 and~2 are conditionally independent
over $[t_{j-1},t_j)$. It is shown that the associated
lower and upper value functions $W^{\pi}$ and $U^{\pi}$ converge
uniformly on compacts to a function $V$, the so-called value in
mixed strategies, as the mesh of $\pi$ tends to zero. This function
$V$ is characterized as the unique viscosity solution of the
associated Hamilton--Jacobi--Bellman--Isaacs equation.
\end{abstract}

% KEYWORDS


\end{frontmatter}
\newpage

%s1 #&#
\section{Introduction}\label{sec1}

In our work, we investigate 2-person zero-sum stochastic differential
games which dynamics\ are defined by a doubly controlled stochastic
differential equation (SDE)
%
%
%e1.1 #&#
%
\begin{eqnarray}
\label{0a} %
dX_{s}^{t,x;u,v}&=&b
\bigl(s,X_{s}^{t,x;u,v},u_{s},v_{s}
\bigr)\,ds\nonumber\\
&&{}+ \sigma \bigl(s,X_{s}^{t,x;u,v},u_{s},v_{s}
\bigr)\,dB_s,\qquad s\in[t,T],
\\
X_{t}^{t,x;u,v}&=&x\in R^d,\nonumber
\end{eqnarray}
driven by a Brownian motion $B$, and endowed with pay-off functionals
defined through a doubly controlled backward stochastic differential
equation (BSDE) (see Section~\ref{sec2} for details) which, in the classical
case, reduces
to
%
%
%e1.2 #&#
%
\begin{equation}
\label{1a} I(t,x;u,v)=E \biggl[\Phi \bigl(X^{t,x;u,
v}_T \bigr)+
\int_t^Tf \bigl(s,X^{t,x;u,v}_s,u_s,v_s
\bigr)\,ds \biggr]
\end{equation}
[see (\ref{DPP-aaa})]. The initial data $(t,x)$ of the game belong to
$[0,T]\times R^d$, and the control processes $u=(u_s)$ and $v=(v_s)$
used by Players~1 and~2, take their values in compact
metric spaces $U$ and $V$, respectively. While the objective of Player
1 is
to maximize the pay-off $I(t,x;u,v)$, that of Player 2 is to minimize it:
Indeed, for Player 2 $I(t,x;u,v)$ represents a cost functional.
However, apart
from rather strong assumptions on the coefficients, for example, that
of independence of the controls $(u, v)$ and of strict
ellipticity for the diffusion coefficient $\sigma\sigma^*(t,x)\ge
\alpha
\cdot I_{R^d}, (t,x)\in[0,T]\times R^d$, for some $\alpha>0$ (refer to
Hamadene, Lepeltier, and Peng~\cite{HLP}), if one wants
to have a dynamic programming principle (DPP) the players can, in
general, not play a game of
the type ``control against control''; they can play, for instance, games
of the type ``nonanticipative strategy against control'' (see, e.g.,
\cite{FS,Buckdahn-Li-2008}) or games of the type ``NAD-strategy
against NAD-strategy'', where NAD stands for nonanticipativity
with delay (see, e.g., \cite{BCR04} and \cite{BCQ11}).

However, a central question in the theory of 2-person zero-sum stochastic
differential games is that of sufficient conditions, under which the
game admits a value, that is, under which the lower and the upper value
functions of the stochastic differential game coincide. In the literature,
since the famous works by Isaacs \cite{I} for the case of deterministic
differential
games and that by Fleming and Souganidis \cite{FS} for stochastic
differential games (see also \cite{FH}), various authors have shown
the equality between the lower
and the upper value functions under the so-called Isaacs condition.

Let us be more precise: Generalizing the pioneering paper on
stochastic differential games by Fleming and Souganidis \cite{FS},
Buckdahn and Li \cite{Buckdahn-Li-2008}, and also Buckdahn, Cardaliaguet
and Quincampoix \cite{BCQ11}, associated the dynamics (\ref{0a}) with
nonlinear cost functionals defined through a BSDE, which was first
introduced by Pardoux and Peng~\cite{PaPe1}:
%
%
%e1.3 #&#
%
\begin{equation}
\label{3a} \hspace*{14pt}\cases{ %
-dY^{t,x; u, v}_s = f
\bigl(s,X^{t,x;u,v}_s,Y^{t,x;u,v}_s,Z^{t,x;u,v}_s,
u_s,v_s \bigr)\,ds -Z^{t,x;u,v}_s
\,dB_s,& \vspace*{2pt}
\cr
Y^{t,x; u, v}_T = \Phi
\bigl(X^{t,x; u, v}_T \bigr),\qquad s\in[t,T].}\hspace*{-20pt}
\end{equation}
They considered as pay-off functional the random variable
(measurable with respect to the information available before the
beginning of the game)
%
%
%e1.4 #&#
%
\begin{equation}
J(t,x;u,v)=Y_t^{t,x; u,v},
\end{equation}
and the lower and the upper value functions for the game
over the time interval $[t,T]$ were introduced, respectively, by putting
%
%
%e1.5 #&#
%
\begin{eqnarray}
W(t,x)&:=& \esssup_{\alpha} \essinf_{\beta} J(t,x;
\alpha,\beta),
\nonumber
\\[-8pt]
\\[-8pt]
\nonumber
U(t,x)&:=& \essinf_{\beta}\esssup_{\alpha}J(t,x; \alpha,\beta)
\qquad (t,x)\in[0,T]\times R^d,
\end{eqnarray}
where $\alpha$ runs the NAD-strategies for Player 1 and $\beta$ those
for Player 2. Given such a couple of admissible NAD-strategies, the
cost functional $J(t,x;\alpha,\beta)$ is defined through the unique
couple of admissible controls $(u,v)$ satisfying $\alpha(v)=u,
\beta(u)=v$, by putting $J(t,x;\alpha,\beta)=J(t,x;u,v)$ (e.g., refer
to~\cite{BCQ11}). We emphasize
that in the above definition the classical case, where $f(s,x,y,z,u,v)=
f(s,x,u,v)$ is independent of $(y,z)$, can be obtained by replacing
$J(t,x;\alpha,\beta)$ by $E[J(t,x;\alpha,\beta)]=I(t,x;\alpha
,\beta)$
[see (\ref{1a})] and the essential supremum and the essential infimum
over a family of
random variables by the supremum and the infimum, respectively; this
does not change the upper and the lower value functions (see
Remark 3.4, \cite{Buckdahn-Li-2008}). The authors showed that, for the
Hamiltonians
%
%
%e1.6 #&#
%
\begin{eqnarray}
H(t,x,y,p,A,u,v)&=&\frac12 \operatorname{tr} \bigl(\sigma
\sigma^*(t,x,u,v)A \bigr) +b(t,x,u,v)p
\nonumber\\
& &{}+f \bigl(t,x,y,p\sigma(t,x,u,v),u,v \bigr),
\nonumber
\\[-8pt]
\\[-8pt]
\nonumber
H^-(t,x,y,p,A)&=&\sup_{u\in U}\inf
_{v\in V}H(t,x,y,p,A,u,v),
\\
H^+(t,x,y,p,A)&=&\inf_{v\in V}\sup_{u\in U}H(t,x,y,p,A,u,v),\nonumber
\end{eqnarray}
$(t,x,y,p,A)\in[0,T]\times R^d\times R\times R^d\times S^d$ ($S^d$ denotes the space of symmetric real matrices of the size
$d\times d$), $W$ and $U$ are the unique viscosity solutions of the
following Hamilton--Jacobi--Bellman--Isaacs (HJBI) equations in the class of
continuous functions with polynomial growth, respectively:
%
%
%e1.7 #&#
%
\begin{eqnarray}
\qquad\frac{\partial}{\partial t}W(t,x)+H^- \bigl(t,x, \bigl(W, \nabla
W,D^2W \bigr) (t,x) \bigr)&=&0,\qquad W(T,x)=\Phi(x),
\nonumber
\\[-8pt]
\\[-8pt]
\nonumber
\frac{\partial}{\partial t}U(t,x)+H^+ \bigl(t,x, \bigl(U,\nabla
U,D^2U \bigr) (t,x) \bigr)&=&0,\qquad U(T,x)=\Phi(x).
\end{eqnarray}
Isaacs condition
says that
%
%
%e1.8 #&#
%
\begin{eqnarray}
H^-(t,x,y,p,A)=H^+(t,x,y,p,A)
\nonumber
\\[-8pt]
\\[-8pt]
\eqntext{(t,x,y,p,A)\in[0,T] \times R^d\times
R\times R^d\times S^d,}
\end{eqnarray}
and under it the both above PDEs\vadjust{\goodbreak} coincide and the uniqueness
of the solution implies that $W(t,x)=U(t,x), (t,x)\in[0,T]\times R^d$,
that is, the game has a value.

But how to get a value, when Isaacs condition is not assumed?
Recently, in \cite{Buckdahn-Li_Quincampoix-2011} the authors
studied deterministic differential games without assuming Isaacs
condition. They considered an adequate notion of mixed strategies
related with a suitable randomization, and were thus able to prove
that such defined upper and lower value functions coincide, and that
this value function defined through mixed strategies
satisfies a Hamilton--Jacobi--Isaacs equation. We also refer to the works
of Chentsov, Krasovskii and Subbotin for the existence of the value
of deterministic differential games \cite{KRSU,SUC}: They studied the
problems of deterministic
differential games without Isaacs condition through positional
strategies but with techniques which differ from those in
\cite{Buckdahn-Li_Quincampoix-2011}. To the authors' best knowledge,
there does not exist any work on the existence of the value
of stochastic differential games without assuming Isaacs condition, it
has been an open problem until now. However, there are also different
recent works studying stochastic differential games without Isaacs'
condition, but without the objective to show the existence of a value
of the game. For instance, Krylov \cite{K1,K2} studied regularity
properties and the dynamic programming principle for the upper value
function of a stochastic differential game over a domain, by starting
from the Isaacs equation; for this he used the idea of \'{S}wi\c{e}ch
\cite{SA} that the viscosity solutions of nondegenerate Isaacs
equations have some regularity properties which can be used for the approach.

In the present work, our objective is to solve this open problem, that
is, to extend the results of
\cite{Buckdahn-Li_Quincampoix-2011} from deterministic
differential games without Isaacs condition to stochastic
differential games. Since this work was heavily inspired by \cite
{Buckdahn-Li_Quincampoix-2011}, we consider the game of the type
``NAD-straegies against
NAD-strategies''. The delay of the nonanticipative strategies
is defined through a partition $\pi=\{0=t_0<t<t_1<\cdots<
t_n=T\}$ of the time interval $[0,T]$. The underlying stochastic
controls for the both players are randomized along the
partition $\pi$ by a hazard which is independent of the
governing Brownian motion, and knowing all information
available at the left time point $t_{j-1}$ of the subintervals
generated by $\pi$, the controls of Players~1 and~2 are
conditionally independent over $[t_{j-1},t_j)$.

While the dynamics are defined by (\ref{0a}), the BSDE defining the pay-off
functional has to take into account that, first, the controls
of the both players are randomized by a hazard independent of the
governing Brownian motion, and second, the both players make
the randomization of their controls conditionally independent of
each other and reveal the information related with only
at the end of each subinterval generated by the partition $\pi$.
This has as consequence that the BSDE has to be considered under
a filtration $\widetilde{\mathbb{F}}^\pi$ which is smaller than the
filtration $\mathbb{F}^\pi$ (but larger than
the Brownian one) for the dynamics
(\ref{0a}); see BSDE (\ref{BSDE}).

With the help of the cost functional defined through our BSDE
we introduce the lower and the upper value functions along a partition
$\pi$, $W^\pi$ and
$U^\pi$. For these, a priori, random fields we prove that they are
deterministic and satisfy along the partition $\pi$, at its points,
the dynamic programming principle. This dynamic programming principle
combined with Peng's BSDE method, refer to Peng~\cite{Pe1}, which we
have to redevelop for our settings here is crucial for the proof that
$W^\pi$ and $U^\pi$ converge uniformly on compacts, as the mesh of
$\pi$ tends to zero, and their limit $V$, the so-called value in
mixed strategies can be characterized as the unique viscosity solution
of the Hamilton--Jacobi--Bellman--Isaacs equation
%
%
%e1.9 #&#
%
\begin{eqnarray}
\label{PDEa} %
\hspace*{20pt}\frac{\partial}{\partial t}V(t,x)+\sup_{\mu\in{\cal P}(U)}
\inf_{\nu\in{\cal P}(V)}H \bigl(t,x, \bigl(V,\nabla V,D^2V \bigr)
(t,x),\mu,\nu \bigr) &=&0,
\nonumber
\\[-8pt]
\\[-8pt]
\nonumber
V(T,x)&=&\Phi(x),
\end{eqnarray}
where
%
%
%e1.10 #&#
%
\begin{eqnarray}
& &H(t,x,y,p,A,\mu,\nu)
\nonumber
\\
&&\qquad=\int_{U\times V} \biggl(\frac12 \operatorname{tr} \bigl(
\sigma \sigma^*(t,x,u,v)A \bigr)+b(t,x,u,v)p
\\
& &\hspace*{58pt}\qquad\quad{}+f \bigl(t,x,y,p\sigma(t,x,u,v),u,v \bigr) \biggr)\mu\otimes
\nu(du\,dv),
\nonumber
\end{eqnarray}
$(t,x,y,p,A)\in[0,T]\times R^d\times R\times R^d\times S^d$. Here
${\cal P}(U)$ denotes the space of all probability measures
on $U$, ${\cal P}(V)$ all on $V$. Since both control state spaces
$U$ and $V$ are supposed to be compact and metric, ${\cal P}(U)$
and ${\cal P}(V)$ are convex and compact, and from the bi-linearity
of $H(t,x,y,p,A,\mu,\nu)$ in $(\mu,\nu)$ we have that for PDE
(\ref{PDEa}) the following Isaacs condition is automatically satisfied:
%
%
%e1.11 #&#
%
\begin{eqnarray}
&&\sup_{\mu\in{\cal P}(U)}\inf_{\nu\in{\cal P}(V)}
H(t,x,y,p,A,\mu,\nu)
\nonumber
\\[-8pt]
\\[-8pt]
\nonumber
&&\qquad=\inf_{\nu\in{\cal P}(V)}\sup_{\mu\in{\cal P}(U)}
H(t,x,y,p,A,\mu, \nu).
\end{eqnarray}
Of course, PDE (\ref{PDEa}) could have been also derived by considering
weak controls, that is, controls with values in ${\cal P}(U)$ and
${\cal P}(U)$,
but our objective has been to work with controls taking values in $U$
and $V$,
respectively, even for the price of a randomization.

Let us point out that the fact that, in our approach, the dynamics and
the BSDE have to be studied under different filtration, means that
unlike in \cite{Buckdahn-Li-2008} and \cite{BCQ11} we are not anymore in
a Markovian framework here for our BSDE. This requires new
approaches, not only for the redevelopment of Peng's BSDE method~\cite{Pe1}
in our settings (Section~\ref{sec4}), but also for the proof that the upper
and the lower value functions are deterministic and H\"{o}lder
continuous with respect to the time parameter.

Let us explain the organization of the paper. In Section~\ref{sec2}, we
introduce the settings for our stochastic differential games,
we define for both players the space of admissible controls
along a partition $\pi$ as well as the notion of NAD-strategies
with respect to $\pi$. Moreover, we introduce the dynamics, the
pay-off functional defined through a BSDE, as well as the upper and
the lower value functions $W^{\pi}$ and $U^{\pi}$ along $\pi$.
In Section~\ref{sec3}, we study properties of $W^{\pi}$ and $U^{\pi}$. We
show, in particular, that they are deterministic continuous
functions which, with respect to the points of the partition
$\pi$, satisfy the dynamic programming principle. In Section~\ref{sec4}, finally, it is shown that, as the mesh of $\pi$ tends to zero,
$W^{\pi}$ and $U^{\pi}$ converge uniformly on compacts to the unique
viscosity solution of the associated Hamilton--Jacobi--Bellman--Isaacs
equation. For this, Peng's BSDE method is redeveloped for our
settings.

%s2 #&#
\section{Preliminaries. Settings of the stochastic differential games}\label{sec2}

Let us begin with introducing the probability space
\[
(\Omega_1,{\cal F}_1,P_1):= \bigl(
\bigl(R^2 \bigr)^{\mathbb{N}},{\cal B} \bigl(R^2
\bigr)^{\otimes\mathbb{N}}, Q_2^{\otimes\mathbb{N}} \bigr),
\]
where $Q_2$ denotes the two-dimensional standard Normal
distribution on the real plane $R^2$ endowed with its Borel
$\sigma$-field ${\cal B}(R^2)$, and $\mathbb{N}$ is the set of all
positive integers. Then, by the above definition,
$\Omega_1=(R^2)^{\mathbb{N}}$ is the space of all $R^2$-valued
sequences $\rho=(\rho_j=(\rho_{j,1},\rho_{j,2}))_{j\ge1}$, and
${\cal F}_1={\cal B} (R^2)^{\otimes\mathbb{N}}$ is the product Borel
$\sigma$-field taken over the sequence of $\sigma$-fields, which all
elements coincide with ${\cal B} (R^2)$, and
$P_1=Q_2^{\otimes\mathbb{N}}$ is the product measure over
$(\Omega_1,{\cal F}_1)$. Let us denote the coordinate mappings on
$\Omega_1$ by $\zeta_j=(\zeta_{j,1},\zeta_{j,2})\dvtx \Omega
_1\rightarrow
R^2$, $j\ge1$:
\[
\zeta_j(\rho)= \bigl(\zeta_{j,1}(\rho),
\zeta_{j,2}(\rho) \bigr)= (\rho_{j,1},\rho_{j,2}),
\qquad \rho= \bigl((\rho_{j,1},\rho_{j,2})
\bigr)_{j\ge1} \in\Omega_1.
\]
We observe that ${\cal F}_1$ coincides with the smallest
$\sigma$-field on $\Omega_1$, with respect to which all coordinate
mappings $\zeta_j,j\ge1$, are measurable.

However, for the study of our stochastic differential games we also
need the classical Wiener space $(\Omega_2,{\cal F}_2,P_2)$, where
$\Omega_2$ is the set of all continuous functions from $[0,T]$ with
values in $R^d$ and starting from zero, endowed with the supremum
norm [i.e., $\Omega_2=C_0( [0,T];R^d)$], and ${\cal F}_2$ is the Borel
$\sigma$-field on $\Omega_2$ completed with respect to the Wiener
measure $P_2$ under which the coordinate process
$B_t(\omega')=\omega'(t), t\in[0,T], \omega'\in\Omega_2$, is a
Brownian motion.

Let us denote by $(\Omega,{\cal F},P)$ the product probability space
\[
(\Omega,{\cal F},P)=(\Omega_1,{\cal F}_1,P_1)
\otimes (\Omega_2, {\cal F}_2,P_2),
\]
which we complete with respect to the probability measure
$P$, and let us extend the coordinate mappings $\zeta$ and $B$ in a
canonical way from $\Omega_1$ and $\Omega_2$, respectively, to~$\Omega$:
\begin{eqnarray}
\zeta_j(\omega):=\zeta_j(\rho),\qquad B_t(
\omega):=B_t\bigl(\omega'\bigr)=\omega'(t),
\nonumber\\
\eqntext{\omega=\bigl(\rho,\omega'\bigr)\in\Omega=\Omega_1
\times\Omega_2, j\ge1,  t\in[0,T].}
\end{eqnarray}
Let us now introduce the filtration with which we work on our
probability space $(\Omega,{\cal F},P)$. By $\mathbb{F}^B=({\cal
F}^B_t)_{t\in[0,T]}$ we denote the filtration generated by the
Brownian motion $B$ and completed by all $P$-null sets. In addition
to the filtration $\mathbb{F}^B$, we also need larger ones, defined
along a partition $\pi=\{0=t_0<t_1<\cdots<t_n=T\}$ of the interval
$[0,T]$. Given such a partition $\pi$, we define $\mathbb{F}^{\pi,i}=
({\cal F}^{\pi,i}_t)_{t\in[0,T]}$, with
\[
{\cal F}^{\pi,i}_t={\cal F}^B_t\vee
\sigma \bigl\{\zeta_\ell= (\zeta_{\ell,1}, \zeta_{\ell,2})
(1\le\ell\le j-1), \zeta_{j,i} \bigr\},
\]
$t\in[t_{j-1},t_j)$, $1\le j\le n, i=1,2$, and we put ${\cal F}^{\pi,i}_T={\cal F}^{\pi,i}_{T-},
i=1,2$. Notice that, for $j=1$, that is, on the time interval
$[t_0,t_1)$, by convention, ${\cal F}_t^{\pi,i}={\cal F}_t^B\vee
\sigma\{
\zeta_{1,i}\}, i=1, 2$. We shall also introduce the filtration
$\mathbb{F}^{\pi}=\mathbb{F}^{\pi,1}\vee\mathbb{F}^{\pi
,2}=({\cal
F}^{\pi}_t={\cal F}^{\pi,1}_t\vee{\cal F}^{\pi,2}_t)_{t\in[0,T]}$,
and we remark that, for $t\in[t_{j-1},t_j)$,
\[
{\cal F}^{\pi}_t={\cal F}^B_t\vee{
\cal H}_j\qquad \mbox{where } {\cal H}_j: = \sigma
\bigl\{\zeta_\ell= (\zeta_{\ell,1}, \zeta_{\ell,2}) (1\le
\ell\le j) \bigr\}.
\]

Finally, we will also need a smaller filtration,
$\widetilde{\mathbb{F}}^\pi=(\widetilde{\cal
F}^{\pi}_t)_{t\in[0,T]}$ with $\widetilde{\cal F}^{\pi}_t:={\cal
F}^B_t\vee{\cal H}_{j-1}$, for $t\in[t_{j-1},t_j), 1\le j\le n$.
Observe that, for all $t\in[t_{j-1},t_j)$, knowing $\widetilde{\cal
F}^\pi_t={\cal F}_t^B\vee{\cal H}_{j-1}$, the $\sigma$-fields ${\cal
F}^{\pi,1}_t$ and ${\cal F}^{\pi,2}_t$ are conditionally
independent.

Let us consider two compact metric spaces $U$ and $V$ as control
state spaces used by the Players 1 and 2, respectively. By ${\cal
P}(U)$ and ${\cal P}(V)$, we denote the space of all probability
measures over $U$ and $V$, endowed with its Borel
$\sigma$-field ${\mathcal{B}}(U)$ and ${\mathcal{B}}(V)$, respectively.
We also observe that it is an immediate consequence
of Skorohod's Representation theorem that the set ${\cal P}(U)$
[resp., ${\cal P}(V)$] coincides with the set of the laws of all
$U$-valued (resp., $V$-valued) random variables defined over
$([0,1], {\cal B}([0,1]),\lambda_1)$ [$\lambda_1$ denotes the
Lebesgue measure on $([0,1], {\cal B}([0,1]))$]. But this latter set
coincides with that of the laws of all random variables defined over
$(R, {\cal B}(R),Q_1)$, where $Q_1$ denotes the standard Normal
distribution over $(R, {\cal B}(R))$. Indeed, denoting by
\[
\Phi_{0,1}(x)=\frac{1}{\sqrt{2\pi}}\int_{-\infty}^x
\exp \biggl\{- \frac{y^2}{2} \biggr\}\,dy,\qquad x\in R,
\]
we have that, for any random variable $\xi$ over $([0,1],
{\cal B}([0,1]),\lambda_1)$, the law of $\xi$ with respect to
$\lambda_1$ coincides with that of $\xi(\Phi_{0,1}(\cdot))\dvtx R\rightarrow
R$ under
$Q_1$. A~consequence is that
\[
{\cal P}(U)=\bigl\{P_\xi\dvtx \xi \mbox{ is } U\mbox{-valued random
variable over } \bigl(\Omega,\sigma\{\zeta_{j,1}\}, P\bigr) \bigr\}
\]
and
\[
{\cal P}(V)= \bigl\{P_\xi\dvtx \xi \mbox{ is } V\mbox{-valued random
variable over } \bigl(\Omega,\sigma\{\zeta_{j,2}\}, P\bigr) \bigr\}
\]
for all $j\ge1.$

Let us now introduce the admissible controls for both players along a
given partition $\pi=\{0=t_0<t_1<\cdots<t_n=T\}$ of the time interval
$[0,T]$.
%
%de2.1 #&#
\begin{definition}[(Admissible controls)]\label{adm.control} Given a
partition $\pi$ of the time interval $[0,T]$ and an initial time
$t\in[0,T]$, the space of admissible controls along the partition
$\pi$ for Player 1 for a game over the time interval $[t,T]$ is the
totality of all $U$-valued $\mathbb{F}^{\pi,1}$-predictable
processes $u=(u_s)_{s\in[t,T]}$ defined over the probability space
$(\Omega, {\cal F},P)$; it is denoted by ${\cal U}_{t,T}^\pi$. For
Player 2 the space of admissible controls along the partition $\pi$
${\cal V}_{t,T}^\pi$ is defined similarly: It is the
collection of all $V$-valued $\mathbb{F}^{\pi,2}$-predictable
processes $v=(v_s)_{s\in[t,T]}$ defined over $(\Omega,{\cal F},P)$.
\end{definition}
After having introduced the spaces of admissible controls, we describe
now the dynamics of our stochastic differential games. For this, we
consider the coefficients
\[
b\dvtx [0,T]\times R^d\times U\times V\rightarrow R^d
\quad\mbox{and}\quad\sigma\dvtx [0,T]\times R^d\times U\times V
\rightarrow R^{d\times d}
\]
which we suppose throughout our work to be bounded,
jointly continuous and Lipschitz in $x\in R^d$, uniformly with
respect to $(t,u,v)\in[0,T]\times U\times V$. Let $\pi$ be a
partition of the time interval $[0,T]$. Then, given arbitrary
initial data $t\in[0,T]$ and $\vartheta\in L^2(\Omega,{\cal F}_t^\pi,
P;R^d)$ as well as admissible control processes $u\in{\cal U}_{t,
T}^\pi$ and $v\in{\cal V}_{t,T}^\pi$, we consider the SDE
%
%
%e2.1 #&#
%
\begin{eqnarray}
\label{SDE} %
dX_{s}^{t,\vartheta;u,v}&=&b
\bigl(s,X_{s}^{t,\vartheta;u,v},u_{s},v_{s}
\bigr)\,ds+ \sigma \bigl(s,X_{s}^{t,\vartheta;u,v},u_{s},v_{s}
\bigr)\,dB_s
\nonumber
\\[-8pt]
\\[-8pt]
 \eqntext{s\in[t,T], X_{t}^{t,\vartheta;u,v}=\vartheta.}
\end{eqnarray}

Under our assumptions on the coefficients $b$ and $\sigma$,
this SDE has a unique strong solution $X^{t,\vartheta;u,v}=(X_{s}^{t,
\vartheta;u,v})_{s\in[t,T]}$ in the space of $R^d$-valued,
$\mathbb{F}^\pi$-adapted continuous processes. Moreover, we have the
following estimates which are by now standard.

For all $p\ge2$, there exists some constant $C_p\in R$ (only
depending on p, on the Lipschitz constants and the bounds of $b$ and
$\sigma$) such that, for all partitions $\pi$ of $[0,T]$, for all
$t\in[0,T], \vartheta,\vartheta'\in L^2(\Omega,{\cal F}_t^\pi,P;
R^d)$ and
all $u\in{\cal U}^{\pi}_{t,T},v\in{\cal V}^{\pi}_{t,T}$, it holds,
$P$-a.s.,
%
%
%e2.2 #&#
%
\begin{eqnarray}
E \Bigl[\sup_{s\in[t,T]}\bigl|X^{t,\vartheta; u, v}_s
-X^{t,\vartheta';u,v}_s\bigr|^p| {{\mathcal{F}}_t^\pi}
\Bigr] & \leq& C_p\bigl|\vartheta-\vartheta'\bigr|^p,
\nonumber
\\[-8pt]
\\[-8pt]
\nonumber
E \Bigl[ \sup_{s\in[t,T]} \bigl|X^{t,\vartheta;u,v}_s\bigr|^p|{{
\mathcal{F}}_t^\pi } \Bigr] & \leq& C_p
\bigl(1+|\vartheta|^p \bigr).
\end{eqnarray}

Let us now come to the pay-off functional which we associate with
the above dynamics of our game. The pay-off functional is a nonlinear,
recursive one, that is, we define it through a backward
stochastic differential equation. For this, we consider the
terminal pay-off function $\Phi\dvtx R^d\rightarrow R$ which we suppose
to be bounded and Lipschitz, as well as the running pay-off function
$f\dvtx [0,T]\times R^d\times R\times R^d\times U\times V\rightarrow R$
which we assume to be jointly continuous and such that
\begin{longlist}[(iii)]
\item[(i)] $f(t,x,y,z,u,v)$ is Lipschitz in $(x,y,z)\in R^d\times
R\times R^d$, uniformly in $(s,u,v)\in[0,T]\times U\times V$;
\item[(ii)] $f(t,x,y,z,u,v)$ is uniformly continuous on $[0,T]\times R^d
\times R\times\overline{B}_K(0)\times U\times V$, for all $K>0$,
where $\overline{B}_K(0)$ denotes the closed ball in $R^d$ centered at $0$
with diameter $K$;
\item[(iii)] $(t,x,y,u,v)\rightarrow f(t,x,y,0,u,v)$ is bounded.
\end{longlist}

Given a partition $\pi$ of the interval $[0,T]$, initial data
$t\in[0,T],  \vartheta\in L^2(\Omega,{\cal F}^\pi_t,\break  P;R^d)$ and
admissible controls $u\in{\cal U}^{\pi}_{t,T}, v\in{\cal
V}^{\pi}_{t,T}$, we consider the following BSDE governed by the
solution $X^{t,\vartheta;u,v}$ of SDE (\ref{SDE}):
%
%
%e2.3 #&#
%
\begin{equation}
\label{BSDE} \cases{dY^{t,\vartheta; u, v}_s = -E \bigl[f
\bigl(s,X^{t,\vartheta; u, v}_s, Y^{t,
\vartheta; u, v}_s,
Z^{t,\vartheta; u, v}_s, u_s, v_s \bigr) |
\widetilde{\cal F}^\pi_s \bigr] \,ds \vspace*{2pt}\cr
\hspace*{54pt}{}+Z^{t,\vartheta;u, v}_s
\,dB_s+dM^{t,\vartheta; u, v}_s, \vspace*{2pt}
\cr
Y^{t,\vartheta; u, v}_T = E \bigl[\Phi \bigl(X^{t,\vartheta; u, v}_T
\bigr) | \widetilde{\cal F}^\pi_{T} \bigr],}
\end{equation}
where $(E[\gamma_s  | \widetilde{\cal
F}^\pi_s])_{s\in[0,T]}$ is understood as
$\widetilde{\mathbb{F}}^\pi$-optional projection of integrable,
measurable processes $\gamma=(\gamma_s)_{s\in[0,T]}$.

We say that $(Y^{t,\vartheta; u, v},Z^{t,\vartheta; u, v},
M^{t,\vartheta; u, v})$ is a solution of this BSDE, if
\begin{longlist}[(iii)]
\item[(i)] $Y^{t,\vartheta; u, v}\in{\cal
S}^2_{\widetilde{\mathbb{F}}^\pi}(t,T;R)$, that is, $Y^{t,\vartheta
; u,
v}=(Y^{t,\vartheta; u, v}_s)_{s\in[t,T]}$ is an
${\widetilde{\mathbb{F}}^\pi}$-adapted c\`{a}dl\`{a}g process which
is square integrable: $E [\sup_{s\in[t,T]}|Y^{t,\vartheta; u,
v}_s|^2 ]<+\infty$;

\item[(ii)] $Z^{t,\vartheta; u, v}\in L_{\widetilde{\mathbb{F}}^\pi}^2(t,T;
R^d)$, that is, $Z^{t,\vartheta; u, v}=(Z^{t,\vartheta; u,
v}_s)_{s\in[t,T]}$
is an $R^d$-valued, ${\widetilde{\mathbb{F}}^\pi}$-predictable
process such that $E [\int_t^T|Z^{t,\vartheta; u,
v}_s|^2\,ds ]<+\infty$;

\item[(iii)] $M^{t,\vartheta; u, v}\in{\cal M}_{\widetilde{\mathbb{F}}^\pi}^2
(t,T;R)$, that is, $M^{t,\vartheta; u, v}=(M^{t,\vartheta; u,
v}_s)_{s\in[t,T]}$ is a square integrable
${\widetilde{\mathbb{F}}^\pi}$-martingale with $M^{t,\vartheta; u,
v}_t=0$. Moreover, $M^{t,\vartheta; u, v}$ is supposed to be orthogonal
to the driving Brownian motion $B$, that is, their joint quadratic
variation process satisfies $[B,M^{t,\vartheta; u, v}]_s=0,
s\in[t,T]$. For the proof of the existence and the uniqueness of the
solution of such BSDE (\ref{BSDE}) it is similar to the classical case,
see also \cite{CFS} and references inside.
\end{longlist}

We have to emphasize here that since the filtration
${\widetilde{\mathbb{F}}^\pi}$ is not the Brownian one, but contains
it strictly, we cannot expect to have a solution of the above BSDE
with vanishing $M^{t,\vartheta; u, v}$. It is by now well known that,
under our assumptions on the coefficients $f$ and $\Phi$, a BSDE of
the above type has a unique solution $(Y^{t,\vartheta; u,
v},Z^{t,\vartheta;
u, v},M^{t,\vartheta; u, v})$. Moreover, considering the special form of
the filtration ${\widetilde{\mathbb{F}}^\pi}$, we can characterize
this solution as follows.

%re2.1 #&#
\begin{remark}\label{remark_bsde}We first observe that on each of the
subintervals $[t_{j-1},t_j)$, $
1\le j\le n$, formed by the partition $\pi=\{0=t_0<t_1<\cdots
<t_n=T\}$, the filtration $\widetilde{\mathbb{F}}^\pi$ coincides
with the Brownian one $({\cal F}^B_s)_{s\in[t_{j-1},t_j)}$ augmented
by the independent $\sigma$-field ${\cal H}_{j-1}$. Hence, on the
interval $[t_{j-1},t_j)$ we have the martingale representation
property for random variables from $L^2(\Omega,\widetilde{\cal
F}_{t_{j}-}^\pi,P)$ with respect to the
$\widetilde{\mathbb{F}}^\pi$-Brownian motion $B$. This has as
consequence that BSDE (\ref{BSDE}) can be solved over the time
intervals $[t_{j-1},t_j)$ with $dM^{t,\vartheta;u,v}_s=0, s\in[t_{j-1},
t_j)$. However, for this $Y^{t,\vartheta;u,v}_{t_{j}-}$ has to be determined
by backward iteration. In order to compute $Y^{t,\vartheta;u,v}_{t_{n}-}$,
we determine from BSDE (\ref{BSDE}) the jump of the c\`{a}dl\`{a}g
process $Y^{t,\vartheta; u, v}$ at time $t_{n}$:
\begin{eqnarray}
\triangle Y^{t,\vartheta; u, v}_{t_{n}} \bigl(:=Y^{t, \vartheta;
u, v}_{t_{n}}-Y^{t,\vartheta; u, v}_{t_{n}-}
\bigr)&=& \triangle M^{t,\vartheta;u, v}_{t_{n}}\nonumber\\
\eqntext{\mbox{that is } Y^{t,\vartheta; u, v}_{t_{n}-}=Y^{t,\vartheta; u, v}_{t_{n}}- \triangle
M^{t,\vartheta; u, v}_{t_{n}}.}
\end{eqnarray}

Taking into account that $ M^{t,\vartheta; u, v}$ is an
$\widetilde{\mathbb{F}}^\pi$-martingale, this yields
\[
Y^{t,\vartheta; u, v}_{t_{n}-}=E \bigl[Y^{t,\vartheta; u,
v}_{t_{n}} |
\widetilde{\cal F}_{t_{n}-}^\pi \bigr]\quad\mbox{and}\quad
\triangle M^{t,\vartheta; u, v}_{t_{n}}=Y^{t,\vartheta; u,
v}_{t_{n}}-E
\bigl[Y^{t,\vartheta; u, v}_{t_{n}} | \widetilde{\cal F}_{t_{n}-}^\pi
\bigr].
\]

Having now $Y^{t,\vartheta; u, v}_{t_{n}-}\in L^2(\Omega,
\widetilde{\cal F}_{t_{n}-}^\pi,P)$, we can consider BSDE
(\ref{BSDE}) over the time interval $[t_{n-1},t_{n})$ like a
classical one, with $dM^{t,\vartheta;u,v}_s=0,
s\in[t_{n-1},t_{n})$. By slving this BSDE over $[t_{n-1},t_{n})$,
we get, in particular, $Y^{t,\vartheta; u, v}_{t_{n-1}}$.
Iterating this argument, we see that
\[
Y^{t,\vartheta; u, v}_{t_{j}-}=E \bigl[Y^{t,\vartheta; u,
v}_{t_{j}} |
\widetilde{\cal F}_{t_{j}-}^\pi \bigr]\quad\mbox{and}\quad
\triangle M^{t,\vartheta; u,
v}_{t_{j}}=Y^{t,\vartheta; u, v}_{t_{j}}-E
\bigl[Y^{t,\vartheta; u,
v}_{t_{j}} | \widetilde{\cal F}_{t_{j}-}^\pi
\bigr],
\]
for all $t_j> t$, and $M^{t,\vartheta; u, v}$ is constant in
the intervals $[t_{j-1}\vee t,
t_j)$, $1\le j\le n$.
\end{remark}

%re2.2 #&#
\begin{remark}\label{classical} In the classical case, where the running
payoff function
$f(s,x,y,z,u,v)$ does not depend on $y$ and on $z$, the solution
$Y^{t,\vartheta;u,v}$ of BSDE~(\ref{BSDE}) takes the simple,
well-known form
\begin{eqnarray}
Y^{t,\vartheta;u,v}_s=E \biggl[\Phi \bigl(X^{t,\vartheta;u,
v}_T
\bigr)+\int_s^Tf \bigl(r,X^{t,\vartheta;u,v}_r,u_r,v_r
\bigr)\,dr | \widetilde {\cal F}_s^\pi \biggr], \nonumber\\
 \eqntext{s
\in[t,T], x\in R^d.}
\end{eqnarray}

From standard estimates for BSDEs of the type of equation
(\ref{BSDE}) we get, for all $p\ge2$, the existence of some
constant $C_p$ depending only $p$ and on the Lipschitz constants and
the bounds of the coefficients, such that, for all partitions $\pi$,
all initial\vadjust{\goodbreak} data $t\in[0,T], \vartheta,\vartheta'\in L^2(\Omega
,{\cal
F}^\pi_t,P;R^d)$ and all $u\in{\cal U}_{t,T},v\in{\cal V}_{t,T}$ it
holds, $P$-a.s.,
%
%
%e2.4 #&#
%
\begin{eqnarray}
\label{BSDE-estimates-1} %
\mathrm{(i)} &&\quad\bigl |Y^{t,\vartheta; u, v}_s
\bigr| \le C_p,\qquad s\in[t,T];
\nonumber
\\
\mathrm{(ii)}&&\quad E \biggl[ \biggl(\int_t^T
\bigl|Z^{t,\vartheta; u, v}_s \bigr|^2\,ds \biggr)^{p/2} \Big|
\widetilde{\cal F}_t^\pi \biggr]\le C_p;
\nonumber\\
\mathrm{(iii)}&&\quad E \biggl[\sup_{s\in[t,T]}
\bigl|Y^{t,\vartheta; u, v}_s- Y^{t,\vartheta'; u,v}_s \bigr|^p\\
&&\hspace*{22pt}{}+ \biggl(\int_t^T \bigl|Z^{t,\vartheta; u, v}_s-Z^{t,\vartheta'; u,
v}_s\bigr|^2
\,ds \biggr)^{p/2} \Big| \widetilde{\cal F}_t^\pi
\biggr]\nonumber\\
&&\quad\qquad \le C_p E \bigl[ \bigl|\vartheta-\vartheta'
\bigr|^p | \widetilde{\cal F}_t^\pi \bigr];
\nonumber
\end{eqnarray}
from where, in particular, for some constant $C\in R$,
%
%
%e2.5 #&#
%
\begin{eqnarray}
\label{BSDE-estimates-2} %
\mathrm{(i)}&&\quad \bigl|Y^{t,\vartheta; u, v}_t
\bigr| \le C,\qquad P\mbox{-a.s.};
\nonumber
\\[-8pt]
\\[-8pt]
\nonumber
\mathrm{(ii)}&&\quad \bigl|Y^{t,\vartheta; u, v}_t-Y^{t,\vartheta'; u, v}_t
\bigr| \le C \bigl(E \bigl[ \bigl|\vartheta-\vartheta' \bigr|^2 |
\widetilde{\cal F}_t^\pi \bigr] \bigr)^{1/2},
\qquad P\mbox{-a.s.}
\end{eqnarray}

For a game, in which the both Players 1 and 2 play along a partition
$\pi$ over a time interval $[t,T]$ and use the admissible controls
$u\in{\cal U}_{t,T}^\pi$ and $v\in{\cal V}_{t,T}^\pi$, we consider
the following pay-off functional:
\[
J^\pi(t,x;u,v)=Y^{t,x; u, v}_t,\qquad (t,x)\in[0,T]
\times R^d, (u,v)\in{\cal U}_{t,T}^\pi\times{
\cal V}_{t,T}^\pi.
\]
However, if we want to study the stochastic differential
game in a general frame, we can not consider games of the type
``control against control'', but we shall study games with
nonanticipative strategies with delay; for a more detailed
discussion the reader is referred to, for example, \cite{BCQ11}.
\end{remark}

Let us introduce the notion of nonanticipative strategies with
delay (NAD-strategies). They differ from the definitions given in
\cite{BCR04} and in \cite{BCQ11} and follow rather the spirit of the
definition given in \cite{Buckdahn-Li_Quincampoix-2011}, but now
extended to the stochastic case.
%
%de2.2 #&#
\begin{definition}[(NAD-strategies along the partition
$\pi$)] \label{NAD}
Let $\pi=\{0=t_0<t_1<\cdots<t_n=T\}$ $(n\ge1)$ an arbitrary partition
of the time interval $[0,T]$ and $t\in[0,T]$. We say that a mapping
$\beta\dvtx
{\cal U}_{t,T}^\pi\longrightarrow{\cal V}_{t,T}^\pi$ is an NAD-strategy
for Player~2 for the game over the time interval $[t,T]$ along the partition
$\pi$, if:

\begin{longlist}[(ii)]
\item[(i)] For all ${\widetilde{\mathbb{F}}^\pi}$-stopping times
$\tau\dvtx  \Omega\rightarrow\pi=\{t_0,t_1,\ldots,t_n\}$ it holds:
Whenever two controls $u,u'\in{\cal U}_{t,T}^\pi$ coincide $ds\,dP$-a.e.
on the stochastic interval $[[t,\tau]]$, then also
$\beta(u)_s=\beta(u')_s,  ds\,dP$-a.e. on $[[t,\tau]]$.

\item[(ii)] For all $0\le j\le n-1$, it holds that, whenever two
controls $u,u'\in{\cal U}_{t,T}^\pi$ coincide $ds\,dP$-a.e. on $[t,t_j]
\times\Omega$, then also $\beta(u)_s=\beta(u')_s, ds\,dP$-a.e. on $[t,
t_{j+1}]\times\Omega$.\vadjust{\goodbreak}
\end{longlist}

The set of all NAD-strategies for Player 2 over $[t,T]$ along
the partition $\pi$ is denoted by ${\cal B}_{t,T}^\pi$.

In an obvious symmetric way we define for Player 1 his set ${\cal
A}_{t, T}^\pi$ of NAD-strategies $\alpha\dvtx {\cal V}_{t,T}^\pi
\longrightarrow{\cal U}_{t,T}^\pi$ over the interval $[t,T]$ along
the partition $\pi$.
\end{definition}

Unlike the definitions in \cite{BCR04} and \cite{BCQ11}, the
delays for which we have this NAD-property (ii) in the above
definition is not considered as arbitrarily small for a given
partition $\pi$, but they are defined by the partition $\pi$. But,
however, in what follows we will study our game as the mesh of the
partition $\pi$ tends to zero.

The following result is crucial and it links our games defined
through a couple of admissible controls with those defined through
NAD-strategies.
%
%le2.1 #&#
\begin{lemma}\label{controls-NAD-strategies}
Let $\pi$ be any partition of the interval $[0,T]$ and $t\in[0,T]$.
Then, for all couples of NAD-strategies $(\alpha,\beta)\in{\cal
A}_{t,T}^\pi\times{\cal B}_{t,T}^\pi$, there is a unique couple of
admissible controls $(u,v)\in{\cal U}_{t,T}^\pi\times{\cal
V}_{t,T}^\pi$ such that $\alpha(v)=u$ and $\beta(u)=v$, $ds\,dP$-a.e.
on $[t,T]\times\Omega$.
\end{lemma}

In the above cited references \cite{BCR04,BCQ11} and
\cite{Buckdahn-Li_Quincampoix-2011} different definitions of
NAD-strategies were given, but the idea of the proof of the above
lemma remains similar. However, let us give it for the convenience of
the reader.

\begin{pf} Let $\pi=\{0=t_0<t_1<\cdots<t_n=T\}$ be a partition of
the interval $[0,T]$, and $(\alpha,\beta) \in\mathcal{A}_{t,T}^{\pi}
\times\mathcal{B}_{t, T}^{\pi}$. Let $t\in[t_i,t_{i+1})$. Then, due
to our definition of NAD strategies, $\alpha(v),\beta(u)$
restricted to $[t,t_{i+1}]$ depend only on $v\in\mathcal{V}_{t,
T}^{\pi}$ and $u\in\mathcal{U}_{t, T}^{\pi}$ restricted to the
interval $[t,t_i]$. But this interval is empty or a singleton, so
that $\alpha(v),\beta(u)$ restricted to the $[t,t_{i+1}]$ do not
depend on $v$ and $u$, respectively. Thus, putting for arbitrary
$u^0\in\mathcal{U}_{t, T}^{\pi},v^0\in\mathcal{V}_{t,T}^{\pi}$,
$u^1:=\alpha(v^0),v^1:=\beta(u^0)$, we get
\[
\alpha\bigl(v^1\bigr)=u^1, \qquad \beta
\bigl(u^1\bigr)=v^1\qquad ds\,dP\mbox{-a.s. on }
[t,t_{i+1}].
\]
Let us suppose now that we have constructed, for $j\ge2$,
$(u^{j-1}, v^{j-1})\in\mathcal{U}_{t, T}^{\pi}\times\mathcal{V}_{t,
t_{l}}^{\Pi}$ such that $\alpha(v^{j-1})=u^{j-1}$ and
$\beta(u^{j-1})=v^{j-1}$, $ds\,dP$-a.s. on $[t,t_{i+j-1}]$. Then we
set $u^{j}:=\alpha(v^{j-1}), v^{j}:=\beta(u^{j-1})$, and,
obviously, $(u^{j},v^{j})\in\mathcal{U}_{t,
T}^{\pi}\times\mathcal{V}_{t, T}^{\pi}$ is such that
$(u^{j},v^{j})=(u^{j-1},v^{j-1})$, $ds\,dP$-a.s. on $[t,t_{i+j-1}]$.
Thus, because of the NAD property [see Definition \ref{NAD}(ii)] of
$\alpha,\beta$, $u^{j}=\alpha(v^{j}), v^{j}=\beta(u^{j})$,
$ds\,dP$-a.s. on $[t,t_{i+j}]$. Consequently, iterating this argument
we obtain the existence of a couple $(u,v)\in{\cal
U}_{t,T}^\pi\times{\cal V}_{t,T}^\pi$ which satisfies the statement
of the lemma. Its uniqueness is an immediate consequence of the
above construction.
\end{pf}

Given a couple of NAD-strategies $(\alpha,\beta)\in{\cal A}_{t,
T}^\pi\times{\cal B}_{t,T}^\pi$ of the both players, the above lemma
allows to define the corresponding dynamics and the corresponding
pay-off functional through those of the associated\vadjust{\goodbreak} admissible
control processes. More precisely, for $(u,v)\in{\cal
U}_{t,T}^\pi\times{\cal V}_{t,T}^\pi$ such that $\alpha(v)=u$ and
$\beta(u)=v$, $ds\,dP$-a.e. on $[t,T]\times\Omega$, we define, for all
$\vartheta\in L^2(\Omega,{\cal F}_t^\pi,P;R^d)$ and $x\in R^d$,
\begin{eqnarray*}
X^{t,\vartheta; \alpha,\beta}&:=&X^{t,\vartheta; u, v}, \\
\bigl(Y^{t,\vartheta; \alpha,\beta},Z^{t,\vartheta; \alpha,\beta
},M^{t,\vartheta;
\alpha,\beta}
\bigr)&:=&\bigl(Y^{t,\vartheta; u, v},Z^{t,\vartheta; u,
v},M^{t,\vartheta; u, v}\bigr),
\\
J^\pi(t,x;\alpha,\beta)&:=&J^\pi(t,x;u,v).
\end{eqnarray*}

After the above preliminary discussion, we are now able to introduce
the upper and the lower value functions for the game over the time
interval $[t,T]$ along a partition $\pi$. We define the \textit{lower
value function along a partition} $\pi$~as
%
%
%e2.6 #&#
%
\begin{equation}
W^\pi(t,x):= \esssup_{\alpha\in{\mathcal
{A}}_{t,T}^\pi} \essinf_{\beta\in{\cal{B}}_{t,T}^\pi}
J^\pi(t,x;\alpha,\beta)
\end{equation}
and the \textit{upper one} as follows:
%
%
%e2.7 #&#
%
\begin{equation}
U^\pi(t,x):= \essinf_{\beta\in{\cal
{B}}_{t,T}^\pi} \esssup_{\alpha\in{\cal{A}^\pi}_{t,T}}J^\pi(t,x;
\alpha,\beta).
\end{equation}

Let us emphasize that the above lower and the upper value functions are
defined as a combination of essential supremum and essential infimum
over a bounded family of $\widetilde{\cal F}_t^\pi$-measurable
random variables $J^\pi(t,x;\alpha,\beta)$. Indeed, due to
(\ref{BSDE-estimates-2})(i),
\[
\bigl|J^\pi(t,x;\alpha,\beta) \bigr|= \bigl|Y^{t,\vartheta; \alpha
,\beta}_t \bigr| \le C,
\qquad P\mbox{-a.s., for all }(\alpha,\beta)\in{\cal A}_{t,T}^\pi
\times{\cal B}_{t,T}^\pi.
\]

Consequently, with the definitions of the essential infimum
and the essential supremum over families of random variables, given in
\cite{D} and \cite{DS} (see also \cite{KS2} for a more detailed discussion),
the upper and the lower value functions $W^\pi(t,x)$ and
$U^\pi(t,x)$ are, a priori, themselves also bounded,
$\widetilde{\cal F}_t^\pi$-measurable random variables. But,
combining arguments from \cite{Buckdahn-Li-2008} and
\cite{Buckdahn-Li_Quincampoix-2011}, we will be able to prove that
they are deterministic. However, for this proof we will have first
to establish a dynamic programming principle.

Let us finish this section with the following estimates for the
lower and the upper value functions, which are an immediate
consequence of the corresponding uniform estimates
(\ref{BSDE-estimates-2}) for the solution of BSDE (\ref{BSDE}).\vspace*{-3pt}

%le2.2 #&#
\begin{lemma}\label{estimates-W,U}
Under our standard assumptions on the coefficients $b,\sigma,f$ and
$\Phi$ there exists a constant $L\in R$ such that, for all partitions
$\pi$ of $[0,T]$ and all $t\in[0,T], x,x'\in R^d$,
%
%
%e2.8 #&#
%
\begin{eqnarray}
\label{estimates_W,U} %
\mathrm{(i)}&&\quad
\bigl|W^{\pi} (t,x) \bigr|+ \bigl|U^{\pi}(t,x) \bigr|\le L,
\nonumber
\\
\mathrm{(ii)} &&\quad \bigl|W^{\pi}(t,x)-W^{\pi}
\bigl(t,x' \bigr)\bigr |+ \bigl|U^{\pi}(t,x)-U^{\pi}
\bigl(t,x' \bigr) \bigr| \leq L\bigl |x-x' \bigr|,\\[-4pt] \eqntext{P
\mbox{-a.s.}}
\end{eqnarray}
\end{lemma}

%s3 #&#
\section{Lower and upper value functions along a partition}\label{sec3}

This section is devoted to the study of properties of the lower
and the upper value functions $W^\pi$ and $U^\pi$ defined along a
partition $\pi$\vadjust{\goodbreak} of the interval $[0,T]$. The main objectives in this
section are to prove that both functions, characterized in the
preceding section as random fields, are in fact deterministic, and
they satisfy a dynamic programming principle along the
partition $\pi$.\vspace*{-2pt}

%th3.1 #&#
\begin{theorem}\label{W,U deterministic} For any partition $\pi$ of
the interval $[0,T]$ and for all $(t,x)\in[0,T]\times R^d$, we have
$W^\pi(t,x)= E[W^\pi(t,x)], U^\pi(t,x)=E[U^\pi(t,x)]$, $P$-a.s.\vspace*{-2pt}
\end{theorem}

%re3.1 #&#
\begin{remark} A consequence of this theorem is that, by identifying
$W^\pi(t, x) :=E[W^\pi(t,x)], U^\pi(t,x):=E[U^\pi(t,x)],
(t,x)\in[0,T]\times R^d$, the lower and the upper value functions along
a partition $\pi$
$W^\pi$ and $U^\pi$ can be regarded as deterministic functions.\vspace*{-2pt}
\end{remark}

The proof of the above theorem is strongly inspired by that of
Proposition~3.1 in~\cite{Buckdahn-Li-2008} and uses heavily the
structure of our underlying probability space $(\Omega,{\cal F},P)$.
We only give the proof for $W^\pi(t,x)$, for some arbitrarily fixed
$(t,x)\in[0,T]\times R^d$. The proof for $U^\pi(t,x)$ is analogous and
won't be given here.

Let the partition $\pi$ of the interval $[0,T]$ be of the form $\pi=
\{0=t_0<t_1<\cdots<t_n=T\}$ and let $1\le j\le n$ be such that $t\in
[t_{j-1},t_{j})$. Recalling that $W^\pi(t,x)$ is an $\widetilde{\cal
F}_t^\pi$-measurable random variable, it follows from the definition
of the $\sigma$-field $\widetilde{\cal F}_t^\pi$ that, $W^\pi(t,x)$ $P$-a.s.
coincides with a measurable functional
$W^\pi(t,x)(\zeta^{(j-1)},B^{(t)})$ of $\zeta^{(j-1)}=(\zeta
_1,\ldots,
\zeta_{j-1})$ of the first $j-1$ components of the coordinate
process $\zeta= (\zeta_\ell)_{\ell\ge1}$ on $\Omega_1$ and the
Brownian motion $B^{(t)}= (B_s)_{s\in[0,t]}$ defined over $\Omega_2$
and restricted to the time interval $[0,t]$.\looseness=-1

Let $H_t$ be the Cameron--Martin space of all absolutely continuous
functions $h\in C([0,T];R^d)$ which derivative $\dot{h}$ is square
integrable and satisfies $\dot{h}_s=0, ds$-a.e. on $[t,T]$, and let
us denote by $\Omega_1^{(j-1)}$ the set of all sequences
$\rho=(\rho_\ell= (\rho_{\ell,1},\rho_{\ell,2}))_{\ell\ge
1}\in\Omega_1$, such that $\rho_\ell=0, \ell\ge j$. Given any
$(a,h)\in\Omega_1^{(j-1)} \times H_t$, we define the transformation
$\tau_{a,h}\dvtx \Omega\rightarrow\Omega$ by putting
$\tau_{a,h}(\rho,\omega'):=(\rho+a,\omega'+h)
(=((\rho_\ell+a_\ell)_{\ell\ge1},\omega'+h) )$,
$(\rho,\omega')\in\Omega=\Omega_1\times\Omega_2$. Such defined
transformation is bijective, $\tau_{a,h}^{-1}=\tau_{-a,-h}$,
$(a,h)\in\Omega_1^{(j-1)} \times H_t$, and its law $P\circ
[\tau_{a,h}]^{-1}$ is equivalent to $P$. Indeed, the law $P\circ
[\tau_{a,h}]^{-1}$ has with respect to $P$ the density
\[
L_{a,h}=\exp \biggl\{\langle a,\zeta\rangle+\int_0^t
\dot{h}_s\,dB_s- \frac12 \biggl(|a|^2+\int
_0^t|\dot{h}_s|^2\,ds
\biggr) \biggr\},
\]
where
\begin{eqnarray*}
\langle a,\zeta\rangle&:=&\sum_{\ell\ge1}a_\ell
\zeta_\ell= \sum_{\ell=1}^{j-1}a_\ell
\zeta_\ell \biggl(=\sum_{1\le\ell\le
j-1,i=1,2}a_{\ell,i}
\zeta_{\ell,i} \biggr)\quad \mbox{and}
\\
|a|^2&=&\sum_{\ell\ge1}|a_\ell|^2=
\sum_{\ell=1}^{j-1} |a_\ell|^2
\biggl(=\sum_{1\le\ell\le j-1,i=1,2}|a_{\ell,i}|^2
\biggr),
\end{eqnarray*}
$a=(a_\ell= (a_{\ell,1},a_{\ell,2}))_{\ell\ge1}\in
\Omega
_1^{(j-1)}$. We observe that the density $L_{a,h}$ is $\widetilde{\cal
F}_t^{\pi}$-measurable and belongs to $L^p(\Omega,{\cal F},P)$, for
all $p\ge1$.

The following lemma is essential for the proof that $W(t,x)$ is
\mbox{deterministic}.\vspace*{-2pt}

%le3.1 #&#
\begin{lemma}\label{xi_deterministic}
Let $\xi\in L^0(\Omega,\widetilde{{\cal F}}_t^\pi,P)$ be a random
variable which, for all
$(a,h)\in\Omega_1^{(j-1)}\times H_t$, is invariant
with respect to all transformations $\tau_{a,h}\dvtx \Omega\rightarrow
\Omega$,
that is, $\xi\circ\tau_{a,h}=\xi$, $P$-a.s. Then, there exists some
deterministic real number $c\in R$, such that $\xi=c, P$-a.s.\vspace*{-2pt}
\end{lemma}
\begin{pf} Let $\xi\in L^0(\Omega,\widetilde{\cal F}_t^\pi,P)$ be invariant
with respect to all transformations
$\tau_{a,h}\dvtx \Omega\rightarrow\Omega$,
$(a,h)\in\Omega_1^{(j-1)}\times H_t$. Then, for all $(a,h)\in
\Omega_1^{(j-1)} \times H_t$ and all bounded Borel functions
$g\dvtx R\rightarrow R$,
%
%
%e3.1 #&#
%
\begin{eqnarray}
&&E \bigl[g(\xi) \bigr]\nonumber \hspace*{-35pt}\\
&&\qquad= E \bigl[g(\xi\circ \tau_{a,h})
\bigr]
\\
&& \qquad = E \biggl[g(\xi)\exp \biggl\{\langle a,\zeta\rangle +\int
_0^t\dot{h}_s\,dB_s
\biggr\} \biggr] \cdot\exp \biggl\{-\frac12 \biggl(|a|^2+\int
_0^t|\dot {h}_s|^2\,ds
\biggr) \biggr\},\nonumber\hspace*{-35pt}
\end{eqnarray}
that is,
%
%
%e3.2 #&#
%
\begin{eqnarray}
& &E \Biggl[g(\xi)\exp \Biggl\{\sum_{\ell=1}^{j-1}a_\ell
\zeta_\ell + \int_0^t
\dot{h}_s \,dB_s \Biggr\} \Biggr]
\nonumber
\\
&&\qquad = E \bigl[g(\xi) \bigr]\cdot\exp \biggl\{\frac12 \biggl(|a|^2+
\int_0^t| \dot{h}_s|^2
\,ds \biggr) \biggr\}
\\
&&\qquad = E \bigl[g(\xi) \bigr]\cdot E \Biggl[\exp \Biggl\{\sum
_{\ell=1}^{j-1}a_\ell \zeta_\ell+
\int_0^t\dot{h}_s
\,dB_s \Biggr\} \Biggr]
\nonumber
\end{eqnarray}
for all $a_\ell\in R^2, 1\le\ell\le j-1$, and all $h\in H_t$, from
where we deduce that $\xi$ is independent of
$(\zeta^{(j-1)}=(\zeta_1,\ldots,\zeta_{j-1}), B^{(t)}=(B_s)_{s
\in[0,t]})$ and, hence also of $\widetilde{\cal F}_t^\pi=\sigma\{
\zeta^{(j-1)},B^{(t)}\}$. But this means that $\xi$ as an $\widetilde
{{\cal F}}_t^\pi$-measurable random variable is independent of itself. The
statement of the lemma follows now easily.\vspace*{-2pt}
\end{pf}

\begin{pf*}{Proof of Theorem \ref{W,U deterministic}} In order to be
able to
conclude our theorem form the above lemma, we only have to show that
the random variable $W^\pi(t,x)$ is invariant with respect to the
transformations $\tau_{a,h}\dvtx  \Omega\rightarrow\Omega$, for all
$(a,h)\in\Omega_1^{(j-1)} \times H_t$. For showing this, we fix
arbitrarily $(a,h)\in\Omega_1^{(j-1)} \times H_t$ and we proceed in
an analogous spirit as that in the proof of Proposition 3.1 in
\cite{Buckdahn-Li-2008}. But, however, the framework is different
here.\vadjust{\goodbreak}

\textit{Step} 1. Given a couple of admissible controls
$(u,v)\in{\cal U}_{t,T}^\pi\times{\cal V}_{t,T}^\pi$, we notice
that also the transformed couple
$(u\circ\tau_{a,h},v\circ\tau_{a,h})$ belongs to ${\cal
U}_{t,T}^\pi\times{\cal V}_{t,T}^\pi$. Indeed, having
$t\in[t_{j-1},t_{j})$,
\begin{eqnarray*}
u_s&=&u_j \bigl(s,(\zeta_1,\ldots,
\zeta_{j-1},\zeta_{j,1}, B_{{\cdot\wedge s}})
\bigr)I_{[t,t_j)}(s)
\\[-2pt]
&&{}+\sum_{\ell=j+1}^nu_\ell
\bigl(s,(\zeta_1,\ldots, \zeta_{\ell-1},\zeta_{\ell,1},B_{{\cdot\wedge s}})
\bigr)I_{[t_{\ell
-1},t_\ell)} (s)\qquad ds\,dP \mbox{-a.e.},
\end{eqnarray*}
for measurable functionals $u_\ell, 1\le\ell\le n$, the
transformed control process $u\circ\tau_{a,h}$ takes the form
%
%
%e3.3 #&#
%
\begin{eqnarray}
&& u_s\circ\tau_{a,h} \nonumber\\[-2pt]
&&\qquad= u_j
\bigl(s,(\zeta_1+a_1,\ldots, \zeta_{j-1}+a_{j-1},
\zeta_{j,1},B_{{\cdot\wedge s}}+h_{\cdot\wedge t}) \bigr)I_{[t,t_j)}(s)
\nonumber
\\[-9pt]
\\[-9pt]
\nonumber
&&\qquad\quad{}+ \sum_{\ell=j+1}^nu_\ell
\bigl(s,(\zeta_1+a_1,\ldots, \zeta_{j-1}+a_{j-1},
\zeta_j,\ldots,\zeta_{\ell-1},\zeta_{\ell,1},\nonumber\\[-2pt]
&&\hspace*{220pt}\qquad{} B_{{\cdot\wedge s}}+h_{\cdot\wedge t}) \bigr)I_{[t_{\ell-1},t_\ell)}
(s),\nonumber
\end{eqnarray}
$ds\,dP\mbox{-a.e.,}$ from where we see that also
$u\circ\tau_{a,h}$ is an admissible control for Player 1; the
symmetric argument shows that $v\circ\tau_{a,h} \in{\cal
V}_{t,T}^\pi$. Applying now the transformation to the forward
equation (\ref{SDE}) and taking into account that the increments of
the Brownian motion after $t$ are not changed by the transformation:
$(B_s-B_t)\circ\tau_{a,h} =B_s-B_t, s\in[t,T]$ (Indeed, recall
that $\dot{h}_s=0, ds$-a.e. on $[t,T]$), we obtain from the
uniqueness of the solution of SDE (\ref{SDE}) that
$X^{t,x;u,v}_s\circ\tau_{a,h}=X^{t,x;u(\tau_{a,h}),v(\tau_{a,h})}_s,
s\in[t,T]$, $P$-a.s. Let us now apply the transformation
$\tau_{a,h}$ to BSDE (\ref{BSDE}). With the argument already used for
its application to the forward SDE we see that BSDE (\ref{BSDE})
becomes
%
%
%e3.4 #&#
%
\begin{eqnarray}
\label{transformed_BSDE} %
&&dY^{t,x; u, v}_s
\circ\tau_{a,h}\nonumber\\
&&\qquad=-E \bigl[f \bigl(s,X^{t,x; u(\tau_{a,h}),
v(\tau_{a,h})}_s,
Y^{t,x; u, v}_s\circ \tau_{a,h}, Z^{t,x; u, v}_s
\circ\tau_{a,h},
\nonumber
\\
& &\qquad\hspace*{164pt}{} u_s(\tau_{a,h}), v_s(
\tau_{a,h}) \bigr) | \widetilde{\cal F}^\pi_s
\bigr] \,ds
\\
&&\hspace*{10pt}\qquad{}+Z^{t,x;u, v}_s\circ \tau_{a,h}
\,dB_s+dM^{t,x; u,v}_s\circ \tau_{a,h},
\nonumber\\
&&Y^{t,x; u, v}_T\circ\tau_{a,h} = E \bigl[\Phi
\bigl(X^{t,x; u(\tau_{a,h}),
v(\tau_{a,h})}_T \bigr) | \widetilde{\cal F}^\pi_{T}
\bigr].
\nonumber
\end{eqnarray}

We remark that
(i) $(Y^{t,x; u, v}\circ\tau_{a,h},Z^{t,x; u, v}\circ\tau_{a,h})
\in{\cal S}^2_{\widetilde{\mathbb{F}}^\pi}(t,T;R)\times
L_{\widetilde{\mathbb{F}}^\pi}^2(t,\break T; R^d)$. Indeed, the
$\widetilde{\mathbb{F}}^\pi$-adaptedness of the transformed process
can be proved directly, and the square integrability follows from
standard $L^p$-estimates for the solutions of BSDEs:
\begin{eqnarray*}
&&E \biggl[\sup_{s\in[t,T]} \bigl|Y^{t,x; u,
v}_s\circ
\tau_{a,h} \bigr|^2+\int_t^T\bigl|Z^{t,x; u, v}_s
\circ\tau _{a,h}\bigr|^2\,ds \biggr]
\\
&&\qquad= E \biggl[ \biggl(\sup_{s\in[t,T]}\bigl|Y^{t,x; u,
v}_s\bigr|^2+
\int_t^T\bigl|Z^{t,x; u, v}_s\bigr|^2
\,ds \biggr)L_{a,h} \biggr]
\\
&&\qquad\le C \bigl(E\bigl[L_{a,h}^2\bigr]
\bigr)^{1/2} \biggl(E \biggl[\sup_{s\in[t,T]}
\bigl|Y^{t,x; u,v}_s\bigr|^4+ \biggl(\int
_t^T\bigl|Z^{t,x; u, v}_s\bigr|^2
\,ds \biggr)^2 \biggr] \biggr)^{1/2}
\\
&&\qquad< +\infty.
\end{eqnarray*}

On the other hand, the fact $L_{a,h}\in
L^2(\Omega,\widetilde{\cal F}_t^\pi,P)$ has as consequence that also
the transformed $(\widetilde{\mathbb{F}}^\pi,P)$-martingale
$M^{t,x;u,v}\circ\tau_{a,h}=(M^{t,x;u,v}_s\circ\tau_{a,h})_{s\in[t,T]}$
is again an $(\widetilde{\mathbb{F}}^\pi,P)$-martingale. Indeed, for
$t\le s\le T$ and $\xi\in L^\infty(\Omega,\widetilde{\cal
F}_s^\pi,P)$, also $\xi\circ\tau_{-a,-h}\in L^\infty(\Omega,
\widetilde{\cal F}_s^\pi,P)$, and
%
%
%e3.5 #&#
%
\begin{eqnarray}
&& E \bigl[ \bigl(M^{t,x;u,v}_T-M^{t,x;u,v}_s
\bigr)\circ\tau_{a,h} \cdot\xi \bigr]
\nonumber
\\
&&\qquad = E \bigl[ \bigl(M^{t,x;u,v}_T-M^{t,x;u,v}_s
\bigr) \cdot\xi\circ\tau_{-a,-h}\cdot L_{a,h} \bigr]
\\
&&\qquad = E \bigl[E \bigl[M^{t,x;u,v}_T-M^{t,x;u,v}_s
| \widetilde{\cal F}_s^\pi \bigr] \cdot\xi\circ
\tau_{-a,-h} L_{a,h} \bigr]=0.
\nonumber
\end{eqnarray}

Consequently, $M^{t,x;u,v}\circ\tau_{a,h}$ is an
$(\widetilde{\mathbb{F}}^\pi,P)$-martingale; its square
integrability follows from an argument similar to that for $(Y^{t,x;
u, v}\circ\tau_{a,h},Z^{t,x; u, v}\circ\tau_{a,h})$, (recall the
explicit representation of $M^{t,x;u,v}$ in terms of $Y^{t,x;u,v}$,
which implies the $L^p$-integrability of $M^{t,x;u,v}$ for all $p\ge
1.$) and its orthogonality to $B$ stems from the fact that it is a
pure jump martingale.

This shows that $(Y^{t,x; u, v}\circ\tau_{a,h},Z^{t,x; u, v}\circ
\tau_{a,h},
M^{t,x;u,v}\circ\tau_{a,h})$ is a solution of BSDE (\ref{BSDE}) with
the couple
of admissible controls $(u(\tau_{a,h}),v(\tau_{a,h}))$. From the
uniqueness of the
solution of this BSDE it then follows that
%
%
%e3.6 #&#
%
\begin{eqnarray}
& & \bigl(Y^{t,x; u, v}\circ \tau_{a,h},Z^{t,x; u, v}\circ
\tau_{a,h}, M^{t,x;u,v} \circ\tau_{a,h} \bigr)
\nonumber
\\[-8pt]
\\[-8pt]
\nonumber
&&\qquad = \bigl(Y^{t,x; u(\tau_{a,h}), v(\tau_{a,h})},Z^{t,x;
u(\tau_{a,h}), v(\tau_{a,h})},M^{t,x;u(\tau_{a,h}),v(\tau_{a,h})}
\bigr),
\end{eqnarray}
and, in particular, it follows that
\[
J^\pi(t,x;u,v)\circ\tau_{a,h}=J^\pi
\bigl(t,x;u(\tau _{a,h}),v(\tau_{a,h})\bigr),\qquad P
\mbox{-a.s.}
\]

\textit{Step} 2. Let us translate in this step the result of
step 1 to couples of NAD strategies. For $\beta\in{\cal B}_{t,T}^\pi$
we define
$\beta_{a,h}(u):=\beta(u(\tau_{-a,-h}))(\tau_{a,h}), u\in{\cal
U}_{t,T}^\pi$.
For such defined mapping $\beta_{a,h}\dvtx {\cal U}_{t,T}^\pi\rightarrow
{\cal V}_{t,T}^\pi$ it can be verified in a straight-forward manner
that it
belongs to ${\cal B}_{t,T}^\pi$. We also observe that $(\beta
_{-a,-h})_{a,h}=\beta$.
A symmetric definition allows to introduce $\alpha_{a,h}\in{\cal
A}_{t,T}^\pi$, for
$\alpha\in{\cal A}_{t,T}^\pi$ and to get $(\alpha
_{-a,-h})_{a,h}=\alpha$.

Given a couple of NAD-strategies $(\alpha,\beta)\in{\cal
A}_{t,T}^\pi
\times{\cal
B}_{t,T}^\pi$, let us denote by $(u,v)\in{\cal U}_{t,T}^\pi\times
{\cal
V}_{t,T}^\pi$ the couple of admissible controls associated with
through Lemma
\ref{controls-NAD-strategies}. Then
\begin{eqnarray*}
\alpha_{a,h} \bigl(v(\tau_{a,h}) \bigr)&=&\alpha(v) (
\tau_{a,h})=u(\tau_{a,h}) \quad\mbox{and}
\\
\beta_{a,h} \bigl(u(\tau_{a,h}) \bigr)&=&\beta(u) (
\tau_{a,h})=v(\tau_{a,h}).
\end{eqnarray*}
Consequently, the couple $(u(\tau_{a,h}),v(\tau_{a,h}))\in
{\cal
U}_{t,T}^\pi\times{\cal V}_{t,T}^\pi$ is associated with $(\alpha
_{a,h},\beta_{a,h})$
through Lemma \ref{controls-NAD-strategies}, and from step 1 we get
%
%
%e3.7 #&#
%
\begin{eqnarray}
J^\pi(t,x;\alpha,\beta)\circ \tau_{a,h} & = &
J^\pi(t,x;u,v)\circ \tau_{a,h}= J^\pi
\bigl(t,x;u(\tau_{a,h}),v(\tau_{a,h}) \bigr)
\nonumber
\\[-8pt]
\\[-8pt]
\nonumber
& = & J^\pi(t,x;\alpha_{a,h},
\beta_{a,h}),\qquad P \mbox{-a.s.}
\end{eqnarray}

\textit{Step} 3. Using the definition of the esssup and
the essinf over a family of random variables as well as the fact
that the transformation $\tau_{a,h}$ is invertible and its law
$P\circ[\tau_{a,h}]^{-1}$ is equivalent to $P$, we show that
%
%
%e3.8 #&#
%
\begin{eqnarray}
&&W^\pi(t,x)\circ\tau_{a,h}\nonumber\\
&&\qquad =  \bigl(
\esssup_{\alpha\in
{\mathcal{A}}_{t,T}^\pi} \essinf_{\beta\in{\cal{B}}_{t,T}^\pi} J^\pi(t,x;\alpha,\beta
) \bigr)\circ\tau_{a,h}
\\
&&\qquad =  \esssup_{\alpha\in{\mathcal{A}}_{t,T}^\pi} \essinf_{\beta\in{\cal{B}}_{t,T}^\pi}
\bigl(J^\pi (t,x;\alpha,\beta)\circ \tau_{a,h} \bigr),\qquad
P \mbox{-a.s.}\nonumber
\end{eqnarray}

Consequently, by combining the results of the previous steps
and by considering that,
thanks to step 2, $\{\alpha_{a,h}, \alpha\in{\cal A}_{t,T}^\pi\}
={\cal
A}_{t,T}^\pi$
and $\{\beta_{a,h}, \beta\in{\cal B}_{t,T}^\pi\}={\cal B}_{t,T}^\pi$,
we obtain
%
%
%e3.9 #&#
%
\begin{eqnarray}
W^\pi(t,x)\circ\tau_{a,h} & = &
\esssup_{\alpha\in{\mathcal
{A}}_{t,T}^\pi} \essinf_{\beta\in{\cal{B}}_{t,T}^\pi} \bigl(J^\pi (t,x;
\alpha,\beta)\circ \tau_{a,h} \bigr)
\nonumber
\\
& = & \esssup_{\alpha\in{\mathcal{A}}_{t,T}^\pi} \essinf_{\beta\in{\cal{B}}_{t,T}^\pi} J^\pi(t,x;
\alpha _{a,h},\beta_{a,h})
\\
& = & W^\pi(t,x),\qquad P\mbox{-a.s.}
\nonumber
\end{eqnarray}

By combining this result with Lemma \ref{xi_deterministic}, we complete
the proof.
\end{pf*}

As an immediate consequence of Lemma \ref{estimates-W,U} and the above
result that
the lower and the upper value functions along a partition are
deterministic, we have
the following result.

%le3.2 #&#
\begin{lemma}\label{estimates+W,U} There exists a constant $L\in R$
which does not
depend on the partition $\pi$ of the interval $[0,T]$, such that, for
all $t\in[0,T]$,
$x,x'\in R^d$,
%
%
%e3.10 #&#
%
\begin{eqnarray}
\mathrm{(i)}&&\quad \bigl|W^{\pi} (t,x) \bigr|+ \bigl|U^{\pi}(t,x)
\bigr| \le L,
\nonumber
\\[-8pt]
\\[-8pt]
\nonumber
\mathrm{(ii)} &&\quad\bigl|W^{\pi}(t,x)-W^{\pi}
\bigl(t,x' \bigr) \bigr|+ \bigl|U^{\pi}(t,x)-U^{\pi}
\bigl(t,x' \bigr)\bigr | \leq L \bigl|x-x' \bigr|.
\end{eqnarray}
\end{lemma}

After having proved that the lower and the upper value functions along
a partition
$\pi$ are deterministic, our objective is now to show that, with
respect to the
points of the partition they satisfy the DPP. A key role will be played
here by the notion of backward stochastic semigroup,
introduced by Peng in \cite{Pe1}.

Given a partition $\pi=\{0=t_0<t_1<\cdots<t_n=T\}$ of the interval $[0,T]$,
initial data $(t,x)\in[0,T)\times R^d$, a positive $\delta<T-t$ and a
couple of\vadjust{\goodbreak}
admissible control processes $(u,v)\in{\cal U}_{t,t+\delta}^\pi
\times
{\cal V}_{t,
t+\delta}^\pi$ as well as a random variable $\eta\in L^2(\Omega
,{\cal
F}^\pi_{t+\delta},
P)$, we define the backward stochastic semigroup
\[
G^{t,x;u,v}_{s,t+\delta}(\eta):=\overline{Y}_s^{u,v},
\qquad s\in[t,t+\delta],
\]
through the BSDE with time horizon $t+\delta$,
%
%
%e3.11 #&#
%
\begin{equation}
\label{semigroup} \cases{d \overline{Y}^{u, v}_s = -E
\bigl[f \bigl(s,X^{t,\vartheta; u, v}_s, \overline{Y}^{u, v}_s,
\overline{Z}^{u, v}_s, u_s, v_s
\bigr) | \widetilde{\cal F}^\pi_s \bigr] \,ds\vspace*{2pt}\cr
\hspace*{39pt}{} +
\overline{Z}^{u, v}_s \,dB_s+d
\overline{M}^{u,
v}_s, \vspace*{2pt}
\cr
\overline{Y}^{u, v}_T
= E \bigl[\eta | \widetilde{\cal F}^\pi_{t+\delta} \bigr],}
\end{equation}
and its unique solution $(\overline{Y}^{u, v}, \overline
{Z}^{u, v},
\overline{M}^{u, v})\in{\cal S}^2_{\widetilde{\mathbb{F}}^\pi
}(t,t+\delta;R)
\times L_{\widetilde{\mathbb{F}}^\pi}^2(t,t+\delta;R^d)\times
{\cal M}_{\widetilde{\mathbb{F}}^\pi}^2(t,t+\delta;R)$ with
$[B,\overline{M}^{u, v}]_s=0, s\in[t,T]$ and $\overline{M}^{u, v}_t=0$,
where $X^{t,\vartheta; u, v}$ is the solution of SDE (\ref{SDE}).

From the discussion made in the frame of Remark \ref{remark_bsde} it becomes
clear that if, for some point $t_j$ of the partition $\pi=\{
0=t_0<t_1<\cdots
<t_n=T\}$, $t_{j-1}\le t<t+\delta=t_j$ and $\eta$ is
$\widetilde{\cal F}^\pi_{t_j-}$-measurable, then $\overline{M}^{u, v}_s=0,
s\in[t,t_j]$.

The properties of the backward stochastic semigroup follow directly
from those of the BSDE through which it is defined, so that we won't discuss
separately here (refer to \cite{Pe1}, or \cite{Buckdahn-Li-2008}). The
notion of backward stochastic semigroup now allows to
study the DPP along a partition $\pi$ of the
time interval $[0,T]$.

%th3.2 #&#
\begin{theorem}\label{DPP} Let $\pi=\{0=t_0<t_1<\cdots<t_n=T\}$ be a
partition of the interval $[0,T]$, and let $t\in[t_i,t_{i+1})$ and
$x\in R^d$. Then, for all $i+1\le j\le n$, $P$-a.s.,
%
%
%e3.12 #&#
%
\begin{eqnarray}
\label{DPP-aa} %
W^{\pi}(t,x) &=& \esssup_{\alpha\in
\mathcal{A}^{\pi}_{t,t_j}}
\essinf_{\beta\in
\mathcal{B}^{\pi}_{t,t_j}}G^{t,x;\alpha,\beta}_{t,t_j} \bigl(W^{\pi
}
\bigl(t_j,X_{t_j}^{t,x;
\alpha,\beta} \bigr) \bigr),
\nonumber
\\[-8pt]
\\[-8pt]
\nonumber
U^{\pi}(t,x) &=& \essinf_{\beta\in\mathcal{B}^\pi_{t,t_j}}\esssup
_{\alpha
\in\mathcal{A}^\pi_{t,t_j}} G^{t,x;\alpha,\beta}_{t,t_j} \bigl(U^{\pi}
\bigl(t_j,X_{t_j}^{t,x;\alpha,\beta} \bigr) \bigr).
\end{eqnarray}
\end{theorem}

%re3.2 #&#
\begin{remark} The space ${\cal U}^\pi_{t,t_j}$ of admissible controls
for Player 1 for games over the time interval $[t,t_j]$ along the
partition $\pi$ is defined as the set of all control processes
$u\in{\cal U}^\pi_{t,T}$ restricted to the time interval $[t,t_j]$; the
space ${\cal V}^\pi_{t,t_j}$ of admissible controls for Player 2 is
defined analogously. The NAD-strategies for Player~2,
$\beta\in\mathcal{B}^\pi_{t,t_j}\dvtx {\cal U}^\pi_{t,t_j}\rightarrow
{\cal
V}^\pi_{t,t_j}$, are defined in the same manner as the NAD-strategies in
$\mathcal{B}^\pi_{t,T}$, with the only difference that we consider
$t_j$ instead $T=t_n$ as terminal horizon. The same is done in the
definition of the set $\mathcal{A}^\pi_{t,t_j}$ of NAD-strategies
for Player~1.
\end{remark}

The proof split into two lemmas for the lower value
function along the partition~$\pi$; it is similar for the upper value
function along the partition $\pi$. Let us fix
arbitrarily a partition\vadjust{\goodbreak} $\pi=\{0=t_0<t_1<\cdots<t_n=T\}$ of the interval
$[0,T]$, and let $t\in[t_i,t_{i+1})$,
$i+1\le j\le n$ and $x\in R^d$. We put
\[
\widetilde{W}_{t_j}^{\pi}(t,x) =\esssup_{\alpha\in\mathcal{A}^\pi_{t,t_j}}
\essinf_{\beta\in\mathcal{B}^\pi_{t,t_j}}G^{t,x;\alpha,
\beta}_{t,t_j}\bigl(W^{\pi}
\bigl(t_j,X_{t_j}^{t,x;\alpha,\beta}\bigr)\bigr).
\]

Obviously, $ \widetilde{W}_{t_j}^{\pi}(t,x)$
is a bounded, $\widetilde{\cal F}_t^\pi$-measurable random variable.

%le3.3 #&#
\begin{lemma}\label{DPP-1} Under the standard assumptions, we
have made on the coefficients it holds that $
\widetilde{W}_{t_j}^{\pi}(t,x)\le W^{\pi}(t,x)$, $P$-a.s.
\end{lemma}
\begin{pf}
\textit{Step} 1. Let us fix an arbitrary
$\varepsilon>0$. Then, we can find $\alpha_1^\varepsilon\in
{\cal A}_{t,t_j}^\pi$ such that
\[
\widetilde{W}_{t_j}^{\pi}(t,x)\le \essinf_{\beta\in\mathcal{B}^\pi_{t,t_j}}G^{t,x;\alpha
_1^\varepsilon,
\beta}_{t,t_j}
\bigl(W^{\pi}\bigl(t_j,X_{t_j}^{t,x;\alpha_1^\varepsilon,\beta}
\bigr)\bigr)+ \varepsilon,\qquad P\mbox{-a.s.}
\]
In order to verify this latter relation, we put
\[
I(\alpha): =\essinf_{\beta\in\mathcal{B}^\pi_{t,t_j}}G^{t,x;\alpha,\beta}_{t,t_j}
\bigl(W^{\pi}\bigl(t_j,X_{t_j}^{t,x;\alpha,\beta}
\bigr)\bigr),\qquad \alpha\in{\cal A}_{t,t_j}^\pi,
\]
and we note that, due to the
properties of the essential supremum over a family of random
variables, there is some sequence $(\alpha^k)_{k\ge1}\subset{\cal
A}_{t,t_j}^\pi$ such that
\[
\widetilde{W}_{t_j}^{\pi}(t,x)=\esssup_{\alpha\in
\mathcal{A}^\pi_{t,t_j}}I(
\alpha)=\sup_{k\ge1}I\bigl(\alpha^k\bigr),\qquad P
\mbox{-a.s.}
\]

Thus, putting $\triangle_k:=\{\widetilde{W}_{t_j}^{\pi}
(t,x)\le I(\alpha^k)+\varepsilon, \widetilde{W}_{t_j}^{\pi}(t,x)
> I(\alpha^\ell)+\varepsilon (1\le\ell\le k-1)\}\in\widetilde
{\cal F}_t^\pi$, $k\ge1$, we define a partition of $\Omega$,
and putting
\[
\alpha^\varepsilon_1(\cdot):=\sum_{k\ge1}
I_{\triangle_k}\alpha^k(\cdot)\dvtx {\cal V}_{t,t_j}^\pi
\rightarrow {\cal U}_{t,t_j}^\pi,
\]
we check easily that\vspace*{-1pt} $\alpha^\varepsilon_1$ is an
NAD-strategy in ${\cal A}_{t,t_j}^\pi$ and that
$\widetilde{W}_{t_j}^{\pi}(t,x)\le\sum_{k\ge1}
I_{\triangle_k}I(\alpha^k)+\varepsilon\le\sum_{k\ge1}
I_{\triangle_k}G^{t,x;\alpha^k,\beta_1}_{t,t_j}(W^{\pi}
(t_j,X_{t_j}^{t,x;\alpha^k,\beta_1}))+\varepsilon$, P-a.s., for all
$\beta_1\in{\cal B}_{t,t_j}^\pi$. Given an
arbitrary $\beta_1\in{\cal B}_{t,t_j}^\pi$, we let $(u^k,v^k)\in
{\cal U}_{t,t_j}^\pi\times{\cal V}_{t,t_j}^\pi$ be such that
$\alpha^k(v^k)=u^k, \beta_1(u^k)=v^k, ds\,dP$-a.e. on
$[t,t_j]\times\Omega$, and we introduce $(u_1,v_1):=\sum_{k\ge
1}I_{\triangle_k}(u^k,v^k)\in{\cal U}_{t,t_j}^\pi\times{\cal
V}_{t,t_j}^\pi$. Then, since for the
$\widetilde{\mathbb{F}}^\pi$-stopping time
$\tau_k=t_jI_{\triangle_k}+tI_{\triangle_k^c}$ the processes $u_1$
and $u^k$ coincide, $ds\,dP$-a.e. on $[[t,\tau_k]]$, also
$\beta_1(u^k)=\beta_1(u_1)$, $ds\,dP$-a.e. on $[[t,\tau_k]]$. Thus,
\[
\beta_1(u_1)=\sum_{k\ge1}I_{\triangle_k}
\beta_1 \bigl(u^k \bigr)= \sum
_{k\ge
1}I_{\triangle_k}v^k=v_1,
\qquad ds\,dP\mbox{-a.e. on} [t,t_j]\times\Omega,
\]
and with a symmetric argument we also have
\[
\alpha^\varepsilon_1(v_1)=\sum
_{k\ge1}I_{\triangle_k}\alpha^k(v_1)
=\sum_{k\ge1}I_{\triangle_k}\alpha^k
\bigl(v^k \bigr)=u_1,\qquad ds\,dP\mbox{-a.e. on }
[t,t_j]\times\Omega.
\]
This shows that the couple $(u_1,v_1)\in{\cal U}_{t,
t_j}^\pi\times{\cal V}_{t,t_j}^\pi$ is associated with
$(\alpha^\varepsilon_1,\beta_1)\in{\cal A}_{t,t_j}^\pi\times{\cal
B}_{t,t_j}^\pi$ by Lemma \ref{controls-NAD-strategies}.
Consequently, from the uniqueness of the solution of SDE (\ref{SDE})
we conclude with a standard argument that
\begin{eqnarray}
\sum_{k\ge1}I_{\triangle_k}X^{t,x;\alpha^k,\beta_1}=\sum
_{k\ge1} I_{\triangle_k}X^{t,x;u^k,v^k}=X^{t,x;u_1,v_1}=
X^{t,x;\alpha_1^\varepsilon,\beta_1} \nonumber\\
\eqntext{\mbox{on } [t,t_j], P\mbox{-a.s.}}
\end{eqnarray}
Similarly, using now the uniqueness of the solution of
BSDE defining the backward stochastic semigroup, we show that
\[
\sum_{k\ge1}I_{\triangle_k} \bigl(
\widetilde{Y}^{t,x;\alpha^k,\beta_1}, \widetilde{Z}^{t,x;\alpha^k,\beta_1},\widetilde{M}^{t,x;
\alpha^k,\beta_1}
\bigr)= \bigl(\widetilde{Y}^{t,x;\alpha_1^\varepsilon,
\beta_1},\widetilde{Z}^{t,x;\alpha_1^\varepsilon,
\beta_1},
\widetilde{M}^{t,x;\alpha_1^\varepsilon,\beta_1} \bigr),
\]
and recalling the definition of the backward stochastic
semigroup, we see that
\[
\sum_{k\ge1}I_{\triangle_k}G^{t,x;\alpha^k,\beta_1}_{t,t_j}
\bigl(W^{\pi} \bigl(t_j,X_{t_j}^{t,x;\alpha^k,\beta_1}
\bigr) \bigr)= G^{t,x;\alpha^\varepsilon_1,\beta_1}_{t,t_j} \bigl(W^{\pi}
\bigl(t_j, X_{t_j}^{t,x;\alpha^\varepsilon_1,\beta_1} \bigr) \bigr).
\]
Consequently, for all $\beta_1\in{\cal B}_{t,t_j}^\pi$,
%
%
%e3.13 #&#
%
\begin{eqnarray}
\widetilde{W}_{t_j}^{\pi}(t,x)& \le& \sum
_{k\ge1} I_{\triangle_k}I \bigl(
\alpha^k \bigr)+\varepsilon
\nonumber
\\
&\le& \sum_{k\ge1}I_{\triangle_k}
G^{t,x;\alpha^k,\beta_1}_{t,t_j} \bigl(W^{\pi} \bigl(t_j,X_{t_j}^{t,x;
\alpha^k,\beta_1}
\bigr) \bigr)+\varepsilon
\\
& = & G^{t,x;\alpha^\varepsilon_1,\beta_1}_{t,t_j} \bigl(W^{\pi}
\bigl(t_j, X_{t_j}^{t,x;\alpha^\varepsilon_1,\beta_1} \bigr) \bigr)+
\varepsilon,\qquad P\mbox{-a.s.}
\nonumber
\end{eqnarray}

Let us make now a special choice of $\beta_1\in{\cal
B}_{t,t_j}^\pi$. Given an arbitrary $\beta\in{\cal B}_{t,T}^\pi$ and
any $u_2\in{\cal U}_{t_j,T}^\pi$, we define for any $u_1\in{\cal
U}_{t,t_j}^\pi$ the process $u_1\oplus u_2:= u_1I_{[t,t_j]}+
u_2I_{(t_j,T]}\in{\cal U}_{t,T}^\pi$, and we put
\[
\beta_1(u_1):=\beta(u_1\oplus
u_2)_{|[t,t_j]},\qquad u_1\in{\cal
U}_{t,t_j}^\pi,
\]
the restriction of $\beta(u_1\oplus u_2)$ to the time
interval $[t,t_j]$. It can be easily verified that such defined
mapping $\beta_1\dvtx {\cal U}_{t,t_j}^\pi\rightarrow{\cal
V}_{t,t_j}^\pi$ belongs to ${\cal B}_{t, t_j}^\pi$, and thanks to
its nonanticipativity property it does not depend on the special
choice of $u_2$. Let us denote by
$(u_1^\varepsilon,v_1^\varepsilon)\in{\cal U}_{t,t_j}^\pi\times
{\cal V}_{t,t_j}^\pi$ the unique couple of control processes
associated with $(\alpha_1^\varepsilon,\beta_1)$ through Lemma
\ref{controls-NAD-strategies}.

\textit{Step} 2. After having proven in step 1 that
\[
\widetilde{W}_{t_j}^{\pi}(t,x)\le G^{t,x;
\alpha^\varepsilon_1,\beta_1}_{t,t_j}
\bigl(W^{\pi}\bigl(t_j,X_{t_j}^{t,
x;\alpha^\varepsilon_1,\beta_1}
\bigr)\bigr)+ \varepsilon,\qquad P\mbox{-a.s.},
\]
let us now estimate the expression $W^{\pi}(t_j,X_{t_j}^{t,
x;\alpha^\varepsilon_1,\beta_1})$ to which the backward stochastic
semigroup is applied at the right-hand side of the above estimate.
For this we consider a Borel partition ${\cal O}_k,k\ge1$, of $R^d$,
consisting of nonempty Borel sets ${\cal O}_k$ with diameter less
or equal to $\varepsilon$, and we fix arbitrarily in each of this
sets ${\cal O}_k$ an element $x_k$. With the arguments already
developed in step 1 we show that, for every $k\ge1$, there is
some $\alpha_2^k\in{\cal A}_{t_j,T}^\pi$ such that
\begin{eqnarray*}
W^{\pi}(t_j,x_k)&=&
\esssup_{\alpha_2 \in{\mathcal{A}}_{t_j,
T}^\pi}\essinf_{\beta_2 \in{\cal{B}}_{t_j,T}^\pi} J^{\pi}(t_j,x_k;
\alpha_2,\beta_2)
\\
&\le& \essinf_{\beta_2 \in{\cal{B}}_{t_j,T}^\pi} J^{\pi}\bigl(t_j,x_k;
\alpha_2^k,\beta_2\bigr)+\varepsilon, \qquad
P\mbox{-a.s.},
\end{eqnarray*}
and putting $\alpha_2^\varepsilon(\cdot):=\sum_{k\ge1}I
\{X^{t,x;\alpha_1^\varepsilon,\beta_1}_{t_j}\in{\cal O}_k\}
\alpha_2^k(\cdot)\dvtx  {\cal V}_{t_j,T}^\pi\rightarrow{\cal U}_{t_j,T}^\pi$
we obtain an NAD-strategy from ${\cal A}_{t_j,T}^\pi$. Indeed, the
sets $\{X^{t,x;\alpha_1^\varepsilon,\beta_1}_{t_j}\in{\cal O}_k\}$,
$k\ge1$, forming a partition of $\Omega$, belong to
\[
{\cal F}^\pi_{t_j-}={\cal F}^B_{t_j}
\vee{\cal H}_{j}= \widetilde{\cal F}^\pi_{t_j}.
\]
(We remark that the relation ${\cal F}^\pi_{s-}=
\widetilde{\cal F}^\pi_{s}$ only holds for points of the partition
$\pi$; this is also the reason, why we do not have a DPP which
does not use the points of the partition $\pi$). Thus, by combining
the arguments developed in step 1 with the Lipschitz property of
$W^{\pi}(t_j,\cdot)$ and $J^{\pi}(t_j,\cdot; \alpha, \beta)$ we can show that,
for all $\beta_2\in{\cal
B}_{t_j,T}^\pi$,
%
%
%e3.14 #&#
%
\begin{eqnarray}
\label{W_1}
&&W^{\pi} \bigl(t_j,X^{t,x;\alpha_1^\varepsilon,\beta_1}_{t_j}
\bigr)\nonumber\\
&&\qquad\le\sum_{k\ge1}I \bigl\{X^{t,x;\alpha_1^\varepsilon, \beta_1}_{t_j}
\in{\cal O}_k \bigr\}W^{\pi}(t_j,x_k)+L
\varepsilon
\nonumber
\\
&&\qquad \le \sum_{k\ge1}I \bigl\{X^{t,x;\alpha_1^\varepsilon,\beta_1
}_{t_j}
\in{\cal O}_k \bigr\}J^{\pi} \bigl(t_j,x_k;
\alpha_2^k,\beta_2 \bigr)+ (L+1)\varepsilon
\\
&&\qquad \le \sum_{k\ge1}I \bigl
\{X^{t,x;\alpha_1^\varepsilon, \beta_1}_{t_j} \in{\cal O}_k \bigr
\}J^{\pi} \bigl(t_j,X^{t,x;\alpha_1^\varepsilon, \beta_1}_{t_j};
\alpha_2^k,\beta_2 \bigr)+ (2L+1)\varepsilon
\nonumber\\
&&\qquad =  J^{\pi} \bigl(t_j,X^{t,x;\alpha_1^\varepsilon, \beta
_1}_{t_j};
\alpha _2^\varepsilon, \beta_2 \bigr)+ (2L+1)
\varepsilon,\qquad P\mbox{-a.s.}
\nonumber
\end{eqnarray}
For our arbitrarily chosen $\beta\in{\cal B}_{t,T}^\pi$ we
put $\beta_2^\varepsilon(u_2):=\beta(u_1^\varepsilon\oplus
u_2)_{|[t_j,T]}\in{\cal V}_{t_j,T}^\pi$, $u_2\in{\cal
U}_{t_j,T}^\pi$. Obviously, $\beta_2^\varepsilon\in{\cal
B}_{t_j,T}^\pi$. Let us denote by
$(u_2^\varepsilon,v_2^\varepsilon)\in{\cal
U}_{t_j,T}^\pi\times{\cal V}_{t_j,T}^\pi$ the unique couple of
control processes associated with
$(\alpha_2^\varepsilon,\beta_2^\varepsilon)$ through Lemma
\ref{controls-NAD-strategies}. Then, defining
$\alpha^\varepsilon\in{\cal A}_{t,T}^\pi$ by setting
\[
\alpha^\varepsilon(v):=\alpha^\varepsilon_1(v_{|[t,t_j]})
\oplus\alpha^\varepsilon_2(v_{|(t_j,T]}),\qquad v\in{\cal
V}_{t,T}^\pi,
\]
we see that, for $(u^\varepsilon,v^\varepsilon
):=(u^\varepsilon_1
\oplus u^\varepsilon_2, v^\varepsilon_1\oplus v^\varepsilon_2)\in
{\cal U}_{t,T}^\pi\times{\cal V}_{t,T}^\pi$,
\begin{eqnarray*}
\alpha^\varepsilon\bigl(v^\varepsilon\bigr)&=&\alpha^\varepsilon
_1\bigl(v^\varepsilon_1\bigr) \oplus
\alpha^\varepsilon_2\bigl(v^\varepsilon_2\bigr)
=u^\varepsilon_1\oplus u^\varepsilon_2=u^\varepsilon,
\\
\beta^\varepsilon\bigl(u^\varepsilon\bigr)&=&\beta^\varepsilon
\bigl(u^\varepsilon _1\oplus u^\varepsilon_2
\bigr)=\beta_1\bigl(u^\varepsilon_1\bigr) \oplus
\beta^\varepsilon_2 \bigl(u^\varepsilon_2
\bigr)=v^\varepsilon_1\oplus v^\varepsilon
_2=v^\varepsilon.
\end{eqnarray*}

Consequently, with the choice $\beta_2=\beta_2^\varepsilon$,
we have
%
%
%e3.15 #&#
%
\begin{eqnarray}
W^{\pi} \bigl(t_j,X^{t,x;\alpha_1^\varepsilon,\beta_1}_{t_j}
\bigr) &\le &J^{\pi} \bigl(t_j,X^{t,x;\alpha_1^\varepsilon, \beta_1}_{t_j};
\alpha _2^\varepsilon, \beta_2^\varepsilon
\bigr)+ (2L+1)\varepsilon
\nonumber
\\
& = & J^{\pi} \bigl(t_j,X^{t,x;u_1^\varepsilon,v_1^\varepsilon
}_{t_j};u_2^\varepsilon,
v_2^\varepsilon \bigr)+ (2L+1)\varepsilon\nonumber \\
&=& Y_{t_j}^{t_j,X^{t,x;u_1^\varepsilon,v_1^\varepsilon}_{t_j};
u_2^\varepsilon,v_2^\varepsilon}+
(2L+1)\varepsilon
\\
& = & Y_{t_j}^{t_j,X^{t,x;u^\varepsilon,v^\varepsilon}_{t_j};
u^\varepsilon,v^\varepsilon}+ (2L+1)\varepsilon\nonumber\\
& =& Y_{t_j}^{t,x;u^\varepsilon,v^\varepsilon}+
(2L+1)\varepsilon,\qquad P\mbox{-a.s.}
\nonumber
\end{eqnarray}
Indeed, the fact that
$X^{t,x;u^\varepsilon_1,v^\varepsilon_1}_{t_j}$ is ${\cal F}^\pi_{t_j-}
=\widetilde{\cal F}^\pi_{t_j}$-measurable, allows to substitute
this random variable at the place of $x'$ in the BSDE for $(Y^{t_j,x';
u^\varepsilon_1,v^\varepsilon_1}_s,Z^{t_j,x';u^\varepsilon_1,
v^\varepsilon_1}_s$, $ M^{t_j,x';u^\varepsilon_1,
v^\varepsilon_1}_s)_{s\in[t_j,T]}$. The uniqueness of the solution of the
resulting BSDE then yields $ Y^{t_j,X^{t,x;u^\varepsilon,
v^\varepsilon};u^\varepsilon,v^\varepsilon}_s=Y_s^{t,x;u^\varepsilon,
v^\varepsilon}, s\in[t_j,T]$.

Combining the above result with that of step 1, and taking into account
the monotonicity and the Lipschitz properties of the backward stochastic
semigroup, which are a direct consequence of the corresponding properties
of the solutions of BSDEs (the proof of them is similar to the
classical case (e.g., refer to Peng \cite{Pe1}), also refer to \cite
{CFS}) we obtain
%
%
%e3.16 #&#
%
\begin{eqnarray}
\label{(2.14)} %
\widetilde{W}_{t_j}^{\pi}(t,x) &
\le& G^{t,x;
\alpha^\varepsilon_1,\beta_1}_{t,t_j} \bigl(W^{\pi} \bigl(t_j,X_{t_j}^{t,
x;\alpha^\varepsilon_1,\beta_1}
\bigr) \bigr)+ \varepsilon
\nonumber
\\
& \le& G^{t,x;\alpha^\varepsilon_1,\beta_1}_{t,t_j} \bigl(Y_{t_j}^{t,x;u^\varepsilon,v^\varepsilon}+
(2L+1)\varepsilon \bigr) + \varepsilon
\nonumber\\
& \le& G^{t,x;u^\varepsilon_1,v^\varepsilon_1}_{t,t_j} \bigl(Y_{t_j}^{t,x;u^\varepsilon,v^\varepsilon}
\bigr)+ C\varepsilon
\nonumber
\\[-8pt]
\\[-8pt]
\nonumber
& = & G^{t,x;u^\varepsilon,v^\varepsilon}_{t,t_j} \bigl(Y_{t_j}^{t,x;
u^\varepsilon,v^\varepsilon}
\bigr)+ C\varepsilon
\nonumber\\
& = & Y_{t}^{t,x;u^\varepsilon,v^\varepsilon}+ C\varepsilon\nonumber
\\
& = & J^{\pi} \bigl(t,x; \alpha^\varepsilon,\beta \bigr)+ C
\varepsilon,\qquad P\mbox{-a.s., for all } \beta\in{\cal B}_{t,T}^\pi.
\nonumber
\end{eqnarray}

Therefore,
%
%
%e3.17 #&#
%
\begin{eqnarray}
\widetilde{W}_{t_j}^{\pi}(t,x)&\le&\esssup_{\alpha
\in
\mathcal{A}^\pi_{t,T}}
\essinf_{\beta\in
\mathcal{B}^\pi_{t,T}}J^{\pi}(t,x;\alpha,\beta)+C\varepsilon
\nonumber
\\[-8pt]
\\[-8pt]
\nonumber
&=&
W^{\pi}(t,x)+C\varepsilon, \qquad P\mbox{-a.s.},
\end{eqnarray}
and considering the
arbitrariness of the choice of $\varepsilon>0$ we can conclude the
proof.
\end{pf}

In order to complete the proof of the DPP, we need still the
following lemma.

%le3.4 #&#
\begin{lemma}\label{DPP-2} Under our standard assumptions it holds
that $ \widetilde{W}_{t_j}^{\pi}(t,x)\ge W^{\pi}(t,x)$, $P$-a.s.
\end{lemma}

\begin{pf} The proof of this lemma uses mainly arguments which
have been already developed in the frame of the proof of the preceding
lemma. For this reason, we give here rather a sketch than
a detailed proof.

Let us begin with fixing an arbitrary $\alpha\in{\cal A}_{t,T}^\pi$.
Given any $v_2\in{\cal V}_{t_j,T}^\pi$ we define $\alpha_1\in{\cal
A}_{t,t_j}^\pi$ by setting $\alpha_1(v_1):=\alpha(v_1\oplus
v_2)_{|[t,t_j]}\in{\cal U}_{t,t_j}^\pi$, for $v_1\in{\cal
V}_{t,t_j}^\pi$. Thanks to the nonanticipativity property of the
elements of ${\cal A}_{t,t_j}^\pi$, $\alpha_1$ does not depend on
the particular choice of $v_2$. From the definition of
$ \widetilde{W}_{t_j}^{\pi}(t,x)$, it follows that
\[
\widetilde{W}_{t_j}^{\pi}(t,x)\ge \essinf_{\beta_1\in\mathcal{B}^\pi_{t,t_j}}G^{t,x;\alpha_1,
\beta_1}_{t,t_j}
\bigl(W^{\pi} \bigl(t_j,X_{t_j}^{t,x;\alpha_1,\beta_1}
\bigr) \bigr),
\]
$P$-a.s., for all $\alpha_1\in\mathcal{A}^\pi_{t,t_j}$, and
from the argument developed in step 1 of the proof of Lemma \ref{DPP-1} we know
that, for an arbitrarily given $\varepsilon>0$ there
exists $\beta_1^\varepsilon\in\mathcal{B}^\pi_{t,t_j}$ (depending
on $\alpha_1\in{\cal A}_{t,t_j}^\pi$) such that
\[
\widetilde{W}_{t_j}^{\pi}(t,x)\ge G^{t,x;
\alpha_1,\beta_1^\varepsilon}_{t,t_j}
\bigl(W^{\pi} \bigl(t_j,X_{t_j}^{t,x;
\alpha_1,\beta_1^\varepsilon}
\bigr) \bigr)-\varepsilon,\qquad P\mbox{-a.s.}
\]

In analogy to step 2 of the proof of Lemma \ref{DPP-1}, we estimate
the expression
$W^{\pi}(t_j,X_{t_j}^{t,x;\alpha_1,\beta_1^\varepsilon})$ to which the
backward stochastic semigroup is applied in the above estimate. For
this, we let $(u_1^\varepsilon,v_1^\varepsilon)\in{\cal
U}_{t,t_j}^\pi\times{\cal V}_{t,t_j}^\pi$ be the unique control
couple associated with $(\alpha_1,\beta_1^\varepsilon)$ through
Lemma \ref{controls-NAD-strategies}, and we define
$\alpha_2^\varepsilon(v_2):=\alpha(v_1^\varepsilon\oplus
v_2)_{[t_j,T]}, v_2\in{\cal V}_{t_j,T}^\pi$. Such defined
mapping $\alpha_2^\varepsilon\dvtx {\cal V}_{t_j,T}^\pi\rightarrow{\cal
U}_{t_j,T}^\pi$ belongs to ${\cal A}_{t_j,T}^\pi$, and using an
adaptation of the argument with the Borel partition ${\cal O}_k, k\ge1$, of $R^d$, from step 2 of the proof of Lemma \ref{DPP-1},
which leads to (\ref{W_1}), we construct an NAD-strategy
$\beta_2^\varepsilon\in{\cal B}_{t_j,T}^\pi$ such that
%
%
%e3.18 #&#
%
\begin{eqnarray}
W^{\pi} \bigl(t_j,X_{t_j}^{t,x;\alpha_1,\beta_1^\varepsilon}
\bigr) & \ge& \essinf_{\beta_2\in
\mathcal{B}^\pi_{t_j,T}}J^{\pi} \bigl(t_j,X_{t_j}^{t,x;\alpha_1,
\beta_1^\varepsilon};
\alpha_2^\varepsilon,\beta_2 \bigr)
\nonumber
\\[-8pt]
\\[-8pt]
\nonumber
& \ge& J^{\pi} \bigl(t_j,X_{t_j}^{t,x;\alpha_1,\beta_1^\varepsilon};
\alpha_2^\varepsilon,\beta_2^\varepsilon \bigr)-
\varepsilon, \qquad P\mbox{-a.s.}
\end{eqnarray}
Letting $(u_2^\varepsilon,v_2^\varepsilon)\in{\cal
U}_{t_j, T}^\pi\times{\cal V}_{t_j,T}^\pi$ be the unique control
couple associated with $(\alpha_2^\varepsilon,\beta_2^\varepsilon)$
through Lemma \ref{controls-NAD-strategies}, we observe that, for
$\beta^\varepsilon\in{\cal B}_{t,T}^\pi$ defined by the relation
$\beta^\varepsilon(u):=\beta^\varepsilon_1(u_{|[t,t_j]})\oplus
\beta^\varepsilon_2(u_{|(t_j,T]}), u\in{\cal U}_{t,T}^\pi$, we
have the couple of controls $u^\varepsilon:=u^\varepsilon_1\oplus
u^\varepsilon_2\in{\cal U}_{t,T}^\pi$,
$v^\varepsilon:=v^\varepsilon_1\oplus v^\varepsilon_2\in{\cal
V}_{t,T}^\pi$ associated with $(\alpha,\beta^\varepsilon)$ through
Lemma~\ref{controls-NAD-strategies}:
\begin{eqnarray*}
\alpha\bigl(v^\varepsilon\bigr)&=&\alpha\bigl(v_1^\varepsilon
\oplus v_2^\varepsilon\bigr) =\alpha_1
\bigl(v_1^\varepsilon\bigr)\oplus\alpha_2^\varepsilon
\bigl(v_2^\varepsilon\bigr) =u_1^\varepsilon
\oplus u_2^\varepsilon=u^\varepsilon,
\\
\beta^\varepsilon\bigl(u^\varepsilon\bigr)&=&\beta^\varepsilon
_1\bigl(u_1^\varepsilon\bigr) \oplus
\beta^\varepsilon_2\bigl(u_2^\varepsilon
\bigr)=v_1^\varepsilon\oplus v_2^\varepsilon=v^\varepsilon.
\end{eqnarray*}

Consequently, thanks to the monotonicity and Lipschitz
properties of the backward stochastic semigroup, we have
%
%
%e3.19 #&#
%
\begin{eqnarray}
\label{(2.17)} %
\widetilde{W}_{t_j}^{\pi}(t,x) &
\ge& G^{t,x;\alpha_1,\beta_1^\varepsilon}_{t,t_j} \bigl(W^{\pi} \bigl(t_j,X_{t_j}^{t,x;\alpha_1,\beta_1^\varepsilon}
\bigr) \bigr)-\varepsilon
\nonumber
\\[-2pt]
& \ge& G^{t,x;\alpha_1,\beta_1^\varepsilon}_{t,t_j} \bigl(J^{\pi}
\bigl(t_j, X_{t_j}^{t,x;\alpha_1,\beta_1^\varepsilon}; \alpha_2^\varepsilon,
\beta_2^\varepsilon \bigr)-\varepsilon \bigr)-\varepsilon
\\[-2pt]
& \ge& G^{t,x;u_1^\varepsilon,v_1^\varepsilon}_{t,t_j} \bigl(Y_{t_j}^{t_j,
X_{t_j}^{t,x;u_1^\varepsilon,v_1^\varepsilon}; u_2^\varepsilon,
v_2^\varepsilon}
\bigr)-C\varepsilon
\nonumber
\\[-2pt]
& = & G^{t,x;u^\varepsilon,v^\varepsilon}_{t,t_j} \bigl(Y_{t_j}^{t_j,
X_{t_j}^{t,x;u^\varepsilon,v^\varepsilon}; u^\varepsilon,
v^\varepsilon}
\bigr)-C\varepsilon
\nonumber\\[-2pt]
& = & G^{t,x;u^\varepsilon,v^\varepsilon}_{t,t_j} \bigl(Y_{t_j}^{t,x;
u^\varepsilon, v^\varepsilon}
\bigr)-C\varepsilon
\nonumber
\\[-2pt]
& = & Y_{t}^{t,x; u^\varepsilon, v^\varepsilon}-C\varepsilon
\nonumber
\\[-2pt]
& = & Y_{t}^{t,x; \alpha,\beta^\varepsilon}-C\varepsilon,\qquad P\mbox{-a.s.}
\nonumber
\end{eqnarray}

We take in the latter estimate first the essential infimum
over $\beta\in{\cal B}_{t,T}^\pi$, and then the essential supremum
over all $\alpha\in{\cal A}_{t,T}^\pi$. Thus, by considering the
arbitrariness of $\varepsilon>0$, we get the statement of the
lemma.\vspace*{-2pt}
\end{pf}

As a consequence of the proof of the DPP, we get the following proposition.\vspace*{-2pt}

%pr3.1 #&#
\begin{proposition}\label{proposition3.1} Under our standard
assumptions, for all $(t,x)\in[0,T]
\times R^d$, it holds
%
%
%e3.20 #&#
%
\begin{eqnarray}
\label{DPP-a} %
W^{\pi}(t,x) &=& \sup_{\alpha\in
\mathcal{A}^\pi_{t,T}}
\inf_{\beta\in
\mathcal{B}^\pi_{t,T}}E \bigl[J^{\pi}(t,x; \alpha,\beta) \bigr],
\nonumber
\\[-8pt]
\\[-8pt]
\nonumber
U^{\pi}(t,x) &=& \inf_{\beta\in\mathcal{B}^\pi_{t,T}}\sup
_{\alpha
\in\mathcal{A}^\pi_{t,T}}E \bigl[J^{\pi}(t,x;\alpha,\beta) \bigr].\vspace*{-2pt}
\end{eqnarray}
\end{proposition}

By combining the above lemma with Remark \ref{classical}, we get the following
result under the classical assumption of a running payoff function not depending
on $(y,z)$:\vspace*{-2pt}

%co3.1 #&#
\begin{corollary} Let us suppose in addition to our standard
assumptions that
the coefficient $f(s,x,y,z,u,v)$ does not depend on $(y,z)$. Then, for all
$(t,x)\in[0,T]\times R^d$,
%
%
%e3.21 #&#
%
\begin{eqnarray}
\label{DPP-aaa} %
 W^{\pi}(t,x) &=& \sup
_{\alpha\in
\mathcal{A}^\pi_{t,T}} \inf_{\beta\in
\mathcal{B}^\pi_{t,T}}E \biggl[\Phi
\bigl(X^{t,x;u,
v}_T \bigr)\nonumber\\[-3pt]
&&\hspace*{98pt}{}+\int_t^Tf
\bigl(s,X^{t,x;u,v}_s,u_s,v_s \bigr)
\,ds \biggr],
\nonumber
\\[-9pt]
\\[-9pt]
\nonumber
U^{\pi}(t,x) &=& \inf_{\beta\in\mathcal{B}^\pi_{t,T}}\sup_{\alpha
\in\mathcal{A}^\pi_{t,T}}E
\biggl[\Phi \bigl(X^{t,x;u,
v}_T \bigr)\\[-3pt]
&&\hspace*{98pt}{}+\int_t^Tf
\bigl(s,X^{t,x;u,v}_s,u_s,v_s \bigr)
\,ds \biggr].\nonumber\vspace*{-2pt}
\end{eqnarray}
\end{corollary}

Now we prove the above Proposition \ref{proposition3.1}.\vspace*{-2pt}

\begin{pf} Let\vspace*{1pt} $(t,x)\in[0,T)\times R^d$, and $t_j\in\pi$ be such that
$t_j\le t<t_{j+1}$. As we have shown in the proof of
the DPP that $\widetilde{W}^\pi_{t_j}(t,x)$ and $W^\pi(t,x)$ coincide,
we see
from (\ref{(2.14)}) that, for every $\varepsilon>0$, there exists
$\alpha^\varepsilon\in{\cal A}_{t,T}^\pi$ such that, for all $\beta
\in
{\cal B}_{t,T}^\pi$,
\[
W^\pi(t,x)\le J^{\pi} \bigl(t,x;\alpha^\varepsilon,
\beta \bigr)+\varepsilon,\qquad P\mbox{-a.s.}
\]
Consequently, taking into account that $W^\pi(t,x)$ is deterministic,
we get $W^\pi(t, x)\le E[J^{\pi}(t,x;\alpha^\varepsilon,\beta
)]+\varepsilon$. By taking
first the infimum over all $\beta\in{\cal B}_{t,T}^\pi$ and after
the supremum
over $\alpha\in{\cal A}_{t,T}^\pi$, we obtain
\[
W^{\pi}(t,x)\le\sup_{\alpha\in\mathcal{A}^\pi_{t,T}} \inf_{\beta\in\mathcal{B}^\pi_{t,T}}E
\bigl[J^{\pi}(t,x;\alpha,\beta) \bigr].
\]
To get the converse relation, we observe that, due to (\ref{(2.17)}),
for every $\varepsilon>0$ and all $\alpha\in{\cal A}_{t,T}^\pi$, there
exists some
$\beta^\varepsilon\in{\cal B}_{t,T}^\pi$ such that
\[
W^\pi(t,x)\ge J^{\pi}\bigl(t,x;\alpha,\beta^\varepsilon
\bigr)-\varepsilon,\qquad P\mbox{-a.s.}
\]
By taking the expectation on both sides of this inequality,
after the
infimum with respect to $\beta^\varepsilon\in{\cal B}_{t,T}^\pi$ and,
at the end, the
supremum over $\alpha\in{\cal A}_{t,T}^\pi$, we obtain that
\[
W^{\pi}(t,x)\ge\sup_{\alpha\in\mathcal{A}^\pi_{t,T}} \inf_{\beta\in\mathcal{B}^\pi_{t,T}}E
\bigl[J^{\pi}(t,x;\alpha,\beta) \bigr].
\]
This proves the statement for $W^{\pi}(t,x)$; that for
$U^{\pi
}(t,x)$ can be proved
similarly.
\end{pf}

At the end of this section, let us still consider the H\"{o}lder
continuity of the lower and the upper value functions along the
partition with respect to the time.

%pr3.2 #&#
\begin{proposition}\label{Hoelder}
Under our standard assumptions there exists a constant~$C$ which is
independent of
the underlying partition $\pi$ of the interval $[0,T]$, such that
%
%
%e3.22 #&#
%
\begin{eqnarray}
\bigl|W^\pi(t,x)-W^\pi(s,x) \bigr|+ \bigl|U^\pi(t,x)-U^\pi(s,x)
\bigr| \le C|t-s|^{1/2},
\nonumber
\\[-8pt]
\\[-8pt]
 \eqntext{s,t\in[0,T], x\in R^d.}
\end{eqnarray}
\end{proposition}

\begin{pf} We restrict ourselves to the proof for $W^\pi$; that for
$U^\pi$ is analogous.\vadjust{\goodbreak}

\textit{Step} 1. Given a partition $\pi$ of the interval $[0,T]$,
let us suppose that $0\le t< s\le T$ and fix arbitrarily $\varepsilon>0$.
From the proof of Proposition \ref{proposition3.1}, we know
that there exists $\alpha^\varepsilon\in{\cal A}_{t,T}^\pi$ such that,
for all
$\beta\in{\cal B}_{t,T}^\pi$,
%
%
%e3.23 #&#
%
\begin{equation}
\label{a1}W^\pi(t,x)\le E \bigl[J^{\pi} \bigl(t,x;\alpha
^\varepsilon,\beta \bigr) \bigr]+\varepsilon.
\end{equation}
For any fixed $v^0\in V$ we let $v^0_1:=v^0I_{[t,s)}$. Then, for
$v_2\in
{\cal V}_{s,
T}^\pi$, $v^0_1\oplus v_2:=v^0I_{[t,s)}+v_2I_{[s,T]}\in{\cal
V}_{t,T}^\pi$, and
$\widetilde{\alpha}^\varepsilon(v_2):=\alpha^\varepsilon
(v^0_1\oplus
v_2)_{|[s,T]}
\in{\cal U}_{s,T}^\pi$. Moreover, it can be easily checked that such defined
mapping $\widetilde{\alpha}^\varepsilon$ belongs to ${\cal
A}_{s,T}^\pi
$. Again from
the proof of Proposition \ref{proposition3.1}, it follows that there is
$\widetilde{\beta}^\varepsilon\in{\cal B}_{s,T}^\pi$ such that
%
%
%e3.24 #&#
%
\begin{equation}
\label{a2}W^\pi(s,x)\ge E \bigl[J^{\pi} \bigl(s,x;
\widetilde { \alpha }^\varepsilon, \widetilde{\beta}^\varepsilon \bigr)
\bigr]- \varepsilon.
\end{equation}

Let $(u_2^\varepsilon,v_2^\varepsilon)\in{\cal U}^\pi
_{s,T}\times
{\cal V}^\pi_{s,T}$ be associated with $(\widetilde{\alpha
}^\varepsilon,
\widetilde{\beta}^\varepsilon)$ through Lemma \ref{controls-NAD-strategies}:
$\widetilde{\alpha}^\varepsilon(v_2^\varepsilon)=u_2^\varepsilon,
\widetilde{\beta}^\varepsilon(u_2^\varepsilon)=v_2^\varepsilon$,
$ds\,dP$-a.e.
on $[s,T]\times\Omega$.

On the other hand, let us define $\beta^\varepsilon(u):=
v_1^0\oplus\widetilde{\beta}^\varepsilon(u_{|[s,T]}), u\in{\cal U}_{t,
T}^\pi$. Obviously, $\beta^\varepsilon\in{\cal B}_{t,T}^\pi$.
Putting $u^\varepsilon:=\alpha^\varepsilon(v_1^0\oplus
v_2^\varepsilon
)\in
{\cal U}_{t,T}^\pi$, we deduce from the fact $u^\varepsilon_{|[s,T]}
=\alpha^\varepsilon(v_1^0\oplus v_2^\varepsilon)_{|[s,T]}=
\widetilde{\alpha}^\varepsilon(v_2^\varepsilon)=u_2^\varepsilon$,
that $(u^\varepsilon,v^\varepsilon:=v_1^0\oplus v_2^\varepsilon)
\in{\cal U}_{t,T}^\pi\times{\cal V}_{t,T}^\pi$ satisfies
%
%
%e3.25 #&#
%
\begin{eqnarray}
\alpha^\varepsilon \bigl(v^\varepsilon \bigr)&=&u^\varepsilon
\quad\mbox{and }
\nonumber
\\[-8pt]
\\[-8pt]
\nonumber
\beta^\varepsilon \bigl(u^\varepsilon \bigr)&=&v_1^0
\oplus\widetilde{\beta }^\varepsilon \bigl(u_2^\varepsilon
\bigr)=v_1^0\oplus v_2^\varepsilon=v^\varepsilon,
\end{eqnarray}
over the interval $[t,T]$, while over the smaller interval
$[s,T]$ it holds
%
%
%e3.26 #&#
%
\begin{eqnarray}
\widetilde{\alpha}^\varepsilon \bigl(v^\varepsilon_{|[s,T]}
\bigr) &=&\widetilde{\alpha}^\varepsilon \bigl(v_2^\varepsilon
\bigr)=u_2^\varepsilon= u^\varepsilon_{|[s,T]}\quad
\mbox{and }
\nonumber
\\[-8pt]
\\[-8pt]
\nonumber
\widetilde{\beta}^\varepsilon \bigl(u^\varepsilon_{|[s,T]}
\bigr) &=&\widetilde{\beta}^\varepsilon \bigl(u^\varepsilon_2
\bigr)=v^\varepsilon_2= v^\varepsilon_{|[s,T]}.
\end{eqnarray}
Consequently, from the relation (\ref{a1}) and (\ref{a2}) it
follows that
%
%
%e3.27 #&#
%
\begin{eqnarray}
\label{b1} %
W^\pi(t,x) &\le& E \bigl[J^{\pi}
\bigl(t,x;u^\varepsilon,v^\varepsilon \bigr) \bigr]+ \varepsilon,
\nonumber
\\[-8pt]
\\[-8pt]
\nonumber
W^\pi(s,x) &\ge& E \bigl[J^{\pi}
\bigl(s,x;u^\varepsilon,v^\varepsilon \bigr) \bigr]-\varepsilon,
\end{eqnarray}
from where
%
%
%e3.28 #&#
%
\begin{eqnarray}
\label{b2} %
& & W^\pi(t,x)-W^\pi(s,x)
\nonumber
\\
&& \qquad\le E \bigl[J^{\pi} \bigl(t,x;u^\varepsilon,v^\varepsilon
\bigr) -J^{\pi} \bigl(s,x;u^\varepsilon,v^\varepsilon \bigr)
\bigr]+2 \varepsilon
\\
&&\qquad \le E \bigl[ \bigl|Y_s^{t,x;u^\varepsilon,v^\varepsilon}- Y_s^{s,x;u^\varepsilon,v^\varepsilon}
\bigr| \bigr]+ \bigl|E \bigl[Y_t^{t,x;
u^\varepsilon,v^\varepsilon}-Y_s^{t,x;u^\varepsilon,v^\varepsilon}
\bigr] \bigr| +2\varepsilon.
\nonumber
\end{eqnarray}
We emphasize that, if $s\notin\pi$, unlike the classical
Markovian case we do not have here that $Y_s^{t,x;u^\varepsilon,
v^\varepsilon}=Y_s^{s,X_s^{t,x;u^\varepsilon,
v^\varepsilon};u^\varepsilon,v^\varepsilon}=J^{\pi
}(s,X_s^{t,x;u^\varepsilon,
v^\varepsilon};u^\varepsilon,v^\varepsilon)$. Indeed, here, if $s\in
(t_{j-1},t_j)$, then $X_s^{t,x;u^\varepsilon,
v^\varepsilon}$ is ${\cal F}^\pi_{s-}$-measurable, where
${\cal F}^\pi_{s-}={\cal F}^B_s\vee{\cal H}_j\supsetneqq{\cal F}^B_s
\vee{\cal H}_{j-1}=\widetilde{\cal F}^\pi_{s-}$, where the BSDE is
considered with respect to the filtration $\widetilde{\mathbb{F}}^\pi$.
However, from the both BSDEs
%
%
%e3.29 #&#
%
\begin{equation}
\cases{ %
dY^{t,x; u, v}_r = -E \bigl[f
\bigl(r,X^{t,x; u, v}_r, Y^{t,
x; u, v}_r,
Z^{t,x; u, v}_r,u_r, v_r \bigr) |
\widetilde{\cal F}^\pi _r \bigr]\,dr
\vspace*{2pt}\cr
\hspace*{52pt}{}+Z^{t,x;u, v}_r \,dB_r+dM^{t,x; u, v}_r,
\vspace*{2pt}
\cr
Y^{t,x; u, v}_T = E \bigl[\Phi
\bigl(X^{t,x; u, v}_T \bigr) | \widetilde{\cal F}^\pi_{T-}
\bigr] }
\end{equation}
and
%
%
%e3.30 #&#
%
\begin{equation}
\cases{ %
dY^{s,x; u, v}_r = -E \bigl[f
\bigl(r,X^{s,x; u, v}_r, Y^{s,
x; u, v}_r,
Z^{s,x; u, v}_r, u_r, v_r \bigr) |
\widetilde{\cal F}^\pi _r \bigr]\,dr\vspace*{2pt}\cr
\hspace*{52pt}{}+Z^{s,x;u, v}_r \,dB_r+dM^{s,x; u, v}_r,
\vspace*{2pt}
\cr
Y^{s,x; u, v}_T = E \bigl[\Phi
\bigl(X^{s,x; u, v}_T \bigr) | \widetilde{\cal
F}^\pi_{T-} \bigr],}
\end{equation}
both studied over the time interval $[s,T]$, we deduce with
standard BSDE estimates that (or, refer to \cite{CFS})
%
%
%e3.31 #&#
%
\begin{eqnarray}
\label{aa1} %
&&E \bigl[ \bigl|Y_s^{t,x;u^\varepsilon,v^\varepsilon}-
Y_s^{s,x;u^\varepsilon,v^\varepsilon} \bigr|^2 \bigr]
\nonumber
\\
&&\qquad \le CE \Bigl[\sup_{r\in[s,T]} \bigl|X^{t,x; u^\varepsilon, v^\varepsilon}_r-X^{s,x; u^\varepsilon,
v^\varepsilon}_r
\bigr|^2 \Bigr]
\\
&&\qquad \le CE \bigl[ \bigl|X^{t,x; u^\varepsilon, v^\varepsilon}_s-x \bigr|^2
\bigr]\le C(s-t)
\nonumber
\end{eqnarray}
(Recall that the coefficients $\sigma$ and $b$ are bounded
and Lipschitz).
Thus, from BSDE (\ref{BSDE}), the boundedness of $f(s,x,y,0,u,v)$,
the Lipschitz continuity of $f(s,x,y,z,u,v)$ in $z$ as well as
(\ref{BSDE-estimates-1}),
%
%
%e3.32 #&#
%
\begin{eqnarray}
\label{b3} %
& & W^\pi(t,x)-W^\pi(s,x)
\nonumber
\\
&&\qquad \le E \bigl[ \bigl|Y_s^{t,x;u^\varepsilon,v^\varepsilon}- Y_s^{s,x;u^\varepsilon,v^\varepsilon}
\bigr| \bigr]+ \bigl|E \bigl[Y_t^{t,x;
u^\varepsilon,v^\varepsilon}-Y_s^{t,x;u^\varepsilon,v^\varepsilon}
\bigr] \bigr| +2\varepsilon
\nonumber
\\
&&\qquad \le C(s-t)^{1/2}+ 2\varepsilon
\nonumber
\\[-8pt]
\\[-8pt]
\nonumber
& &\qquad\quad{}+ (s-t)^{1/2} \biggl(E \biggl[\int_t^s
\bigl|f \bigl(r,X_r^{t,x;u^\varepsilon
,v^\varepsilon}, Y_r^{t,x;u^\varepsilon,v^\varepsilon},Z_r^{t,x;u^\varepsilon,
v^\varepsilon},\\
&&\hspace*{253pt}{}u^\varepsilon_r,v^\varepsilon_r
\bigr) \bigr|^2\,dr \biggr] \biggr)^{1/2}
\nonumber
\\
&&\qquad \le C|s-t|^{1/2}+2\varepsilon,
\nonumber
\end{eqnarray}
for some constant $C$ not depending on $\pi$ and on
$\varepsilon$. Thus, in virtue of the arbitrariness of $\varepsilon>0$
we have
\[
W^\pi(t,x)-W^\pi(s,x)\le C|s-t|^{1/2}.
\]

\textit{Step} 2. Now, for the same partition
$\pi$, and the case $0\le t< s\le T$, we make a lower
estimate for $W^\pi(t,x)-W^\pi(s,x)$. For this we notice that, for
arbitrarily given $\varepsilon>0$ we can find
$\widetilde{\alpha}^\varepsilon\in{\cal A}_{s,T}^\pi$ such that,
for all
$\widetilde{\beta}\in{\cal B}_{s,T}^\pi$,
%
%
%e3.33 #&#
%
\begin{equation}
\label{c1}W^\pi(s,x)\le E \bigl[J^\pi \bigl(s,x;
\widetilde {\alpha }^\varepsilon, \widetilde{\beta} \bigr) \bigr]+\varepsilon.
\end{equation}
For any fixed $u^0\in U$ we put $u^0_1:=u^0I_{[t,s)}$, and we define
$\alpha^\varepsilon\in{\cal A}_{t,T}^\pi$ by setting
$\alpha^\varepsilon(v):=u_1^0\oplus\widetilde{\alpha}^\varepsilon
(v_{|[s,T]}), v\in{\cal V}_{t,T}^\pi$. Let $\beta^\varepsilon\in
{\cal B}_{t,T}^\pi$ such that
%
%
%e3.34 #&#
%
\begin{equation}
\label{c2}W^\pi(t,x)\ge E \bigl[J^\pi \bigl(t,x;\alpha
^\varepsilon, \beta^\varepsilon \bigr) \bigr]-\varepsilon,
\end{equation}
and let $(u^\varepsilon,v^\varepsilon)\in{\cal
U}_{t,T}^\pi
\times{\cal V}_{t,T}^\pi$ be associated with $(\alpha^\varepsilon,
\beta^\varepsilon)$ through Lemma \ref{controls-NAD-strategies}. On
the other hand, by defining $\widetilde{\beta}^\varepsilon\in
{\cal B}_{s,T}^\pi$ by putting $\widetilde{\beta}^\varepsilon(u_2)=
\beta^\varepsilon(u_1^0\oplus u_2)_{|[s,T]}, u_2\in{\cal
U}_{s,T}^\pi
$, it can be easily verified
that $(u^\varepsilon_{|[s,T]},v^\varepsilon_{|[s,T]})\in{\cal
U}_{s,T}^\pi
\times{\cal V}_{s,T}^\pi$ is associated with
$(\widetilde{\alpha}^\varepsilon,\widetilde{\beta}^\varepsilon)$
in the
sense of Lemma \ref{controls-NAD-strategies}. Consequently,
%
%
%e3.35 #&#
%
\begin{eqnarray}
\label{b1} %
W^\pi(s,x) &\le& E \bigl[J^\pi
\bigl(s,x;u^\varepsilon,v^\varepsilon \bigr) \bigr]+ \varepsilon,
\nonumber
\\[-8pt]
\\[-8pt]
\nonumber
W^\pi(t,x) &\ge& E \bigl[J^\pi
\bigl(t,x;u^\varepsilon,v^\varepsilon \bigr) \bigr]-\varepsilon,
\end{eqnarray}
and we can proceed now in analogy to step 1 to
deduce that
\[
W^\pi(t,x)-W^\pi(s,x)\ge-C|s-t|^{1/2}.
\]
Combining this result with that of step 1 we complete the proof.
\end{pf}

%s4 #&#
\section{Value in mixed strategies and associated HJB--Isaacs equation}\label{sec4}

The objective of this section is to study the limit of the lower and
the upper value functions $W^\pi$ and $U^\pi$ along a partition
$\pi$, when the mesh of the partition $\pi$ tends to zero, and to
show that both $W^\pi$ and $U^\pi$ converge uniformly on compacts to
the same limit function $V$ which is the unique viscosity solution
of the following Hamilton--Jacobi--Bellman--Isaac equation
%
%
%e4.1 #&#
%
\begin{equation}
\label{HJBI}\quad\hspace*{8pt}\cases{ %
\displaystyle\frac{\partial}{\partial t}V(t,x)+H \bigl(t,x,
\bigl(V,DV,D^2V \bigr) (t,x) \bigr)=0&\quad $(t,x)\in[0,T)\times
R^d$, \vspace*{2pt}
\cr
V(T,x)=\Phi(x),&\quad $x\in R^d$,
}\hspace*{-6pt}
\end{equation}
with Hamiltonian
%
%
%e4.2 #&#
%
\begin{eqnarray}
&&H(t,x,y,p,A) \nonumber\\
&&\qquad= \sup_{\mu\in{\cal P}(U)} \inf
_{\nu\in{\cal P}(V)}
\nonumber
\\[-8pt]
\\[-8pt]
\nonumber
&&\quad\qquad{}\times\int_{U\times V} \biggl(
\frac{1}{2}\operatorname{tr} \bigl(\sigma\sigma^{T}(t,x, u, v)A \bigr)+ b(t, x, u, v)p
\\
&&\hspace*{102pt}{}+ f \bigl(t, x,y,p\cdot\sigma(t,x,u, v),u, v \bigr)
\biggr)\mu\otimes \nu(du\,dv),\nonumber
\end{eqnarray}
$(t,x,y,p,A)\in[0,T]\times R^d\times R\times R^d\times S^d$,
where $S^d$ denotes the space of symmetric matrices from
$R^{d\times d}$. For this we need the following supplementary assumption
which is coherent with our standard assumptions on the coefficients
$\sigma,b$ and~$f$.

%co4.1 #&#
\begin{condition}\label{condition} We suppose that either
\begin{itemize}
\item $\sigma(s,x,u,v)=\sigma(s,x), (s,x,u,v)
\in[0,T]\times R^d\times U\times V$ is independent of the controls;
or
\item $f(s,x,y,z,u,v)$ is linear in $z$:
\[
f(s,x,y,z,u,v)=f_0(s,x,y,u,v)+f_1(s)z,
\]
$(s,x,y,z,u,v)\in[0,T]\times R^d\times R\times R^d\times
U\times V$, where $f_0=(f_0(s,x,y,u,v))\dvtx\break   [0,T]\times R^d\times R\times
U\times V \rightarrow R$ bounded, jointly continuous and
Lipschitz in $(x,y)$, uniformly with respect to $(s,u,v)$, and
$f_1\dvtx [0,T]\rightarrow R^d$ is continuous.\vspace*{-2pt}
\end{itemize}
\end{condition}

More precisely, we have the following
theorem.\vspace*{-2pt}
%
%th4.1 #&#
\begin{theorem}\label{main result} Under our standard assumptions on the
coefficients $\sigma,b,f$ and $\Phi$ as well as Condition \ref{condition},
we have the existence of a bounded, continuous function $V\dvtx [0,T]\times
R^d\rightarrow R$ such that, for every sequence of partitions $\pi_n$,
$n\ge1$, of the interval $[0,T]$ with mesh $|\pi_n|\rightarrow0$, as
$n\rightarrow+\infty$, $W^{\pi_n}\rightarrow V$, and $U^{\pi_n}
\rightarrow V$, uniformly on compacts, as $n\rightarrow+\infty$.
Moreover, $V$ is the viscosity solution of PDE~(\ref{HJBI}), unique in the
class of continuous functions with polynomial growth.\vspace*{-2pt}
\end{theorem}

For the convenience of the reader, we recall briefly the definition of
a viscosity solution, which we give directly for PDE (\ref{HJBI}). The reader
interested in a more detailed description of the concept of viscosity solution
is referred to the overview paper by Crandall, Ishii and Lions \cite{CIL}.\vspace*{-2pt}

%de4.1 #&#
\begin{definition}\label{de4.1} A function $V\in C([0,T]\times R^d)$ is said to be:

(i) a viscosity subsolution of PDE (\ref{HJBI}), if, first,
$V(T,x)\le
\Phi(x), x\in R^d$, and if, second, for any $(t,x)\in[0,T)\times R^d$
and any test function $\varphi\in C^{1,2}([0,T]\times R^d)$ such that
$V-\varphi$ achieves a local maximum at $(t,x)$, it holds
%
%
%e4.3 #&#
%
\begin{equation}
\frac{\partial}{\partial t}\varphi (t,x)+H \bigl(t,x, \bigl(\varphi, \nabla
\varphi,D^2\varphi \bigr) (t,x) \bigr)\ge0;
\end{equation}
(ii) a viscosity supersolution of PDE (\ref{HJBI}), if, first,
$V(T,x)\ge
\Phi(x), x\in R^d$, and if, second, for any $(t,x)\in[0,T)\times R^d$
and any test function $\varphi\in C^{1,2}([0,T]\times R^d)$ such that
$V-\varphi$ achieves a local minimum at $(t,x)$, it holds
%
%
%e4.4 #&#
%
\begin{equation}
\frac{\partial}{\partial t}\varphi (t,x)+H \bigl(t,x, \bigl(\varphi, \nabla
\varphi,D^2\varphi \bigr) (t,x) \bigr)\le0;
\end{equation}
(iii) a viscosity solution of (\ref{HJBI}) if it is both a viscosity
sub- but also a viscosity supersolution of (\ref{HJBI}).\vspace*{-2pt}
\end{definition}

%re4.1 #&#
\begin{remark}Let us point out that in Definition \ref{de4.1} the
space\break
$C^{1,2}([0,T]\times R^d)$ of the test functions can be replaced by any
subspace containing $C^\infty([0,T]\times R^d)$,\vadjust{\goodbreak} as long as one can
show the uniqueness with the help of $C^\infty$-test functions, as, for
instance, done in \cite{CIL}. Thus, our uniqueness results allows to
restrict to a class of test functions, more adapted for our
computations, the space $C^{3}([0,T]\times R^d)$ of functions which are
three times continuous differentiable with respect to $(t,x)$. On the
other hand, taking into account the uniform boundedness of the
functions $W^\pi, U^\pi$ and, hence, also of $V$, the standard
argument of changing a test function $\varphi\in C^{3}([0,T]\times
R^d)$ such that $V-\varphi$ achieves a local extremum at $(t,x)$, at
the exterior of a small ball around $(t,x)$, allows to consider only
test functions $\varphi\in C^{3}_{\ell,b}([0,T]\times R^d)$, that is,
$C^{3}$-functions with bounded derivatives of orders 1, 2 and 3 (and
which themselves have, consequently, a linear growth).
\end{remark}

Following the arguments developed, for example, in Str\"{o}mberg \cite
{ST} Theorem~5, we have the following comparison principle.

%pr4.1 #&#
\begin{proposition}\label{comparison principle} Let us suppose our
standard assumptions on the coefficients $\sigma,b,f$ and $\Phi$, and
let $V_1, V_2\dvtx [0,T]\times R^d\rightarrow R$ be continuous functions having a growth
not exceeding that of $\exp\{\gamma|x|\}$, for some $\gamma>0$.
Then, if $V_1$ is a viscosity subsolution and $V_2$ a viscosity
supersolution of (\ref{HJBI}), we have $V_1(t,x)\le V_2(t,x), (t,x)
\in[0,T]\times R^d$.
\end{proposition}

%re4.2 #&#
\begin{remark}Let us emphasize that the condition of exponential growth is
optimal for the uniqueness of the continuous viscosity solution, as
long as $\sigma$ is bounded; this is the case due to our assumptions.
However, the assumption of bounded coefficients and so, in particular,
that of $\sigma$, has been imposed in order to simplify our argument.
Our approach can be extended without major difficulties to coefficients
$\sigma$ of linear growth. In this case the class of continuous
functions $V$ within which one has the uniqueness of the viscosity
solution is smaller than that of the above Proposition \ref{comparison principle}; it's that
of $V$ such that, for some $\gamma>0$,
%
%
%e4.5 #&#
%
\begin{eqnarray}
\lim_{|x|\rightarrow+\infty}V(t,x)\exp \bigl\{-\gamma \bigl(\log \bigl(|x|+1 \bigr)
\bigr)^2 \bigr\} =0
\nonumber
\\[-8pt]
\\[-8pt]
\eqntext{\mbox{uniformly with respect to }t\in[0,T];}
\end{eqnarray}
see, for example, \cite{Buckdahn-Li-2008}.
\end{remark}

As a direct consequence of this comparison principle, we have
the following corollary.
%
%co4.1 #&#
\begin{corollary} PDE (\ref{HJBI}) has at most one continuous viscosity
solution $V\dvtx [0,T]\times R^d\rightarrow R$ with exponential growth,
that is, satisfying the condition that, for suitable $\gamma>0$,
%
%
%e4.6 #&#
%
\begin{eqnarray}
\lim_{|x|\rightarrow+\infty}V(t,x)\exp\{-\gamma |x|\} =0
\nonumber
\\[-8pt]
\\[-8pt]
\eqntext{\mbox{uniformly with respect to }t\in[0,T].}
\end{eqnarray}
In particular, uniqueness holds within the class of continuous functions
with polynomial growth.\vspace*{-2pt}
\end{corollary}

All what follows will be devoted to the proof of Theorem \ref{main result}.
The proof will be given through a sequel of auxiliary results.

Let us begin by choosing an arbitrary sequence of partitions $\pi_n:=
\{0=t_0^n<t_1^n<\cdots<t_{N_n}^n=T\}$, $n\ge1$, of the interval $[0,T]$
such that $|\pi_n|:=\sup_{1\le i\le N_n}(t_i-t_{i-1})\rightarrow0$,
as $n
\rightarrow+\infty$. Then, from Lemma \ref{Hoelder} and Proposition
\ref{Hoelder}, we see that the family of functions $(W^{\pi_n},U^{\pi_n}),
n\ge1$, is uniformly Lipschitz in $x$, uniformly with respect to $t$,
and H\"{o}lder continuous in $t$, uniformly with respect to $x$. Consequently,
the following result follows from the Arzel\`{a}--Ascoli theorem
combined with
a standard diagonalization argument.\vspace*{-2pt}

%le4.1 #&#
\begin{lemma}\label{limit_W_pi_n} There exists a subsequence of partitions,
which we denote again
by $(\pi_n)_{n\ge1}$, as well as bounded continuous functions $W, U\dvtx [0,T]
\times R^d\rightarrow R$ such that $(W^{\pi_n},U^{\pi_n})\rightarrow(W,U)$,
uniformly on compacts in $[0,T]\times R^d$. Moreover,
%
%
%e4.7 #&#
%
\begin{equation}
\label{estimate_W,U}\qquad \bigl|W(t,x)-W \bigl(t',x'
\bigr) \bigr|+ \bigl|U(t,x)-U \bigl(t',x' \bigr) \bigr|\le C \bigl(
\bigl|t-t' \bigr|^{1/2}+ \bigl|x-x'\bigr | \bigr),
\end{equation}
$(t,x), (t',x')\in[0,T]\times R^d$, where $C$ is a constant
which does not depend on the choice of the
sequence of partitions $\pi_n, n\ge1$.\vspace*{-2pt}
\end{lemma}

Although the functions $W,U$ given by the above lemma depend a priori
on the
choice of the sequence of partitions $\pi_n,n\ge1$, as well as on the
subsequence with respect to which $(W^{\pi_n},U^{\pi_n})$ converges,
we will
show later that $W,U$ are universal and coincide even.

Inspired by the approach in \cite{Buckdahn-Li-2008} we put,
for some arbitrarily chosen but fixed $\varphi\in C^3_{\ell,b} ([0,T]
\times
{\mathbb{R}}^d)$,
%
%
%e4.8 #&#
%
\begin{eqnarray}
\qquad F(s,x,y,z,u, v)&=&f \bigl(s, x, y+\varphi(s,x), z+
D\varphi(s,x)\cdot
\sigma (s,x,u, v),u, v \bigr)
\nonumber
\\[-9pt]
\\[-9pt]
\nonumber
& &{}+ {\cal L}\varphi(s,x,u,v),
\end{eqnarray}
$(s,x,y,z,u, v)\in[0,T] \times{\mathbb{R}}^d \times
{\mathbb{R}} \times{\mathbb{R}}^d \times U \times V$, where
%
%
%e4.9 #&#
%
\begin{eqnarray}
&&{\cal L}\varphi(s,x,u,v)
\nonumber
\\[-9pt]
\\[-9pt]
\nonumber
&&\qquad:=\frac{\partial}{\partial
s}\varphi (s,x) + \frac{1}{2}
\operatorname{tr} \bigl(\sigma\sigma^{T}(s,x, u, v)D^2
\varphi \bigr)+ D \varphi\cdot b(s,x,u, v).
\end{eqnarray}

Let us now fix arbitrarily $(t,x)\in[0,T)\times R^d$. Given an arbitrary
partition $\pi=\{0=t_0<t_1<\cdots<t_n=T\}$, we let $1\le j\le n$
be such that $t<t_j$. Let us investigate the following
BSDE defined on the interval
$[t,t_j] \dvtx $
%
%
%e4.10 #&#
%
\begin{eqnarray}
\label{BSDE1appr} %
&&dY^{1,u,v}_s = -E \bigl[{
\cal L}\varphi \bigl(s,X^{t,x;u,v}_s,u_s,v_s
\bigr)|\widetilde{\cal F}^\pi_{s} \bigr]\,ds\nonumber\\[-2pt]
&&\hspace*{44pt}{}-E \bigl[f
\bigl(s,X^{t,x;u,v}_s, Y^{1,u,v}_s +E
\bigl[\varphi \bigl(s,X^{t,x;u,v}_s \bigr)|\widetilde{\cal
F}^\pi _{s} \bigr],Z^{1,u,v}_s\nonumber\\
&&\hspace*{67pt}{}+E
\bigl[\nabla\varphi \bigl(s,X^{t,
x;u,v}_s \bigr)\sigma
\bigl(s,X^{t,x;u,v}_s,u_s,v_s \bigr) |
\widetilde{\cal F}^\pi _{s} \bigr], u_s,v_s
\bigr)|\widetilde{\cal F}^\pi_{s} \bigr]\,ds
\\
&&\hspace*{44pt}{}+Z^{1,u,v}_s\,dB_s+dM^{1,u,v}_s,\nonumber\\
\eqntext{ Y^{1,u,v}_{t_j} = 0,  M^{1,u,v} \mbox{ martingale orthogonal to } B,
M^{1,u,v}_t=0,}
\end{eqnarray}
where the process $X^{t,x;u,v}$ is the unique solution of SDE
(\ref{SDE}) and $(u,v)\in{\cal U}_{t, t_j}^\pi\times{\cal V}_{t,
t_j}^\pi$.

It can be easily verified that (or, refer to \cite{CFS}), under our
standard assumptions on the coefficients
$\sigma,b$ and $f$, the above BSDE has a unique solution $(Y^{1,u,v},Z^{1,u,
v},M^{1,u,v})$ over the time interval $[t,t_j]$.

We have the following relation between the solution $Y^{1,u,v}$
and the backward stochastic semigroup $G^{t,x;u,v}_{s,t_j} [\varphi(t_j,
X^{t,x;u,v}_{t_j})]$:\vspace*{-2pt}

%le4.2 #&#
\begin{lemma}\label{lemma_BSDE1} For every $s\in[t,t_j]$, it holds
%
%
%e4.11 #&#
%
\begin{equation}\qquad
Y^{1,u,v}_s = G^{t,x;u,v}_{s,t_j} \bigl[
\varphi \bigl(t_j,X^{t,x;u,v}_{t_j} \bigr) \bigr] -E
\bigl[\varphi \bigl(s,X^{t,x;u,v}_s \bigr) | \widetilde{\cal
F}^\pi_{s} \bigr], \qquad P\mbox{-a.s.,}
\end{equation}
and in particular, for $s=t$,
%
%
%e4.12 #&#
%
\begin{equation}
Y^{1,u,v}_t = G^{t,x;u,v}_{t,t_j} \bigl[\varphi
\bigl(t_j,X^{t,x;u,v}_{t_j} \bigr) \bigr] -
\varphi(t,x), \qquad P\mbox{-a.s.}\vspace*{-2pt}
\end{equation}
\end{lemma}

\begin{pf} Recall that $G^{t,x;u,v}_{s,t_j}[\varphi(t_j, X^{t,x;
u,v}_{t_j})]$ is defined through the BSDE
%
%
%e4.13 #&#
%
\begin{equation}
\label{ab1}\qquad \cases{ %
dY^{u,v}_s = -E \bigl[f
\bigl(s,X^{t,x;u,v}_s, Y^{u,v}_s,Z^{u,v}_s,u_s,v_s
\bigr) |\widetilde{\cal F}_s^\pi \bigr]\,ds
\vspace*{2pt}\cr
\hspace*{38pt}{}+Z^{u,v}_s \,dB_s +dM^{u,v}_s,
\vspace*{2pt}
\cr
Y^{u,v}_{t_j} = E \bigl[\varphi
\bigl(t_j,X^{t,x;u,v}_{t_j} \bigr) | \widetilde{\cal
F}^\pi_{t_j} \bigr],\qquad s \in[t,t_j],
\vspace*{2pt}
\cr
M^{u,v} \mbox{ square integrable martingale, orthogonal
to }B, M^{u,v}_t=0, }
\end{equation}
by the relation:
%
%
%e4.14 #&#
%
\begin{equation}
G^{t,x;u,v}_{s,t_j} \bigl[\varphi \bigl(t_j,X^{t,x;u,v}_{t_j}
\bigr) \bigr] =Y^{u,v}_s, \qquad s\in[t,t_j].
\end{equation}

We notice that, since $X^{t,x;u,v}$ is $\mathbb{F}^\pi$-adapted, we have
%
%
%e4.15 #&#
%
\begin{equation}
E \bigl[\varphi \bigl(s,X^{t,x;u,v}_s \bigr) | \widetilde{\cal
F}^\pi_{s} \bigr]= E \bigl[\varphi \bigl(s,X^{t,x;u,v}_s
\bigr) | {\cal F}^B_{T}\vee{\cal H}_{\ell-1}
\bigr],
\end{equation}
$s\in[t\vee t_{\ell-1},t\vee t_\ell), 1\le\ell\le j$. Hence,
with the help of the It\^{o} formula we obtain on each
interval $[t\vee t_{\ell-1},t\vee t_\ell), 1\le\ell\le j$,
%
%
%e4.16 #&#
%
\begin{eqnarray}
&&dE \bigl[\varphi \bigl(s,X^{t,x;u,v}_s \bigr) |
\widetilde{\cal F}^\pi_{s} \bigr]\nonumber\\
&&\qquad =E \bigl[{\cal L}
\varphi \bigl(s,X^{t,x;u,v}_s,u_s,v_s
\bigr) | \widetilde{\cal F}^\pi_{s} \bigr]\,ds
\\
& &\qquad\quad {}+E \bigl[\nabla\varphi \bigl(s,X^{t,
x;u,v}_s
\bigr)\sigma \bigl(s,X^{t,x;u,v}_s,u_s,v_s
\bigr) | \widetilde{\cal F}^\pi _{s} \bigr]
\,dB_s.\nonumber
\end{eqnarray}

Let us put
\[
M_s:=\sum_{\ell\dvtx  t<t_\ell\le s}\triangle E \bigl[
\varphi \bigl(t_\ell, X^{t,x;u,v}_{t_\ell} \bigr)|
\widetilde{\cal F}^\pi_{t_\ell} \bigr],\qquad s
\in[t,t_j],
\]
with
\[
\triangle E\bigl[\varphi\bigl(t_\ell, X^{t,x;u,v}_{t_\ell}
\bigr)|\widetilde{\cal F}^\pi_{t_\ell}\bigr]= E\bigl[\varphi
\bigl(t_\ell,X^{t,x;u,v}_{t_\ell}\bigr)|\widetilde{\cal
F}^\pi _{t_\ell}\bigr]- E\bigl[\varphi\bigl(t_\ell,X^{t,x;u,v}_{t_\ell}
\bigr)|\widetilde{\cal F}^\pi _{t_\ell-}\bigr].
\]

Obviously, $M$ is a pure jump martingale with respect to the
filtration $\mathbb{F}^B$ and, hence, orthogonal to $B$, and
%
%
%e4.17 #&#
%
\begin{eqnarray}
& &dE \bigl[\varphi \bigl(s,X^{t,x;u,v}_s \bigr) |
\widetilde{\cal F}^\pi_{s} \bigr]
\nonumber
\\
&&\qquad= E \bigl[{\cal L}\varphi \bigl(s,X^{t,x;u,v}_s,u_s,v_s
\bigr) | \widetilde{\cal F}^\pi_{s} \bigr]\,ds
\nonumber
\\[-8pt]
\\[-8pt]
\nonumber
& &\qquad\quad{} +E \bigl[\nabla\varphi \bigl(s,X^{t,
x;u,v}_s
\bigr)\sigma \bigl(s,X^{t,x;u,v}_s,u_s,v_s
\bigr) | \widetilde{\cal F}^\pi _{s} \bigr]
\,dB_s+dM_s, \\
\eqntext{s\in[t,t_j].}
\end{eqnarray}
Consequently, $(Y^{u,v}_s-E[\varphi
(s,X^{t,x;u,v}_s)|\widetilde{\cal F}^\pi_{s}],Z^{u,v}_s-
E [\nabla\varphi(s,X^{t,x;u,v}_s)\sigma(s, X^{t,x;u,v}_s,\break  u_s,v_s)
|\widetilde{\cal F}^\pi_{s} ],
M^{u,v}_s-M_s), t\le s\le t_j$, is a solution
of BSDE (\ref{BSDE1appr}). From its uniqueness, we can conclude the
statement of the lemma.
\end{pf}

Let us now simplify the preceding BSDE (\ref{BSDE1appr}) by replacing the
process $X^{t,x;u,v}$ by its initial value $x$. Then BSDE (\ref{BSDE1appr})
takes the form
%
%
%e4.18 #&#
%
\begin{equation}
\label{BSDE2appr}\qquad  \cases{ %
dY^{2,u,v}_s = -E
\bigl[F \bigl(s,x,Y^{2,u,v}_s,Z^{2,u,v}_s,u_s,
v_s \bigr)| \widetilde{\cal F}^\pi_s \bigr]
\,ds\vspace*{2pt}\cr
\hspace*{45pt}{}+Z^{2,u,v}_s \,dB_s+dM_s^{2,u,v},
\vspace*{2pt}
\cr
Y^{2,u,v}_{t_j} = 0, \qquad s
\in[t,t_j], \vspace*{2pt}
\cr
M^{u,v} \mbox{ square
integrable martingale, orthogonal to }B, M^{u,v}_t=0,}
\end{equation}
where $(u,v)\in{\cal U}_{t,t_j}^\pi\times{\cal
V}_{t,t_j}^\pi
$. As in
the discussion of BSDE (\ref{BSDE1appr}) we see that the above BSDE has
a unique solution. From the BSDEs (\ref{BSDE1appr}) and
(\ref{BSDE2appr}), we have
the following lemma.
%
%le4.3 #&#
\begin{lemma}\label{lemma_BSDE2} For every $(u,v)\in{\cal
U}_{t,t_j}^\pi\times{\cal V}_{t,
t_j}^\pi$ we have
%
%
%e4.19 #&#
%
\begin{equation}
\bigl|Y^{1,u,v}_t-Y^{2,u,v}_t \bigr| \leq
C(t_j-t)^{{3}/{2}},\qquad P\mbox{-a.s.},
\end{equation}
where C is independent of the control processes $u$ and $v$, but also
independent of the
partition $\pi$.
\end{lemma}

\begin{pf} Let $(u,v)\in{\cal U}_{t,t_j}^\pi\times{\cal V}_{t,
t_j}^\pi$. Then, for all $s\in[t,t_j]$, thanks to Condition~\ref{condition},
%
%
%e4.20 #&#
%
\begin{eqnarray}
& & E \bigl[{\cal L}\varphi(s,x,u_s,v_s)
\nonumber\\
&&\hspace*{12pt}{}+f \bigl(s,x,y+\varphi(s,x),z+ E \bigl[\nabla\varphi(s,x) \sigma(s,x,u_s,v_s)|
\widetilde{\cal F}^\pi_s \bigr],u_s,
v_s \bigr)|\widetilde{\cal F}^\pi_s \bigr]
\nonumber
\\
&&\qquad= E \bigl[{\cal L}\varphi(s,x,u_s,v_s)\\
&&\hspace*{44pt}{} +f
\bigl(s,x,y+\varphi(s,x),z+ \nabla\varphi(s,x)\sigma(s,x,u_s,v_s),u_s,
v_s \bigr)|\widetilde{\cal F}^\pi_s \bigr]
\nonumber\\
& &\qquad = E \bigl[F(s,x,y,z,u_s,v_s)|\widetilde{\cal
F}^\pi_s \bigr],\qquad P\mbox{-a.s.}
\nonumber
\end{eqnarray}
Consequently, we have to compare the solution of BSDE
(\ref{BSDE1appr})
%
%
%e4.21 #&#
%
\begin{eqnarray}
&&dY^{1,u,v}_s\nonumber\\
&&\qquad =-E \bigl[{\cal L}\varphi
\bigl(s,X^{t,x;u,v}_s,u_s,v_s \bigr)\nonumber\\
&&\qquad\hspace*{30pt}{}+f
\bigl(s,X^{t,x;u,v}_s,Y^{1,u,v}_s+E \bigl[
\varphi \bigl(s,X^{t,x;u,v}_s \bigr)|\widetilde{\cal
F}^\pi_{s} \bigr], Z^{1,u,v}_s
\\
&&\qquad\hspace*{51pt}{}+E
\bigl[\nabla\varphi \bigl(s,X^{t,
x;u,v}_s \bigr)\sigma
\bigl(s,X^{t,x;u,v}_s,u_s,v_s \bigr) |
\widetilde{\cal F}^\pi _{s} \bigr], u_s,v_s
\bigr)|\widetilde{\cal F}^\pi_{s} \bigr]\,ds
\nonumber\\
&&\hspace*{10pt}\qquad{}+Z^{1,u,v}_s\,dB_s+dM^{1,u,v}_s,\qquad
 Y^{1,u,v}_{t_j} =0,\nonumber
\end{eqnarray}
with that of BSDE (\ref{BSDE2appr}) which can be rewritten as
%
%
%e4.22 #&#
%
\begin{eqnarray}
dY^{2,u,v}_s & = & -E \bigl[{\cal L}
\varphi(s,x,u_s,v_s)\nonumber\\
&&\hspace*{19pt}{}+f \bigl(s,x,Y^{2,u,v}_s
+E \bigl[\varphi(s,x)|\widetilde{\cal F}^\pi_{s} \bigr],
Z^{2,u,v}_s
\nonumber
\\[-8pt]
\\[-8pt]
\nonumber
&&\hspace*{41pt}{} +E \bigl[\nabla\varphi(s,x)
\sigma(s,x,u_s,v_s) |\widetilde{\cal
F}^\pi_{s} \bigr], u_s,v_s
\bigr)|\widetilde{\cal F}^\pi_{s} \bigr]\,ds
\\
&&{}+Z^{2,u,v}_s\,dB_s+ dM^{2,u,v}_s,
\qquad
Y^{2,u,v}_{t_j} =  0,\nonumber
\end{eqnarray}
and from BSDE standard estimates we deduce
%
%
%e4.23 #&#
%
\begin{eqnarray}
&& \bigl|Y^{1,u,v}_t-Y^{2,u,v}_t
\bigr|^2+E \biggl[\int_t^{t_j}
\bigl|Z^{1,u,v}_r-Z^{2,u,v}_r
\bigr|^2\,dr\Big| \widetilde{\cal F}_t^\pi \biggr]\nonumber\\
&&\quad{} +E
\biggl[\sum_{\ell\le j; t<t_\ell} \bigl|\triangle M^{1,u,v}_{t_\ell}-
\triangle M^{2,u,v}_{t_\ell
} \bigr|^2\Big| \widetilde {\cal
F}_t^\pi \biggr]
\nonumber
\\[-8pt]
\\[-8pt]
\nonumber
&&\qquad \le CE \biggl[ \biggl(\int_t^{t_j}\bigl|X^{t,x;u,v}_r-
x\bigr|\,dr \biggr)^2\Big|\widetilde{\cal F}_t^\pi
\biggr]
\\
&&\qquad \le C(t_j-t)^3,\qquad P\mbox{-a.s.},\nonumber
\end{eqnarray}
%+
where the constant $C$ depends only on the boundedness
and Lipschitz constants of the coefficients and the derivatives of
$\varphi$, but not on $j$ nor the considered partition~$\pi$.
\end{pf}

Let us now state the following crucial lemma which, although
inspired by Lemma 4.3 in \cite{Buckdahn-Li-2008}, differs heavily
because of the different framework studied here.

%le4.4 #&#
\begin{lemma}\label{lemma_BSDE3} Let $Y^0=(Y^0_s)_{s\in[t,t_j]}$ denote
the unique solution of
the following ordinary backward differential equation:
%
%
%e4.24 #&#
%
\begin{equation}
\cases{- {\dot{Y}}^0_s = F_0
\bigl(s,x,Y^0_s,0 \bigr), &\quad $s \in[t,t_j]$,
\vspace*{2pt}
\cr
Y^0_{t_j} = 0, }
\end{equation}
where, for $ (s,y,z)\in[t,t_j]\times R\times R^d$,
%
%
%e4.25 #&#
%
\begin{eqnarray}
\label{qq} %
F_0(s,x,y,z)&: =& \sup_{\mu\in{\cal P}(U)}
\Bigl(\inf_{v\in V}F(s,x,y,z,\mu,v) \Bigr)
\nonumber
\\
& = & \sup_{\mu\in{\cal P}(U)} \Bigl(\inf_{\nu\in{\cal
P}(V)}F(s,x,y,z,
\mu,\nu) \Bigr) \\
&&{}\times\Bigl(=\inf_{\nu\in{\cal P}(V)}\sup_{\mu\in{\cal
P}(U)}F(s,x,y,z,
\mu,\nu) \Bigr).\nonumber
\end{eqnarray}

Then, for all $s\in[t,t_j]$, $P$-a.s.,
%
%
%e4.26 #&#
%
\begin{equation}
Y^0_s =  \esssup_{u \in{\cal U}_{t, t_j}^\pi}
\essinf_{v \in
{\cal V}_{t,t_j}^\pi} Y^{2,u,v}_s
= \essinf_{v \in
{\cal V}_{t,t_j}^\pi}\esssup_{u \in{\cal U}_{t, t_j}^\pi}Y^{2,u,v}_s.
\end{equation}
\end{lemma}

\begin{pf}\textit{Step} 1. Given $(u,v)\in{\cal U}_{t,t_j}^\pi
\times{\cal V}_{t,t_j}^\pi$, let
$(Y^{2,u,v}, Z^{2,u,v}, M^{2,u,v})$ be the unique solution of BSDE
(\ref{BSDE2appr}). We recall that, for all $s\in[t\vee t_{\ell-1},
t_\ell), (1\le\ell\le j)$, $(Y^{2,u,v}_s, Z^{2,u,v}_s,M^{2,u,v}_s)$
is $\widetilde{\cal F}^\pi_s(={\cal F}_s^B\vee{\cal H}_{\ell
-1})$-measurable,
$u_s$ is ${\cal
F}^{\pi,1}_s(={\cal F}_s^B\vee{\cal H}_{\ell-1}\vee\sigma\{\zeta
_{\ell,
1}\})$-measurable and $v_s$ is ${\cal F}^{\pi,2}_s(={\cal
F}_s^B\vee{\cal H}_{\ell-1}\vee\sigma\{\zeta_{\ell, 2}\})$-measurable.
Consequently, knowing $\widetilde{\cal F}^\pi_s$, $u_s$ and $v_s$
are conditionally independent, and defining
\begin{eqnarray}
\mu_s^u(A):=P\bigl\{u_s\in A|
\widetilde{\cal F}^\pi_s\bigr\} , \qquad
\nu_s^v(B):=P\bigl\{v_s\in B | \widetilde{
\cal F}^\pi_s\bigr\}, \nonumber\\
 \eqntext{A\in{\cal B} (U), B
\in{\cal B}(V),}
\end{eqnarray}
we have
%
%
%e4.27 #&#
%
\begin{eqnarray}
& & E \bigl[F \bigl(s,x,Y^{2,u,v}_s
,Z^{2,u,v}_s,u_s, v_s \bigr) |
\widetilde{\cal F}^\pi_s \bigr]
\nonumber
\\
&&\qquad = F \bigl(s,x,Y^{2,u,v}_s
,Z^{2,u,v}_s,\mu^u_s,
\nu^v_s \bigr) \\
&&\qquad\quad{}\times\biggl(:= \int_{U\times V}F
\bigl(s,x,Y^{2,u,v}_s ,Z^{2,u,v}_s,u',
v' \bigr) \mu^u_s\otimes
\nu^v_s \bigl(du'\,dv'
\bigr) \biggr).\nonumber
\end{eqnarray}

Indeed, this relation can be easily checked by considering
first instead of $F(s,x,Y^{2,u,v}_s,Z^{2,u,v}_s, u_s, v_s)$
integrands of the form $\xi_s f_1(u_s)f_2(v_s)$, $\xi_s\in L^\infty
(\Omega$, $\widetilde{\cal F}^\pi_s,P)$ and $f_1,f_2$
bounded Borel functions over $U$ and $V$, respectively, and applying
later a Monotonic Class theorem.

Hence, with the notation
$ F(s,x,y,z,\mu,v):=\int_UF(s,x,y,z,u',v)\mu(du')$, $
\mu\in{\cal P}(U)$, and with putting
\[
F_1(s,x,y,z,\mu):=\inf_{v\in V} F(s,x,y,z,\mu,v)
\Bigl(=\inf_{\nu\in{\cal P}(V)}F(s,x,y,z,\mu,\nu ) \Bigr),
\]
$(s,y,z,\mu)\in[0,T]\times R\times R^d\times{\cal P}(U)$,
we obtain
%
%
%e4.28 #&#
%
\begin{eqnarray}
&& E \bigl[F \bigl(s,x,Y^{2,u,v}_s
,Z^{2,u,v}_s,u_s, v_s \bigr) |
\widetilde{\cal F}^\pi_s \bigr]
\nonumber\\
&&\qquad= \int_{V} F \bigl(s,x,Y^{2,u,v}_s,Z^{2,u,v}_s,
\mu_s^u, v' \bigr)\nu^v_s
\bigl(dv' \bigr)
\\
& &\qquad\ge F_1 \bigl(s,x,Y^{2,u,v}_s,Z^{2,u,v}_s,
\mu_s^u \bigr), \qquad ds\,dP\mbox{-a.e.}\nonumber
\end{eqnarray}

Consequently, denoting by $(Y^{3,u},Z^{3,u},M^{3,u})\in{\cal
S}^2_{\widetilde{\mathbb{F}}^\pi}(t,t_j; R)\times
L_{\widetilde{\mathbb{F}}^\pi}^2(t,t_j;\break R^d)\times
{\cal M}^2_{\widetilde{\mathbb{F}}^\pi}(t,t_j;R)$ the unique solution
of the
BSDE
%
%
%e4.29 #&#
%
\begin{equation}
\label{BSDE3appr}\qquad  \cases{ %
dY^{3,u}_s =
-F_1 \bigl(s,x,Y^{3,u}_s,Z^{3,u}_s,
\mu^u_s \bigr)\,ds\vspace*{2pt}\cr
\hspace*{38pt}{} +Z^{3,u}_s
\,dB_s+dM^{3,u}_s,\qquad s
\in[t,t_j], \vspace*{2pt}
\cr
Y^{3,u}_{t_j} = 0,
\vspace*{2pt}
\cr
M^{3,u} \mbox{ square integrable martingale, orthogonal
to }B, M^{3,u}_t=0,}
\end{equation}
we deduce from the comparison theorem for BSDEs (refer to
\cite{CFS}, for classical case it can be referred to \cite{Pe1}, or
\cite{Buckdahn-Li-2008}) that $Y^{2,u,v}_s
\ge Y^{3,u}_s, s\in[t,t_j]$, $P$-a.s., for all $v\in{\cal V}^\pi_{t,t_j}$.
For this, we observe that $F_1(s,x,y,z,\mu)$ is a jointly continuous
function over $[0,T]\times R^d\times R\times R^d\times{\cal P}(U)$,
which is
Lipschitz in $(y,z)$, uniformly with respect to $(s,x,\mu)$. Thus, taking
into account the arbitrariness of $v\in{\cal V}_{t,t_j}^\pi$, we deduce
%
%
%e4.30 #&#
%
\begin{equation}
\label{BSDE3_inequ} Y^{3,u}_s\le
\essinf_{v\in
{\cal
V}_{t,t_j}^\pi}Y^{2,u,v}_s, \qquad P\mbox{-a.s, } s
\in[t,t_j].
\end{equation}

Let us show that we have even equality in the above inequality.
For this we observe that, since the function $F$ is continuous over $[t,t_j]
\times R^d\times R\times R^d\times{\cal P}(U)\times V$, there exists a Borel
measurable function $v^*\dvtx [t,t_j]\times R\times R^d\times{\cal
P}(U)\rightarrow
V$ such that
\[
F_1(s,x,y,z,\mu)=\inf_{v\in V}F(s,x,y,z,\mu,v)= F
\bigl(s,x,y,z,\mu, v^*(s,y,z,\mu)\bigr),
\]
$(s,y,z,\mu)\in[t,t_j]\times R\times R^d\times{\cal P}(U)$.
With the help of this measurable function, we introduce the control
process $v_s^*:=v^*(s,Y_s^{3,u},Z_s^{3,u},\mu_s^u), s\in[t,t_j]$. We notice
that $v^*=(v^*_s)_{s\in[t,t_j]}$ belongs to ${\cal V}_{t,t_j}^\pi$ and
is even
$\widetilde{\mathbb{F}}^\pi$-adapted. Thus,
%
%
%e4.31 #&#
%
\begin{eqnarray}
&&E \bigl[F \bigl(s,x,Y^{3,u}_s
,Z^{3,u}_s,u_s, v_s^* \bigr) |
\widetilde{\cal F}^\pi_s \bigr]\nonumber\\
&&\qquad =  F
\bigl(s,x,Y^{3,u}_s,Z^{3,u}_s,
\mu^u_s, v_s^* \bigr)
\\
&&\qquad =  F_1 \bigl(s,x,Y^{3,u}_s
,Z^{3,u}_s,\mu^u_s \bigr), \qquad ds
\,dP \mbox{-a.e.},\nonumber
\end{eqnarray}
from where we see that $(Y^{3,u},Z^{3,u},M^{3,u})$ is a
solution of BSDE
(\ref{BSDE2appr}) driven by the couple $(u,v^*)\in{\cal
U}_{t,t_j}^\pi
\times
{\cal V}_{t,t_j}^\pi$ of admissible controls. Consequently, the
uniqueness of the
solution of BSDE (\ref{BSDE2appr}) yields that $Y^{2,u,v^*}_s=Y^{3,u}_s,
s\in[t,t_j]$, and from~(\ref{BSDE3_inequ}) we obtain:
%
%
%e4.32 #&#
%
\begin{equation}
\label{BSDE3-1_inequ} Y^{3,u}_s=
\essinf_{v\in
{\cal
V}_{t,t_j}^\pi}Y^{2,u,v}_s, \qquad P\mbox{-a.s, } s
\in[t,t_j], u\in{\cal U}_{t,t_j}^\pi.
\end{equation}

\textit{Step} 2. We begin with showing the latter relation
in (\ref{qq}). For this end we remark that, for all $(s,y,z)$, the
function $(\mu,\nu)\rightarrow F(s,x,y,z,\mu,\nu)=\int_{U}\int_{ V}
F(s,x,y,z,u,v)\nu(dv)\mu(du)$, $(\mu,\nu)\in{\cal P}(U)\times
{\cal P}(V)$,
is bi-linear and,\break  hence, concave-convex in $(\mu,\nu)$ belonging
to the cross product ${\cal P}(U)\times{\cal P}(V)$ of two convex compact
spaces. Consequently, this mapping admits a saddle point, and it
follows in
particular that the order of $\sup_{\mu\in{\cal P}(U)}$ and
$\inf_{\nu\in{\cal P}(V)}$ is exchangeable without changing the value
of $F_0(s,x,y,z)$.

Let us now consider an arbitrary $u\in{\cal U}_{t,t_j}^\pi$. From the
definition of the function $F_0(s,x,y,z)$ and that of
$F_1(s,x,y,z,\mu)$, we have
%
%
%e4.33 #&#
%
\begin{eqnarray}
&&F_0(s,x,y,z)\nonumber\\
&&\qquad  =  \sup_{\mu\in{\cal P}(U)}F_1(s,x,y,z,
\mu)
\\
&&\qquad \ge F_1 \bigl(s,x,y,z,\mu^u_s \bigr),\qquad
(s,y,z) \in[t,t_j]\times R\times R^d,  u\in{\cal
U}_{t,t_j}^\pi.\nonumber
\end{eqnarray}

Consequently, since $(Y^0,0)$ can be regarded as the solution
of the BSDE
\[
dY^0_s=-F_0 \bigl(s,x,Y^0_s,0
\bigr)\,ds+0\cdot \,dB_s,\qquad s\in[t,t_j],
Y^0_{t_j}=0,
\]
we get from the comparison theorem for BSDEs that $Y^0_s
\ge Y^{3,u}_s, s\in[t,t_j], P$-a.s. Hence, in view of the
arbitrariness of the choice of $u\in{\cal U}_{t,t_j}^\pi$, it follows that
%
%
%e4.34 #&#
%
\begin{equation}
\label{Y_0} Y^{0}_s\ge\esssup_{u \in{\cal U}_{t,
t_j}^\pi}Y^{3,u}_s,
\qquad P\mbox{-a.s.}, s\in[t,t_j].
\end{equation}

It remains to prove that we have even equality in this latter
relation. For this end, we notice that thanks to the uniform continuity of
the function $(s,y,\mu)\rightarrow F_1(s,x,y,0,\mu)$ over
$[t,t_j]\times R
\times{\cal P}(U)$ [we note that $x$ in $F_1(s,x,y$, $0,\mu)$ is fixed],
and the compactness of ${\cal P}(U)$ endowed with the topology
generated by
the weak convergence, we have the existence of a Borel measurable selection
$\mu^*=(\mu^*(s,y))\dvtx [t,t_j]\times R\rightarrow{\cal P}(U)$
such that
\[
F_0(s,x,y,0)=F_1 \bigl(s,x,y,0,\mu^*(s,y) \bigr), \qquad (s,y)
\in[t,t_j]\times R.
\]

Again from the uniform continuity of $(s,y,\mu)\rightarrow
F_1(s,x,y,0,\mu)$, we get that, for
arbitrarily given $\varepsilon>0$ there is some
$\delta(=\delta_\varepsilon)>0$ such that $|F_1(s,x,\break y,0,\mu)-F_1(s',x,y',
0,\mu)|\le\varepsilon$, for all $\mu\in{\cal P}(U)$ and all $(s,y),
(s',y')$ with $|(s,y)-(s',y')|\le\delta$. Let $(\Delta_\ell)_{\ell
\ge1}$ be a Borel partition of the set $[t,t_j]\times R$, composed of
nonempty sets $\Delta_\ell$ with diameter less than or equal to
$\delta$.
For every $\ell\ge1$, let us fix arbitrarily an element $(s_\ell
,y_\ell)$
of $\Delta_\ell$, and let us put $\mu_\ell:=\mu^*(s_\ell,y_\ell)$.
Moreover, let us consider an independent sequence of random variables
$\xi_\ell\in L^0(\Omega,\sigma\{\zeta_{j,1}\},P;U)$ such that,
for all $\ell\ge1$, the law $P\circ[\xi_\ell]^{-1}$ coincides with~$\mu_\ell$.

With the above introduced quantities, we define the control process
\[
u_s^*:=\sum_{\ell\ge1}I \bigl\{
\bigl(s,Y^0_s \bigr)\in \Delta_\ell \bigr\}
\cdot \xi_\ell, \qquad s\in[t,t_j].
\]

Such defined process belongs, obviously, to ${\cal U}_{t,
t_j}^\pi$. Moreover, we observe that, for all $s\in[t,t_j]$, $u^*_s$ is
$\sigma\{\zeta_{j,1}\}$-measurable and, consequently, independent of
$\widetilde{\cal F}_s^\pi$. Hence, for all $A\in{\cal B}(U)$,
%
%
%e4.35 #&#
%
\begin{eqnarray}
\mu^{u^*}_s(A) & = & P \bigl
\{u^*_s\in A | \widetilde{\cal F}_s^\pi
\bigr\}= P \bigl\{u^*_s\in A \bigr\}
\nonumber
\\[-8pt]
\\[-8pt]
\nonumber
& = & \sum_{\ell\ge1}I \bigl\{
\bigl(s,Y^0_s \bigr)\in\Delta_\ell \bigr\}
\mu_\ell(A).
\end{eqnarray}

It follows that $\mu^{u^*}_s= \sum_{\ell\ge1}I\{(s,Y^0_s)
\in\Delta_\ell\}\mu_\ell$. Hence, due to our choice of the partition
$\Delta_\ell, \ell\ge1$,
%
%
%e4.36 #&#
%
\begin{eqnarray}
\label{F_0_F_1} %
F_0 \bigl(s,x,Y_s^0,0 \bigr) & \le&
\varepsilon+\sum_{\ell\ge1} I \bigl\{
\bigl(s,Y^0_s \bigr)\in\Delta_\ell \bigr
\}F_0(s_\ell,x,y_\ell,0)
\nonumber\\
& = & \varepsilon+\sum_{\ell\ge1} I \bigl\{
\bigl(s,Y^0_s \bigr)\in\Delta_\ell \bigr
\}F_1(s_\ell,x,y_\ell,0,\mu_\ell)
\\
& = & \varepsilon+\sum_{\ell\ge1} I \bigl\{
\bigl(s,Y^0_s \bigr)\in\Delta_\ell \bigr
\}F_1 \bigl(s_\ell,x,y_\ell,0,
\mu_s^{u^*} \bigr)
\nonumber\\
& \le& 2\varepsilon+F_1 \bigl(s,x,Y^0_s,0,
\mu_s^{u^*} \bigr), \qquad s\in[t,t_j].\nonumber
\end{eqnarray}

Let us compare now $Y^0$ with the solution $(Y^{3,u^*},Z^{3,u^*})$ of BSDE
(\ref{BSDE3appr}) controlled by $u^*\in{\cal U}_{t,t_j}^\pi$. Obviously,
\begin{eqnarray*}
d \bigl(Y^0_s-Y^{3,u^*}_s \bigr)&=&-
\bigl(F_0 \bigl(s,x,Y^0_s,0
\bigr)-F_1 \bigl(s,x,Y^{3,u^*}_s,
Z^{3,u^*}_s,\mu_s^{u^*} \bigr) \bigr)
\,ds\\
&&{}-Z^{3,u^*}_s\,dB_s-dM^{3,u^*}_s,
\end{eqnarray*}
$s\in[t,t_j], Y^0_{t_j}-Y^{3,u^*}_{t_j}=0$,
and from the It\^{o} formula,
%
%
%e4.37 #&#
%
\begin{eqnarray}
& &d \bigl( \bigl(Y^0_s-Y^{3,u^*}_s
\bigr)^+ \bigr)^2
\nonumber\\
&&\qquad = -2 \bigl(Y^0_s-Y^{3,u^*}_s
\bigr)^+ \bigl(F_0 \bigl(s,x,Y^0_s,0
\bigr)-F_1 \bigl(s,x,Y^{3,u^*}_s,Z^{3,u^*}_s,
\mu_s^{u^*} \bigr) \bigr)\,ds
\nonumber
\\[-8pt]
\\[-8pt]
\nonumber
& &\qquad\quad{} +\bigl|Z^{3,u^*}_s\bigr|^2I \bigl
\{Y^0_s-Y^{3,u^*}_s>0 \bigr\}\,ds -
2 \bigl(Y^0_s-Y^{3,u^*}_s
\bigr)^+Z^{3,u^*}_s\,dB_s
\\
& & \qquad\quad{}+ I \bigl\{Y^0_s-Y^{3,u^*}_s>0
\bigr\}\,d \bigl[M^{3,u^*} \bigr]_s -2 \bigl(Y^0_s-Y^{3,u^*}_s
\bigr)^+dM^{3,u^*}_s,\nonumber
\end{eqnarray}
and from standard estimates combined with (\ref{F_0_F_1})
we get
%
%
%e4.38 #&#
%
\begin{eqnarray}
& & \bigl( \bigl(Y^0_s-Y^{3,u^*}_s
\bigr)^+ \bigr)^2+E \biggl[\int_s^{t_j}
\bigl|Z^{3,u^*}_r\bigr|^2I \bigl\{Y^0_r-Y^{3,u^*}_r>0
\bigr\}\,dr \nonumber\\
&&\hspace*{101pt}{}+\int_{(s,t_j]} I \bigl\{Y^0_r-Y^{3,u^*}_r>0
\bigr\}\,d \bigl[M^{3,u^*} \bigr]_r \Big| \widetilde{\cal
F}_s^\pi \biggr]
\nonumber\\
&&\qquad = 2 E \biggl[\int_s^{t_j}
\bigl(Y^0_r-Y^{3,u^*}_r \bigr)^+
\bigl(F_0 \bigl(r,x,Y^0_r,0 \bigr)
\nonumber\\
&&\hspace*{69pt}{}-F_1 \bigl(r,x,Y^{3,u^*}_r,Z^{3,u^*}_r,
\mu_r^{u^*} \bigr) \bigr)\,dr\Big | \widetilde{\cal
F}_s^\pi \biggr]
\nonumber\\
&&\qquad \le 2 E \biggl[\int_s^{t_j}
\bigl(Y^0_r-Y^{3,u^*}_r \bigr)^+
\bigl(2\varepsilon+F_1 \bigl(r,x,Y^0_r,0,
\mu_r^{u^*} \bigr)\\
&&\hspace*{82pt}{} -F_1 \bigl(r,x,Y^{3,u^*}_r,
Z^{3,u^*}_r,\mu_r^{u^*} \bigr)
\bigr)\,dr \Big| \widetilde{\cal F}_s^\pi \biggr]
\nonumber\\
&&\qquad \le 2 E \biggl[\int_s^{t_j}
\bigl(Y^0_r-Y^{3,u^*}_r \bigr)^ +
\bigl(2\varepsilon+C\bigl|Y^0_r-Y^{3,u^*}_r\bigr|+C\bigl|Z^{3,u^*}_r\bigr|
\bigr)\,dr\Big | \widetilde{\cal F}_s^\pi \biggr]
\nonumber\\
&& \qquad \le\varepsilon^2+CE \biggl[\int_s^{t_j}
\bigl( \bigl(Y^0_r- Y^{3,u^*}_r
\bigr)^+ \bigr)^2\,dr\Big | \widetilde{\cal F}_s^\pi
\biggr] \nonumber\\
&&\qquad\quad{}+\frac12 E \biggl[\int_s^{t_j}\bigl|Z^{3,u^*}_r\bigr|^2I
\bigl\{ Y^0_r-Y^{3,u^*}_r>0 \bigr
\}\,dr \Big| \widetilde{\cal F}_s^\pi \biggr].\nonumber
\end{eqnarray}

Hence, from Gronwall's lemma, we see that,
for some constant $C$ independent of $\varepsilon$,
$(Y^0_s-Y^{3,u^*}_s)^+\le C\varepsilon$, $s\in[t,t_j]$, that is,
\[
Y^0_s\le Y^{3,u^*}_s+C\varepsilon,
s\in[t,t_j], \qquad P\mbox{-a.s.}
\]

This latter relation together with (\ref{Y_0}) yields
\[
Y^{0}_s=\esssup_{u \in{\cal U}_{t,
t_j}^\pi}Y^{3,u}_s,\qquad
P\mbox{-a.s.}, s\in[t,t_j].
\]

Recalling the result of step 1 we can conclude the first relation
of the lemma. The second one follows by a symmetric argument.
\end{pf}

After the above auxiliary lemmas, we are now able to
characterize the functions~$W$ and $U$ introduced by Lemma
\ref{limit_W_pi_n} as viscosity solution of PDE (\ref{HJBI}).

%le4.5 #&#
\begin{lemma}\label{viscosity solution}
The functions $W, U\dvtx [0,T]\times R^d\rightarrow R$ coincide and solve
PDE~(\ref{HJBI}) in viscosity sense.
\end{lemma}

\begin{pf} \textit{Step} 1. Let us show in this step that
the function $W$ introduced in Lemma \ref{limit_W_pi_n} as the
uniform limit on compacts of a suitable sequence of lower value
functions $W^{\pi_n}, n\ge1$, is a viscosity supersolution of
(\ref{HJBI}).

For this, we fix arbitrarily $(t,x)\in[0,T)\times R^d$ and we let
$\varphi\in C_{\ell,b}^{3}([0,T]\times R^d)$ be such that $W-\varphi
\ge W(t,x)-\varphi(t,x)=0$ on $[0,T)\times R^d$. Let $\rho>0$ be
arbitrarily small and $K>0$ sufficiently large. Since $W^{\pi_n},
n\ge1$, converges uniformly on compacts to $W$, there is some
$n_{\rho,K}\ge1$ such that, for all $n\ge n_{\rho,K}$, $|W(s,x')
-W^{\pi_n}(s,x')|\le\rho$, for every $(s,x')\in[0,T]\times R^d$ with
$|x'-x|\le K$. Then it follows from the DPP (Theorem \ref{DPP}) that,
for all $n\ge n_{\rho,K}$ and every $t_j^n\in\pi_n$ with $t<t_j^n\le T$,\vspace*{1pt}
%
%
%e4.39 #&#
%
\begin{eqnarray}
\varphi(t,x)+\rho &=&  W(t,x)+\rho\nonumber\\[1pt]
&\ge &W^{\pi_n}(t,x)
\\[1pt]
&=&\esssup_{\alpha\in\mathcal{A}^{\pi_n}_{t,t_j^n}} \essinf _{\beta\in
\mathcal{B}^{\pi_n}_{t,t_j^n}}G^{t,x;\alpha,\beta}_{t,t_j^n}
\bigl(W^{\pi_n} \bigl(t_j^n,X_{t_j^n}^{t,x;\alpha,\beta}
\bigr) \bigr).\nonumber\vspace*{1pt}
\end{eqnarray}

On the other hand, taking into account that the functions
$W^{\pi_n}, n\ge1$, are bounded, uniformly with respect to $n\ge1$ and $W$ is bounded,
we have, for some constant $C_0$ (independent of $n$),\vspace*{1pt}
%
%
%e4.40 #&#
%
\begin{eqnarray}
&& W^{\pi_n} \bigl(t_j^n,X_{t_j^n}^{t,x;\alpha,\beta}
\bigr) \nonumber\\[1pt]
&&\qquad \ge W \bigl(t_j^n,X_{t_j^n}^{t,
x;\alpha,\beta}
\bigr)-\rho-2C_0I \bigl\{\bigl|X_{t_j^n}^{t,x;\alpha,\beta}-x\bigr|>K \bigr
\}
\\[1pt]
&&\qquad \ge \varphi \bigl(t_j^n,X_{t_j^n}^{t,x;\alpha,
\beta}
\bigr)-\rho-2C_0I \bigl\{\bigl|X_{t_j^n}^{t,x;\alpha,\beta}-x\bigr|>K \bigr
\},\nonumber\vspace*{1pt}
\end{eqnarray}
for all $\alpha\in\mathcal{A}^{\pi_n}_{t,t_j^n}, \beta
\in
\mathcal{B}^{\pi_n}_{t,t_j^n}$, and from the comparison theorem as well
as BSDE standard estimates (refer to \cite{CFS}) applied
to the BSDE defining our backward stochastic semigroup we obtain\vspace*{1pt}
%
%e4.41 #&#
%
\begin{eqnarray}
\varphi(t,x)+\rho&\ge&\esssup_{\alpha\in\mathcal{A}^{\pi_n}_{t,t_j^n}}
\essinf_{\beta\in\mathcal{B}^{\pi_n}_{t,t_j^n}}G^{t,x;\alpha
,\beta}_{t,
t_j^n}
\bigl(W^{\pi_n} \bigl(t_j^n,X_{t_j^n}^{t,x;\alpha,\beta}
\bigr) \bigr)
\nonumber\\[1pt]
& \ge&\esssup_{\alpha\in\mathcal{A}^{\pi_n}_{t,t_j^n}} \essinf _{\beta
\in
\mathcal{B}^{\pi_n}_{t,t_j^n}}\nonumber\\[1pt]
&&{}\times G^{t,x;\alpha,\beta
}_{t,t_j^n}
\bigl(\varphi \bigl(t_j^n, X_{t_j^n}^{t,x;\alpha,\beta}
\bigr)-\rho -2C_0I \bigl\{\bigl|X_{t_j^n}^{t,x;\alpha,
\beta}-x\bigr|>K
\bigr\} \bigr)
\nonumber
\\[-8pt]
\\[-8pt]
\nonumber
& \ge&\esssup_{\alpha\in\mathcal{A}^{\pi_n}_{t,t_j^n}} \essinf _{\beta
\in
\mathcal{B}^{\pi_n}_{t,t_j^n}} G^{t,x;\alpha,\beta
}_{t,t_j^n}
\bigl(\varphi \bigl(t_j^n, X_{t_j^n}^{t,x;\alpha,\beta}
\bigr) \bigr) \nonumber\\[1pt]
&&{}-\esssup_{\alpha\in\mathcal{A}^{\pi_n}_{t,t_j^n},\beta\in
\mathcal{B}^{\pi_n}_{t,t_j^n}}\nonumber\\[1pt]
&&{}\times L \bigl(E \bigl[ \bigl(
\rho+2C_0I \bigl\{ \bigl|X_{t_j^n}^{t,x;\alpha,
\beta}-x\bigr|>K \bigr\}
\bigr)^2|\widetilde{\cal F}_t^{\pi_n} \bigr]
\bigr)^{1/2},\nonumber
\end{eqnarray}
where the constant $L$ depends only on the coefficient $f$.
However, since
%
%
%e4.42 #&#
%
\begin{eqnarray}
& & E \bigl[ \bigl(\rho+2C_0I \bigl\{\bigl|X_{t_j^n}^{t,x;\alpha,
\beta}-x\bigr|>K
\bigr\} \bigr)^2|\widetilde{ \cal F}_t^{\pi_n}
\bigr]
\nonumber\\
&&\qquad\le 2\rho^2+8C_0^2\frac{1}{K^2}E
\bigl[\bigl|X_{t_j^n}^{t,x;\alpha,
\beta}-x\bigr|^2|\widetilde{\cal
F}_t^{\pi_n} \bigr]
\\
&&\qquad\le 2\rho^2+\frac{C}{K^2}\qquad (\alpha,\beta)\in
\mathcal{A}^{\pi_n}_{t,
t_j^n}\times\mathcal{B}^{\pi_n}_{t,t_j^n},
n\ge1\nonumber
\end{eqnarray}
(Recall that the coefficients $\sigma$ and $b$ of the
dynamics of
the game are bounded), we get for $K:=1/\rho$, for all $n\ge n_\rho:=
n_{\rho,K}$,
%
%
%e4.43 #&#
%
\begin{equation}
\qquad
\varphi(t,x)+C\rho\ge\esssup_{\alpha\in\mathcal{A}^{\pi_n}_{t,t_j^n}} \essinf_{\beta\in\mathcal{B}^{\pi_n}_{t,t_j^n}}G^{t,x;\alpha
,\beta}_{t,
t_j^n}
\bigl(\varphi \bigl(t_j^n,X_{t_j^n}^{t,x;\alpha,\beta}
\bigr) \bigr),
\end{equation}
where $C\in R$ is a constant independent of $\rho$, $n$ and $t_j^n$.
From the latter estimate, we deduce with the help of Lemmas
\ref{lemma_BSDE1} and \ref{lemma_BSDE2} that
%
%
%e4.44 #&#
%
\begin{eqnarray}
C\rho&\ge&\esssup_{\alpha\in\mathcal{A}^{\pi_n}_{t,t_j^n}} \essinf_{\beta\in\mathcal{B}^{\pi_n}_{t,t_j^n}}
\bigl(G^{t,x;\alpha
,\beta}_{t,
t_j^n} \bigl(\varphi \bigl(t_j^n,X_{t_j^n}^{t,x;\alpha,\beta}
\bigr) \bigr)-\varphi (t,x) \bigr)
\nonumber\\
&= & \esssup_{\alpha\in\mathcal{A}^{\pi_n}_{t,t_j^n}} \essinf_{\beta\in\mathcal{B}^{\pi_n}_{t,t_j^n}}Y_t^{1,\alpha
,\beta}
\\
&\ge& \esssup_{\alpha\in\mathcal{A}^{\pi_n}_{t,t_j^n}} \essinf_{\beta\in\mathcal{B}^{\pi_n}_{t,t_j^n}}Y_t^{2,\alpha
,\beta}
-C \bigl(t_j^n-t \bigr)^{3/2},\qquad P\mbox{-a.s.}\nonumber
\end{eqnarray}

Of course, as before, the quantities $Y_t^{1,\alpha,\beta},
Y_t^{2,\alpha,\beta}$ have to be understood as $Y_t^{1,u,v},Y_t^{2,u,v}$
for $(u,v)\in{\cal U}^{\pi_n}_{t,t_j^n}\times{\cal V}^{\pi_n}_{t,t_j^n}$
associated with $(\alpha,\beta)\in{\cal A}^{\pi_n}_{t,t_j^n}\times
{\cal
B}^{\pi_n}_{t,t_j^n}$ through Lemma \ref{controls-NAD-strategies}.
Moreover, they are defined by Lemmas \ref{lemma_BSDE1} and
\ref{lemma_BSDE2} for $t_j=t_j^n$, that is, they depend on the choice of
$t_j^n\in\pi_n$ and so, in particular, $n\ge n_\rho$.
Obviously, since ${\cal U}^{\pi_n}_{t,t_j^n}$ can be regarded as a subset
of ${\cal A}^{\pi_n}_{t,t_j^n}$ by identifying $u\in
{\cal U}^{\pi_n}_{t,t_j^n}$ with the NAD strategy $\alpha^u(v):=u,
v\in
{\cal V}^{\pi_n}_{t,t_j^n}$,
%
%
%e4.45 #&#
%
\begin{eqnarray}
\label{e1} %
C\rho+C \bigl(t_j^n-t
\bigr)^{3/2} &\ge& \esssup_{\alpha\in\mathcal{A}^{\pi
_n}_{t,t_j^n}} \essinf_{\beta\in\mathcal{B}^{\pi_n}_{t,t_j^n}}Y_t^{2,\alpha
,\beta}
\nonumber\\
&\ge& \esssup_{u\in\mathcal{U}^{\pi_n}_{t,t_j^n}} \essinf_{\beta\in\mathcal{B}^{\pi_n}_{t,t_j^n}}Y_t^{2,u,\beta
(u)}
\nonumber
\\[-8pt]
\\[-8pt]
\nonumber
&\ge& \esssup_{u\in\mathcal{U}^{\pi_n}_{t,t_j^n}} \essinf_{v\in\mathcal{V}^{\pi_n}_{t,t_j^n}}Y_t^{2,u,v}
\\
& = & Y^0_t,\qquad P\mbox{-a.s.}, n\ge n_\rho,\nonumber
\end{eqnarray}
where the latter equality was stated in Lemma \ref{lemma_BSDE3}.
Remark that here, of course, $Y^0$ is defined by Lemma \ref{lemma_BSDE3}
for $t_j^n$. Since
\[
Y^0_s=\int_s^{t_j^n}F_0
\bigl(r,x,Y_r^0,0 \bigr)\,dr, \qquad s\in
\bigl[t,t_j^n \bigr]
\]
and $F_0(r,x,y,0)$ is bounded, continuous, and Lipschitz in
$y$, uniformly with respect to $r$, it follows that
$|Y^0_s|\le C(t_j^n-t), s\in[t,t_j^n]$, and
%
%
%e4.46 #&#
%
\begin{eqnarray}
\label{e2} %
\frac{1}{t_j^n-t}{Y^0_t}&= &
\frac{1}{t_j^n-t} \int_t^{t_j^n}F_0
\bigl(r,x,Y_r^0,0 \bigr)\,dr
\nonumber\\
&\ge&\frac{1}{t_j^n-t}\int_t^{t_j^n}
\bigl(F_0(r,x,0,0)-L\bigl|Y_r^0\bigr| \bigr)\,dr
\\
&\ge&\frac{1}{t_j^n-t}\int_t^{t_j^n}F_0(r,x,0,0)
\,dr-C \bigl(t_j^n-t \bigr).\nonumber
\end{eqnarray}

Let $\rho\le(T-t)^{3/2}$. Since the mesh $|\pi_n|$ of the
partition $\pi_n$ converges to zero as $n\rightarrow+\infty$, we can find
for $n\ge n_\rho$ large enough some $t_j^n\in\pi_n, t_j^n>t$, such that
$(t_j^n-t)^{3/2}/2\le\rho\le(t_j^n-t)^{3/2}$. Consequently, for
$n\ge n_\rho$ large enough we can conclude from (\ref{e1}) and (\ref{e2})
that
\[
C \bigl(t_j^n-t \bigr)^{1/2}\ge
\frac{1}{t_j^n-t}{Y^0_t} \ge\frac{1}{t_j^n-t}\int
_t^{t_j^n}F_0(r,x,0,0)\,dr-C
\bigl(t_j^n-t \bigr).
\]
Thus, taking the limit as $\rho\rightarrow0$ (and, hence,
$n\rightarrow+\infty$ and $t_j^n-t\rightarrow0$), we obtain $F_0(t,x,0,0)
\le0$. But recalling the definition of $F_0$ from Lemma \ref{lemma_BSDE3},
we see that
%
%
%e4.47 #&#
%
\begin{eqnarray}
0&\ge& F_0(t,x,0,0)=\sup_{\mu\in{\cal P}(U)}\inf
_{\nu\in{\cal P}(V)} F(t,x,y,z,\mu,\nu)
\nonumber\\
& =& \sup_{\mu\in{\cal P}(U)}\inf_{\nu\in{\cal P}(V)}\nonumber\\
&&{}\times\int
_{U\times V} \biggl( \frac{\partial}{\partial t}\varphi(t,x) +
\frac{1}{2}\operatorname{tr} \bigl(\sigma \sigma^{T} (t,x, u, v)D^2
\varphi \bigr)
\nonumber
\\[-8pt]
\\[-8pt]
\nonumber
&&\hspace*{44pt}{}+ D\varphi.b(t, x, u, v)\\
&&\hspace*{44pt}{}+ f \bigl(t, x,\varphi(t,x),D\varphi
(t,x)\cdot\sigma (t,x,u, v),u, v \bigr) \biggr)\mu\otimes\nu(du\,dv)
\nonumber\\
& = &\frac{\partial}{\partial t}\varphi(t,x)+H \bigl(t,x, \bigl(\varphi, D
\varphi,D^2\varphi \bigr) (t,x) \bigr).\nonumber
\end{eqnarray}
Therefore, $W$ is a viscosity supersolution of PDE (\ref{HJBI}).

\textit{Step} 2. With an argument symmetric to that developed
in step 1 we show that $U$ is a viscosity subsolution of PDE (\ref{HJBI}).
Since both $W$ and $U$ are bounded continuous solutions, $W$ is a viscosity
supersolution and $U$ is a viscosity subsolution of (\ref{HJBI}), it follows
from the comparison principle (Proposition \ref{comparison principle})
that $W\ge U$ on $[0,T]\times R^d$. On the other hand, $(W,U)$ is the pointwise
limit over the sequence $(W^{\pi_n},U^{\pi_n}), n\ge1$, where the lower
value function $W^{\pi_n}$ along the partition $\pi_n$
is less than or equal to the upper one $U^{\pi_n}$, for all $n\ge1$.
Consequently, $W$ and $U$ coincide, and both are viscosity solutions of
PDE (\ref{HJBI}). Again from the comparison principle it follows that
this viscosity solution $W=U=V$ is the unique one inside the class of
continuous unions with at most polynomial growth.
\end{pf}

The above lemma allows now to prove Theorem \ref{main result}.

\begin{pf} From our above discussion, we have seen that for any arbitrary
sequence of partitions $\pi_n, n\ge1$, with $|\pi_n|\rightarrow0$, as
$n\rightarrow+\infty$, there is a subsequence which, abusing
notation, we
have also denoted by $\pi_n, n\ge1$, such that $W^{\pi_n}$ as well as
$U^{\pi_n}$ converge uniformly on compacts to the unique viscosity
solution $V$ of PDE (\ref{HJBI}) (uniqueness in the class of continuous
functions with polynomial growth); see Lemma \ref{viscosity solution}.
Consequently, the limit $V$ does not depend on the special choice of the
sequence of partitions $\pi_n, n\ge1$. Consequently, $W^{\pi_n}$ as
well as $U^{\pi_n}$ converge uniformly on compacts to the unique viscosity
solution $V$, for all sequence of partitions $\pi_n, n\ge1$ with mesh
$|\pi_n|\rightarrow0$, as $n\rightarrow+\infty$. The proof is complete.
\end{pf}

\section*{Acknowledgements}
The authors would like to thank the anonymous Associate Editor and the
anonymous referee for their valuable comments and suggestions from
which the manuscript greatly has benefited.

% imsref loaded by akundreckaite, 2013-10-21 09:44:43
% imsref loaded by akundreckaite, 2013-10-21 10:32:45


%suskaldyti doi

\printaddresses

\end{document}